\documentclass[AMA,STIX1COL]{WileyNJD-v2}
\usepackage[T1]{fontenc}
\usepackage[utf8]{inputenc}
\usepackage{units}
\usepackage{amsmath}
\usepackage{amsthm}
\usepackage{amssymb}
\usepackage{stmaryrd}
\usepackage{graphicx}

\makeatletter

\articletype{Research Article}
\received{February 2020}
\revised{}
\accepted{} 
\jnlcitation{\cname{%
\author{A. D. Kercher},
\author{A. Corrigan}, and
\author{D. A. Kessler}}
(\cyear{2020}),
\ctitle{The Moving Discontinuous Galerkin Finite Element Method with Interface Condition Enforcement for Compressible Viscous Flows},
\cjournal{IJNMF},
\cvol{2020;00:1--6}.} 

\providecommand{\tabularnewline}{\\}

\numberwithin{equation}{section}
\numberwithin{figure}{section}

\usepackage{tikz}
\usetikzlibrary{shapes,arrows}

\@ifundefined{showcaptionsetup}{}{%
 \PassOptionsToPackage{caption=false}{subfig}}
\usepackage{subfig}
\makeatother

\begin{document}
\title{The Moving Discontinuous Galerkin Finite Element Method with Interface
Condition Enforcement for Compressible Viscous Flows}

\author[lcp]{Andrew D. Kercher}
\author[lcp]{Andrew Corrigan}
\author[lcp]{David A. Kessler}
\authormark{Kercher \textsc{et al}}
\address[lcp]{Laboratories for Computational Physics and Fluid Dynamics,  U.S. Naval Research Laboratory, 4555 Overlook Ave SW, Washington, DC 20375}
\corres{*Andrew D. Kercher, \email{andrew.kercher@nrl.navy.mil}} 
\global\long\def\middlebar{\,\middle|\,}%
\global\long\def\cof{\operatorname{cof}}%
\global\long\def\det{\operatorname{det}}%
\global\long\def\adj{\operatorname{adj}}%
\global\long\def\average#1{\left\{  \!\!\left\{  #1\right\}  \!\!\right\}  }%
\global\long\def\expnumber#1#2{{#1}\mathrm{e}{#2}}%

\abstract{The moving discontinuous Galerkin finite element method with interface
condition enforcement (MDG-ICE) is applied to the case of viscous
flows. This method uses a weak formulation that separately enforces
the conservation law, constitutive law, and the corresponding interface
conditions in order to provide the means to detect interfaces or under-resolved
flow features. To satisfy the resulting overdetermined weak formulation,
the discrete domain geometry is introduced as a variable, so that
the method implicitly fits a priori unknown interfaces and moves the
grid to resolve sharp, but smooth, gradients, achieving a form of
anisotropic curvilinear $r$-adaptivity. This approach avoids introducing
low-order errors that arise using shock capturing, artificial dissipation,
or limiting. The utility of this approach is demonstrated with its
application to a series of test problems culminating with the compressible
Navier-Stokes solution to a Mach 5 viscous bow shock for a Reynolds
number of $10^{5}$ in two-dimensional space. Time accurate solutions
of unsteady problems are obtained via a space-time formulation, in
which the unsteady problem is formulated as a higher dimensional steady
space-time problem. The method is shown to accurately resolve and
transport viscous structures without relying on numerical dissipation
for stabilization.}

\maketitle
\footnotetext{\textbf{Abbreviations:} MDG-ICE, Moving Discontinuous Galerkin Finite Element Method with Interface Condition Enforcement \\ \\ \hspace*{140pt}Distribution A. Approved for public release: distribution unlimited.} 

\section{Introduction\label{sec:Introduction}}

The discontinuous Galerkin (DG) method~\cite{Bas97,Bas97_2,Coc98,Coc00,Arn02,Har02,Fid05,Hes07,Per08,Har13},
has become a popular method for simulating flow fields corresponding
to a wide range of physical phenomena, from low speed incompressible
flows~\cite{Liu00,Bas07,Rhe12} to chemically reacting compressible
Navier-Stokes flows~\cite{Lv14,Joh19_SCITECH,Joh19}, due to its
ability to achieve high-order accuracy on unstructured grids and its
natural support for local polynomial, $p$, adaptivity. However, the
solutions are known to contain oscillations in under-resolved regions
of the flow, e.g., shocks, material interfaces, and boundary layers,
where stabilization, usually in the form of shock capturing, or limiting,
is required. These ad hoc methods often lead to inconsistent discretizations
that are no longer capable of achieving high-order accuracy~\cite{Wan13}.
This lack of robustness and accuracy when treating discontinuities
and regions with sharp, but smooth, gradients is one of the main obstacles
preventing widespread adoption of high-order methods for the simulation
of complex high-speed turbulent flow fields.

The Moving Discontinuous Galerkin Method with Interface Condition
Enforcement (MDG-ICE) was introduced by the present authors~\cite{Cor17,Cor18}
as a high-order method for computing solutions to inviscid flow problems,
even in the presence of discontinuous interfaces. The method accurately
and stably computes flows with a priori unknown interfaces without
relying on shock capturing. In order to detect a priori unknown interfaces,
MDG-ICE uses a weak formulation that enforces the conservation law
and its interface condition separately, while treating the discrete
domain geometry as a variable. Thus, in contrast to a standard DG
method, MDG-ICE has both the means to detect via interface condition
enforcement and satisfy via grid movement the conservation law and
its associated interface condition. Using this approach, not only
does the grid move to fit interfaces, but also to better resolve smooth
regions of the flow.

In this work, we apply MDG-ICE to the case of viscous flows, and therefore
extend the capability of MDG-ICE to move the grid to resolve otherwise
under-resolved flow features such as boundary layers and viscous shocks.
In addition to enforcing the conservation law and interface (Rankine-Hugoniot)
condition, for the case of viscous flows, we separately enforce a
constitutive law and the corresponding interface condition, which
constrains the continuity of the state variable across the interface.
Thus, in contrast to a standard DG method, MDG-ICE implicitly achieves
a form of anisotropic curvilinear $r$-adaptivity via satisfaction
of the weak formulation.

We study the utility of this approach by solving problems involving
linear advection-diffusion, unsteady Burgers flow, and steady compressible
Navier-Stokes flow. The ability of the method to move the grid in
order to resolve boundary layer profiles is studied in the context
of one-dimensional linear advection-diffusion, where convergence under
polynomial refinement is also considered. The problem of space-time
viscous shock formation for Burgers flow is considered in order to
assess the ability of the method to accurately resolve and transport
viscous shocks without relying on shock capturing or limiting. Lastly,
a Mach 5 compressible Navier-Stokes flow over a cylindrical blunt
body in two dimensions is studied for a series of increasing Reynolds
numbers to assess the ability of MDG-ICE to simultaneously resolve
multiple viscous structures, i.e., a viscous shock and boundary layer,
via anisotropic curvilinear $r$-adaptivity.

\subsection{Background\label{sec:Background}}

In prior work, MDG-ICE was shown to be a consistent discretization
for discontinuous flows and is therefore capable of using high-order
approximations to achieve extremely accurate solutions for problems
containing discontinuous interfaces on relatively coarse grids~\cite{Cor18}.
The previously presented test cases demonstrated that MDG-ICE can
be used to compute both steady and unsteady flows with a priori unknown
interface topology and point singularities using higher-order elements
in arbitrary-dimensional spaces. For example, MDG-ICE was applied
to fit intersecting oblique planar shocks in three dimensions. The
ability to fit steady shocks extends to unsteady flows using a space-time
formulation, cf. Lowrie et al.~\cite{Low95,Low98}, that was applied
to compute the solution to a space-time inviscid Burgers shock formation
problem, where a continuous temporal initial condition steepens to
form a shock, while later work presented proof-of-concept results
for unsteady flow in three- and four-dimensional space-time~\cite{Cor19_SCITECH}.
More recently, MDG-ICE was applied to shocked compressible flow problems
of increased complexity, including transonic flow over a smooth bump
over which an attached curved shock forms for which optimal-order
convergence was verified~\cite{Cor18_AVIATION,Cor19_AVIATION}.

Earlier attempts at aligning the grid with discontinuous interfaces
present in the flow field have resulted in mixed success, cf. Moretti~\cite{Mor02}
and Salas~\cite{Sal09,Sal11}. These earlier $explicit$ shock fitting,
or tracking, approaches were capable of attaining high-order accuracy
in the presence of shocks, but the general applicability of such methods
is limited. A specialized discretization and implementation strategy
is required at discontinuous interfaces making it difficult to handle
discontinuities whose topologies are unknown a priori or whose topologies
evolve in time. In contrast, MDG-ICE is an $implicit$ shock fitting
method, automatically fitting a priori unknown interfaces, their interactions,
and their evolving trajectory in space-time.

Another promising form of implicit shock tracking, or fitting, is
the optimization-based, $r$-adaptive, approach proposed independently
by Zahr and Persson~\cite{Zah17,Zah18_AIAA,Zah18}, which has been
used to compute very accurate solutions to discontinuous inviscid
flows on coarse grids without the use of artificial stabilization.
This approach retains a standard discontinuous Galerkin method as
the state equation, while employing an objective function to detect
and fit interfaces present in the flow. Recently, Zahr et. al.~\cite{Zah19,Zah20_SCITECH}
have extended their implicit shock tracking framework with the addition
of a new objective function based on an enriched test space and an
improved SQP solver. Furthermore, the regularization introduced by
the current authors~\cite{Cor17,Cor18,Cor18_AVIATION,Cor19_AVIATION,Cor19_SCITECH}
has been modified to include a factor proportional to the inverse
volume of an element that accounts for variations in the element sizes
of a given grid. This may be beneficial for obtaining solutions on
highly nonuniform grids.

In the case of viscous flows, regions with sharp, but smooth, gradients
present a unique set of challenges. The resolution required to achieve
high-order convergence, or at a minimum, achieve stability, is such
that computations on uniformly refined grids are prohibitively expensive.
Therefore, the local resolution must be selectively increased in certain
regions of the flow. Identifying these regions is not always obvious
and striking a balance between computational feasibility and accuracy
is an equally challenging task. Traditionally, overcoming these challenges
was viewed as an a priori grid design problem with solutions including
anisotropic grid generation~\cite{Tam00,Pai01,Geo02,Bot04,Li05,Jon06,Dob08,Com10,Los10,Los11}
and boundary layer grid generation~\cite{Loh93,Loh00,Pir94,Mar96,Gar00,Bot02,Ito02,Ito06,Aub09}. 

A complementary approach to problem-specific grid generation is a
posteriori grid adaptation~\cite{Fid11,Yan12,Har13,Car20}. This
is an iterative process in which regions of interest are locally refined
with the goal of reducing the discretization error. Anisotropic grid
adaptation, which combines a posteriori grid adaptation with anisotropic
grid generation, has been shown to successfully enhance accuracy for
a range of aerodynamic applications as reviewed by Alauzet and Loseille~\cite{Ala16}.
MDG-ICE seeks to achieve similar anisotropic grid adaptation as an
intrinsic part of the solver, such that the region of anisotropic
refinement evolves with the flow field solution, thereby avoiding
grid coarsening as the viscous layer is more accurately resolved.

In the case of least-squares (LS) methods, the residual is a natural
indicator of the discretization error~\cite{Jia87,Car97}. In particular,
for the discontinuous Petrov-Galerkin (DPG) method introduced by Demkowicz
and Gopalakrishnan~\cite{Dem10,Dem11,Dem12,Gop13,Gop14,Car16,Dem18},
in which the ultra-weak formulation corresponds to the best approximation
in the polynomial space, a posteriori grid adaptation, for both the
grid resolution $h$ and the polynomial degree $p$, is driven by
the built-in error representation function in the form of the Riesz
representation of the residual~\cite{Dem15_20}. In addition to such
a posteriori $hp$-adaptivity strategies, MDG-ICE achieves a form
of in situ $r$-adaptivity~\cite{Mil81,Mil81_2,Gel81,Ban93,Boc96,Roe02,San16}
where the resolution of the flow is continuously improved through
repositioning of the grid points. For a review of $r$-adaptivity
the reader is referred to the work of Budd et al.~\cite{Bud09} and
the survey of Huang and Russell~\cite{Hua10}. 

\section{Moving Discontinuous Galerkin Method with Interface Condition Enforcement
for Compressible Viscous Flows\label{sec:Formulation}}

In this section we develop the formulation of MDG-ICE for compressible
viscous flows. We assume that $\Omega\subset\mathbb{R}^{d}$ is a
given domain, which may be either a spatial domain $\Omega\subset\mathbb{R}^{d=d_{x}}$
or a space-time domain $\Omega\subset\mathbb{R}^{d=d_{x}+1}$. In
many cases, the space-time domain is defined in terms of a fixed spatial
domain $\Omega_{x}\subset\mathbb{R}^{d_{x}}$ and time interval $T\subset\left\{ t\in\mathbb{R}:t>0\right\} $
by $\Omega=\Omega_{x}\times T$. In the remainder of this work, we
assume that $\Omega$ is partitioned by $\mathcal{T}$, consisting
of disjoint sub-domains or cells $\kappa$, so that $\overline{\Omega}=\cup_{\kappa\in\mathcal{T}}\overline{\kappa}$,
with interfaces $\epsilon$, composing a set $\mathcal{E}$ so that
$\cup_{\epsilon\in\mathcal{E}}\epsilon=\cup_{\kappa\in\mathcal{T}}\partial\kappa$.
Furthermore, we assume that each interface $\epsilon$ is oriented
so that a unit normal $n:\epsilon\rightarrow\mathbb{R}^{d}$ is defined.
In order to account for space-time problems, we also consider the
spatial normal $n_{x}:\epsilon\rightarrow\mathbb{R}^{d_{x}}$, which
is defined such that $\left(n_{x,1},\ldots n_{x,d_{x}}\right)=\left(n_{1},\ldots n_{d_{x}}\right)$.

\subsection{Governing equations\label{subsec:Governing-equations-Primal-Formulation}}

Consider a nonlinear conservation law governing the behavior of smooth,
$\mathbb{R}^{m}$-valued, functions $y$,

\begin{align}
\nabla\cdot\mathcal{F}\left(y,\nabla_{x}y\right)=0 & \textup{ in }\Omega,\label{eq:conservation-strong-primal}
\end{align}
in terms of a given flux function, $\mathcal{F}:\mathbb{R}^{m}\times\mathbb{R}^{m\times d_{x}}\rightarrow\mathbb{R}^{m\times d}$
that depends on the flow state variable $y$ and its $d_{x}$-dimensional
spatial gradient,
\begin{equation}
\nabla_{x}y=\left(\frac{\partial y}{\partial x_{1}},\ldots,\frac{\partial y}{\partial x_{d_{x}}}\right).\label{eq:spatial-gradient}
\end{equation}
The flux function is assumed to be defined in terms of a spatial flux
function $\mathcal{F}^{x}:\mathbb{R}^{m}\times\mathbb{R}^{m\times d_{x}}\rightarrow\mathbb{R}^{m\times d_{x}}$
that itself is defined in terms of a convective flux, depending on
the state variable only, and the viscous, or diffusive flux, which
also depends on the spatial gradient of the state variable,
\begin{equation}
\mathcal{F}^{x}\left(y,\nabla_{x}y\right)=\mathcal{F}^{c}\left(y\right)-\mathcal{F}^{v}\left(y,\nabla_{x}y\right).\label{eq:navier-stokes-spatial-flux}
\end{equation}
In the case of a spatial domain, $d=d_{x}$, the flux function $\mathcal{F}$
coincides with the spatial flux, 
\begin{equation}
\mathcal{F}\left(y,\nabla_{x}y\right)=\mathcal{F}^{x}\left(y,\nabla_{x}y\right),\label{eq:spatial-flux}
\end{equation}
so that the divergence operator in Equation~(\ref{eq:conservation-strong-primal})
is defined as the spatial divergence operator
\begin{equation}
\nabla\cdot\mathcal{F}\left(y,\nabla_{x}y\right)=\nabla_{x}\cdot\mathcal{F}^{x}\left(y,\nabla_{x}y\right)=\frac{\partial}{\partial x_{1}}\mathcal{F}_{1}^{x}\left(y,\nabla_{x}y\right)+\ldots+\frac{\partial}{\partial x_{d_{x}}}\mathcal{F}_{d_{x}}^{x}\left(y,\nabla_{x}y\right).\label{eq:spatial-divergence}
\end{equation}
Otherwise, in the case of a space-time domain, $d=d_{x}+1$, the space-time
flux incorporates the state variable as the temporal flux component,
\begin{equation}
\mathcal{F}\left(y,\nabla_{x}y\right)=\left(\mathcal{F}_{1}^{x}\left(y,\nabla_{x}y\right),\ldots,\mathcal{F}_{d_{x}}^{x}\left(y,\nabla_{x}y\right),y\right),\label{eq:spacetime-flux}
\end{equation}
so that the divergence operator in~(\ref{eq:conservation-strong-primal})
is defined as the space-time divergence operator
\begin{equation}
\nabla\cdot\mathcal{F}\left(y,\nabla_{x}y\right)=\nabla_{x}\cdot\mathcal{F}^{x}\left(y,\nabla_{x}y\right)+\frac{\partial}{\partial t}y.\label{eq:spacetime-divergence}
\end{equation}

In this work, we consider conservation laws corresponding to linear
advection-diffusion, space-time viscous Burgers, and compressible
Navier-Stokes flow as detailed in the following sections.

\subsubsection{Linear advection-diffusion\label{subsec:Linear-advection-diffusion}}

Linear advection-diffusion involves a single-component flow state
variable $y:\Omega\rightarrow\mathbb{R}^{1}$ with a linear diffusive
flux,
\begin{equation}
\mathcal{F}^{v}\left(y,\nabla_{x}y\right)=\epsilon\nabla_{x}y,\label{eq:linear-viscous-flux}
\end{equation}
that is independent of the state $y$, where the coefficient $\epsilon$
corresponds to mass diffusivity. The convective flux is given as 
\begin{equation}
\mathcal{F}^{c}\left(y\right)=\left(v_{1}y,\ldots,v_{d_{x}}y\right),\label{eq:linear-convective-flux}
\end{equation}
where $\left(v_{1},\ldots,v_{d_{x}}\right)\in\mathbb{R}^{d_{x}}$
is a prescribed spatial velocity that in the present setting is assumed
to be spatially uniform. The corresponding spatial flux is given by
\begin{equation}
\mathcal{F}^{x}\left(y,\nabla_{x}y\right)=\left(\left(v_{1}y,\ldots,v_{d_{x}}y\right)-\epsilon\nabla_{x}y\right),\label{eq:spatial-linear-advection-flux}
\end{equation}

\subsubsection{One-dimensional Burgers flow\label{subsec:One-dimensional-Burgers-flow}}

As in the case of linear advection-diffusion, one-dimensional Burgers
flow involves a single-component flow state variable $y:\Omega\rightarrow\mathbb{R}^{1}$
with a linear viscous flux,
\begin{equation}
\mathcal{F}^{v}\left(y,\nabla_{x}y\right)=\epsilon\nabla_{x}y,\label{eq:linear-viscous-flux-Burgers}
\end{equation}
which is independent of the state $y$, where the coefficient, $\epsilon$,
corresponds to viscosity. The convective flux is given as 
\begin{equation}
\mathcal{F}^{c}\left(y\right)=\left(\frac{1}{2}y^{2}\right),\label{eq:convective-burgers-flux}
\end{equation}
so that the one-dimensional spatial flux is given by
\begin{equation}
\mathcal{F}^{x}\left(y,\nabla_{x}y\right)=\left(\frac{1}{2}y^{2}-\epsilon\nabla_{x}y\right),\label{eq:spatial-burgers-flux}
\end{equation}

\subsubsection{Compressible Navier-Stokes flow\label{subsec:Compressible-Navier-Stokes-flow}}

For compressible Navier-Stokes flow, the state variable $y:\Omega\rightarrow\mathbb{R}^{m}$,
where $m=d_{x}+2$, is given by

\begin{equation}
y=\left(\rho,\rho v_{1},\ldots,\rho v_{d_{x}},\rho E\right).\label{eq:navier-stokes-state}
\end{equation}
The $i$-th spatial component of the convective flux, $\mathcal{F}^{c}:\mathbb{R}^{m}\rightarrow\mathbb{R}^{m\times d_{x}}$,
is
\begin{equation}
\mathcal{F}_{i}^{c}\left(y\right)=\left(\rho v_{i},\rho v_{i}v_{1}+p\delta_{i1},\ldots,\rho v_{i}v_{d_{x}}+p\delta_{id_{x}},\rho Hv_{i}\right),\label{eq:ns-convective-flux-spatial-component}
\end{equation}
where $\delta_{ij}$ is the Kronecker delta, $\rho:\Omega\rightarrow\mathbb{R}_{+}$
is density, $\left(v_{1},\ldots,v_{d_{x}}\right):\mathbb{R}^{m}\rightarrow\mathbb{R}^{d_{x}}$
is velocity, $\rho E:\Omega\rightarrow\mathbb{R}_{+}$ is stagnation
energy per unit volume, and 
\begin{equation}
H=\left(\rho E+p\right)/\rho\label{eq:stagnation-enthalpy}
\end{equation}
 is stagnation enthalpy, where $H:\mathbb{R}^{m}\rightarrow\mathbb{R}_{+}$.
Assuming the fluid is a perfect gas, the pressure $p:\mathbb{R}^{m}\rightarrow\mathbb{R}_{+}$
is defined as 
\begin{equation}
p=\left(\gamma-1\right)\left(\rho E-\frac{1}{2}\sum_{i=1}^{d_{x}}\rho v_{i}v_{i}\right),\label{eq:pressure}
\end{equation}
where the ratio of specific heats for air is given as $\gamma=1.4$.
The $i$-th spatial component of the viscous flux is given by
\begin{equation}
\mathcal{F}_{i}^{\nu}\left(y,\nabla_{x}y\right)=\left(0,\tau_{1i},\ldots,\tau_{d_{x}i},\sum_{j=1}^{d_{x}}\tau_{ij}v_{j}-q_{i}\right),\label{eq:navier-stokes-viscous-flux-spatial-component}
\end{equation}

\noindent where $q:\mathbb{R}^{m}\times\mathbb{R}^{m\times d_{x}}\rightarrow\mathbb{R}^{d_{x}}$
is the thermal heat flux, $\tau:\mathbb{R}^{m}\times\mathbb{R}^{m\times d_{x}}\rightarrow\mathbb{R}^{d_{x}\times d_{x}}$
is the viscous stress tensor. The $i$-th spatial component of the
thermal heat flux is given by

\noindent 
\begin{equation}
q_{i}=-k\frac{\partial T}{\partial x_{i}},\label{eq:navier-stokes-thermal-heat-flux-component}
\end{equation}
where $T:\mathbb{R}^{m}\rightarrow\mathbb{R}_{+}$ is the temperature
and $k$ is thermal conductivity. The temperature $T$ is defined
as

\begin{equation}
T=\frac{p}{R\rho},\label{eq:temperature}
\end{equation}

\noindent where $R=287$ is the mixed specific gas constant for air.
The $i$-th spatial component of the viscous stress tensor is given
by

\noindent 
\begin{equation}
\tau_{i}=\mu\left(\frac{\partial v_{1}}{\partial x_{i}}+\frac{\partial v_{i}}{\partial x_{1}}-\delta_{i1}\frac{2}{3}\sum_{j=1}^{d_{x}}\frac{\partial v_{j}}{\partial x_{j}},\ldots,\frac{\partial v_{d_{x}}}{\partial x_{i}}+\frac{\partial v_{i}}{\partial x_{d_{x}}}-\delta_{id_{x}}\frac{2}{3}\sum_{j=1}^{d_{x}}\frac{\partial v_{j}}{\partial x_{j}}\right),\label{eq:navier-stokes-viscous-stress-tensor-component}
\end{equation}

\noindent where $\mu$ is the dynamic viscosity coefficient.

\subsection{Interface conditions for viscous flow\label{sec:visous-strong-formulation}}

The viscous conservation laws described in the previous sections require
a constraint on the continuity of the state variable, $y$, across
an interface, in addition to the interface condition considered in
our previous work~\cite{Cor18}, which enforced the continuity of
the normal flux across an interface. In order to deduce the interface
conditions governing viscous flow, we revisit the derivation of the
DG formulation for viscous flow, cf.~Arnold et al.~\cite{Arn02}.
This discussion follows Section~6.3 of Hartmann and Leicht~\cite{Har13}
and restricts the presentation to a viscous flux. Upon deducing the
governing interface conditions, we will reintroduce the convective
flux in Section~\ref{subsec:strong-form}. Here, we consider a spatial
conservation law,
\begin{align}
-\nabla_{x}\cdot\left(\mathcal{F}^{v}\left(y,\nabla_{x}y\right)\right)=0 & \textup{ in }\kappa\qquad\forall\kappa\in\mathcal{T},\label{eq:viscous-operator-strong}
\end{align}
defined in terms of a given viscous flux function $\mathcal{F}^{v}:\mathbb{R}^{m}\times\mathbb{R}^{m\times d_{x}}\rightarrow\mathbb{R}^{m\times d_{x}}$,
for piecewise smooth functions $y$ and their spatial gradients $\nabla_{x}y$.
We introduce an $\mathbb{R}^{m\times d_{x}}$-valued auxiliary variable
$\sigma$ and rewrite~(\ref{eq:viscous-operator-strong}) as a first-order
system of equations
\begin{align}
-\nabla_{x}\cdot\sigma=0 & \textup{ in }\kappa\qquad\forall\kappa\in\mathcal{T},\label{eq:viscous-operator-strong-fos}\\
\sigma-G\left(y\right)\nabla_{x}y=0 & \textup{ in }\kappa\qquad\forall\kappa\in\mathcal{T}.\label{eq:viscous-constitutive-strong-fos}
\end{align}
We assume here that $\mathcal{F}^{v}$ is linear with respect to its
gradient argument so that
\begin{equation}
G\left(y\right)\nabla_{x}y=\mathcal{F}^{v}\left(y,\nabla_{x}y\right)=\mathcal{F}_{\nabla_{x}y}^{v}\left(y,\nabla_{x}y\right)\nabla_{x}y\label{eq:homogeneity-tensor}
\end{equation}
where $G\left(y\right)\in\mathbb{R}^{m\times d_{x}\times m\times d_{x}}$
is a tensor of rank 4 that is referred to as the \emph{homogeneity
tensor}~\cite{Har13}.

We integrate~(\ref{eq:viscous-operator-strong-fos}) and~(\ref{eq:viscous-constitutive-strong-fos})
against separate test functions and upon an application of integration
by parts arrive at the following weak formulation : find $\left(y,\sigma\right)\in Y\times\Sigma$
such that

\begin{align}
0= & +\sum_{\kappa\in\mathcal{T}}\left(\sigma,\nabla_{x}v\right)_{\kappa}\nonumber \\
 & -\sum_{\kappa\in\mathcal{T}}\left(\sigma\cdot n_{x},v\right)_{\partial\kappa}\nonumber \\
 & +\sum_{\kappa\in\mathcal{T}}\left(\sigma,\tau\right)_{\kappa}\nonumber \\
 & +\sum_{\kappa\in\mathcal{T}}\left(y,\nabla_{x}\cdot\left(G(y)^{\top}\tau\right)\right)_{\kappa}\nonumber \\
 & -\sum_{\kappa\in\mathcal{T}}\left(y\otimes n_{x},G(y)^{\top}\tau\right)_{\partial\kappa}\qquad\forall\left(v,\tau\right)\in V_{y}\times V_{\sigma},\label{eq:weak-formulation-viscous-operator-fos-integrated-by-parts}
\end{align}
where the solution spaces $Y\times\Sigma$ and test spaces $V_{y}\times V_{\sigma}$
are broken Sobolev spaces. Since $y$ and $\sigma$ are multi-valued
across element interfaces, in a DG formulation, they are substituted
with single-valued functions of their traces,
\begin{align}
\hat{\sigma}= & \hat{\sigma}\left(\sigma^{+},\sigma^{-}\right),\label{eq:auxiliary-numerical-flux}\\
\hat{y}= & \hat{y}\left(y^{+},y^{-}\right),\label{eq:state-numerical-flux}
\end{align}
cf.~Table~3.1 of Arnold et al.~\cite{Arn02} for various definitions
of both $\hat{\sigma}$ and $\hat{y}$. After another application
of integration by parts and transposition of the homogeneity tensor,
we obtain: find $\left(y,\sigma\right)\in Y\times\Sigma$ such that
\begin{align}
0= & -\sum_{\kappa\in\mathcal{T}}\left(\nabla_{x}\cdot\sigma,v\right)_{\kappa}\nonumber \\
 & +\sum_{\kappa\in\mathcal{T}}\left(\left(\sigma-\hat{\sigma}\right)\cdot n_{x},v\right)_{\partial\kappa}\nonumber \\
 & +\sum_{\kappa\in\mathcal{T}}\left(\sigma-G(y)\nabla_{x}y,\tau\right)_{\kappa}\nonumber \\
 & -\sum_{\kappa\in\mathcal{T}}\left(G(y)\left(\left(\hat{y}-y\right)\otimes n_{x}\right),\tau\right)_{\partial\kappa}\qquad\forall\left(v,\tau\right)\in V_{y}\times V_{\sigma}.\label{eq:weak-formulation-viscous-operator-fos}
\end{align}
Finally, the auxiliary variable $\sigma$ is substituted with $G\left(y\right)\nabla_{x}y$,
the tensor-valued test function $\tau$ is substituted with $\nabla_{x}v$,
so that upon a final application of integration by parts we obtain
a DG primal formulation: find $y\in Y$ such that
\begin{align}
0= & \sum_{\kappa\in\mathcal{T}}\left(G(y)\nabla_{x}y,\nabla_{x}v\right)_{\kappa}\nonumber \\
 & -\sum_{\kappa\in\mathcal{T}}\left(\hat{\sigma}\cdot n_{x},v\right)_{\partial\kappa}\nonumber \\
 & +\sum_{\kappa\in\mathcal{T}}\left(G(y)\left(\left(\hat{y}-y\right)\otimes n_{x}\right),\nabla_{x}v\right)_{\partial\kappa}\qquad\forall v\in V_{y},\label{eq:weak-formulation-viscous-operator-fos-1}
\end{align}
cf.~Equation~(254) and Section~6.6 in the work of Hartmann and
Leicht~\cite{Har13}.

In contrast, we propose an MDG-ICE formulation that retains the auxiliary
variable and instead makes a different substitution: the test functions
$v$ and $\tau$ that appear in the surface integrals of~(\ref{eq:weak-formulation-viscous-operator-fos})
are substituted with separate test functions $w_{y}\in W_{y}$ and
$w_{\sigma}\in W_{\sigma}$ from the single-valued trace spaces of
$V_{y}$ and $V_{\sigma}$. Upon accumulating contributions from adjacent
elements in~(\ref{eq:weak-formulation-viscous-operator-fos}) to
each interface, we obtain: find $\left(y,\sigma\right)\in Y\times\Sigma$
\begin{align}
0= & -\sum_{\kappa\in\mathcal{T}}\left(\nabla_{x}\cdot\sigma,v\right)_{\kappa}\nonumber \\
 & +\sum_{\epsilon\in\mathcal{E}}\left(\left\llbracket n_{x}\cdot\sigma\right\rrbracket ,w_{y}\right)_{\epsilon}\nonumber \\
 & +\sum_{\kappa\in\mathcal{T}}\left(\sigma-G(y)\nabla_{x}y,\tau\right)_{\kappa}\nonumber \\
 & -\sum_{\epsilon\in\mathcal{E}}\left(\average{G\left(y\right)}\left\llbracket y\otimes n_{x}\right\rrbracket ,w_{\sigma}\right)_{\epsilon}\qquad\forall\left(v,\tau,w_{y},w_{\sigma}\right)\in V_{y}\times V_{\sigma}\times W_{y}\times W_{\sigma}.\label{eq:weak-formulation-viscous-operator-fos-interface}
\end{align}
We make use of the relationship
\begin{align}
 & \left(\sigma^{+}-\hat{\sigma}\right)\cdot n_{x}^{+}+\left(\sigma^{-}-\hat{\sigma}\right)\cdot n_{x}^{-}\nonumber \\
= & \left(\sigma^{+}-\hat{\sigma}\right)\cdot n_{x}^{+}-\left(\sigma^{-}-\hat{\sigma}\right)\cdot n_{x}^{+}\nonumber \\
= & \left(\sigma^{+}-\sigma^{-}\right)\cdot n_{x}^{+}\nonumber \\
= & \left\llbracket n_{x}\cdot\sigma\right\rrbracket ,\label{eq:weak-formulation-viscous-oprator-fos-jump-condition}
\end{align}
so that contributions from $\hat{\sigma}$ vanish, and
\begin{align}
 & G(y^{+})\left(\left(\hat{y}-y^{+}\right)\otimes n_{x}^{+}\right)+G(y^{-})\left(\left(\hat{y}-y^{-}\right)\otimes n_{x}^{-}\right)\nonumber \\
= & G(y^{+})\left(\frac{1}{2}\left(y^{-}-y^{+}\right)\otimes n_{x}^{+}\right)+G(y^{-})\left(\frac{1}{2}\left(y^{+}-y^{-}\right)\otimes n_{x}^{-}\right)\nonumber \\
= & \frac{1}{2}\left(G(y^{+})\left(\left(y^{-}-y^{+}\right)\otimes n_{x}^{+}\right)+G(y_{h}^{-})\left(\left(y^{+}-y^{-}\right)\otimes n_{x}^{-}\right)\right)\nonumber \\
= & \frac{1}{2}\left(G(y^{+})+G(y^{-})\right)\left(y^{-}-y^{+}\right)\otimes n_{x}^{+}\nonumber \\
= & -\average{G\left(y\right)}\left\llbracket y\otimes n_{x}\right\rrbracket ,\label{eq:weak-formulation-viscous-oprator-fos-constitutive-interface-condition}
\end{align}
on interior interfaces~(\ref{eq:interior-interfaces}), where we
define $\hat{y}=\average y$ , a common choice among the various DG
discretizations~\cite{Arn02,Har13}.

From~(\ref{eq:weak-formulation-viscous-operator-fos-interface}),
we deduce the strong form of the viscous interface conditions to be
\begin{align}
\left\llbracket n_{x}\cdot\sigma\right\rrbracket =0 & \textup{ on }\epsilon\qquad\forall\epsilon\in\mathcal{E},\label{eq:viscous-interface-condition-strong}\\
\average{G\left(y\right)}\left\llbracket y\otimes n_{x}\right\rrbracket =0 & \textup{ on }\epsilon\qquad\forall\epsilon\in\mathcal{E}.\label{eq:viscous-interface-condition-state-strong}
\end{align}
The first interface condition, Equation~(\ref{eq:viscous-interface-condition-strong})
is the jump or Rankine-Hugoniot condition~\cite{Maj12} that ensures
continuity of the normal flux at the interface and will balance with
the jump in the normal convective flux in Equation~(\ref{eq:interface-condition-strong-viscous}).
The second interface condition~(\ref{eq:viscous-interface-condition-state-strong})
corresponds to the constitutive law~(\ref{eq:constitutive-strong-viscous})
and enforces a constraint on the continuity of the state variable
at the interface.

\subsection{Formulation in physical space with fixed geometry\label{subsec:formulation-physical-space}}

Having deduced the interface conditions that arise in the case of
viscous flow, we reintroduce the convective flux and write the second
order system~(\ref{eq:conservation-strong-primal}) as a system of
first-order equations, incorporating the additional interface conditions~(\ref{eq:viscous-interface-condition-strong})
and~(\ref{eq:viscous-interface-condition-state-strong}).

\subsubsection{Strong formulation\label{subsec:strong-form}}

Consider a nonlinear conservation law, generalized constitutive law,
and their corresponding interface conditions,

\begin{align}
\nabla\cdot\mathcal{F}\left(y,\sigma\right)=0 & \textup{ in }\kappa\qquad\forall\kappa\in\mathcal{T},\label{eq:conservation-strong-viscous}\\
\sigma-G(y)\nabla_{x}y=0 & \textup{ in }\kappa\qquad\forall\kappa\in\mathcal{T},\label{eq:constitutive-strong-viscous}\\
\left\llbracket n\cdot\mathcal{F}\left(y,\sigma\right)\right\rrbracket =0 & \textup{ on }\epsilon\qquad\forall\epsilon\in\mathcal{E},\label{eq:interface-condition-strong-viscous}\\
\average{G\left(y\right)}\left\llbracket y\otimes n_{x}\right\rrbracket =0 & \textup{ on }\epsilon\qquad\forall\epsilon\in\mathcal{E},\label{eq:interface-condition-state-strong-viscous}
\end{align}
governing the flow state variable $y$ and auxiliary variable $\sigma$.
The interface condition~(\ref{eq:interface-condition-strong-viscous})
corresponding to the conservation law~(\ref{eq:conservation-strong-viscous})
is the jump or Rankine-Hugoniot condition~\cite{Maj12}, which now
accounts for both the convective and viscous flux, ensuring continuity
of the normal flux at the interface. The interface condition~(\ref{eq:interface-condition-state-strong-viscous})
corresponding to the constitutive law~(\ref{eq:constitutive-strong-viscous})
is unmodified from~(\ref{eq:viscous-interface-condition-state-strong})
by the inclusion of the convective flux. 

The flux $\mathcal{F}\left(y,\sigma\right)$ is defined in terms of
the spatial flux $\mathcal{F}^{x}\left(y,\sigma\right)$ analogously
to~(\ref{eq:spatial-flux}) or~(\ref{eq:spacetime-flux}). The spatial
flux $\mathcal{F}^{x}\left(y,\sigma\right)$ is defined as
\begin{equation}
\mathcal{F}^{x}\left(y,\sigma\right)=\mathcal{F}^{c}\left(y\right)-\mathcal{\tilde{F}}^{v}\left(y,\sigma\right),\label{eq:flux-definition-in-terms-of-auxiliary-variable-flux-formulation}
\end{equation}
where $\tilde{\mathcal{F}}^{v}:\mathbb{R}^{m}\times\mathbb{R}^{m\times d_{x}}\rightarrow\mathbb{R}^{m\times d_{x}}$
is the modified viscous flux defined consistently with the primal
formulation of Section~\ref{subsec:Governing-equations-Primal-Formulation},
\begin{equation}
\mathcal{\tilde{F}}^{v}\left(y,G\left(y\right)\nabla_{x}y\right)=\mathcal{F}^{v}\left(y,\nabla_{x}y\right),\label{eq:modified-viscous-flux-consistency}
\end{equation}
 and $G\left(y\right)\in\mathbb{R}^{m\times d_{x}\times m\times d_{x}}$
is now a generalized \emph{constitutive tensor} that depends on the
specific choice of constitutive law.

One approach to defining the constitutive law is to define a \emph{gradient
formulation}, where the constitutive tensor $G\left(y\right)$ is
taken as the identity,
\begin{equation}
G\left(y\right)\nabla_{x}y=\nabla_{x}y,\label{eq:gradient-formulation-tensor}
\end{equation}
while the viscous flux remains unmodified,
\begin{equation}
\mathcal{\tilde{F}}^{v}\left(y,\sigma\right)=\mathcal{F}^{v}\left(y,\sigma\right).\label{eq:gradient-formulation-modified-viscous-flux}
\end{equation}
The gradient formulation has been used in the context of local discontinuous
Galerkin~\cite{Per06} and hybridized discontinuous Galerkin methods~\cite{Per10}.
This formulation results in a constitutive law~(\ref{eq:constitutive-strong-viscous})
and corresponding interface condition~(\ref{eq:interface-condition-state-strong-viscous})
that are linear with respect to the state variable and do not introduce
a coupling between flow variable components~\cite{Per06}.  In this
case, the interface condition~(\ref{eq:interface-condition-state-strong-viscous})
reduces to 
\begin{equation}
\left\llbracket y\otimes n_{x}\right\rrbracket =0,\label{eq:gradient-formulation-interface-condition}
\end{equation}
which implies $\left\llbracket y\right\rrbracket =0$ at spatial interfaces,
an interface condition arising in the context of elliptic interface
problems~\cite{Han02,Mas12} that directly enforces the continuity
of the state variable. While this choice is reasonable if the solution
is smooth, this approach would not be appropriate for flows that contain
discontinuities in the state variable, such as problems with inviscid
sub-systems, cf. Mott et al.~\cite{Mot20}.

An alternative approach is to define a \emph{flux formulation}, as
in Section~(\ref{sec:visous-strong-formulation}), where the constitutive
tensor $G\left(y\right)$ is defined to be the homogeneity tensor~(\ref{eq:homogeneity-tensor})
so that
\begin{equation}
G\left(y\right)\nabla_{x}y=\mathcal{F}^{v}\left(y,\nabla_{x}y\right)=\mathcal{F}_{\nabla_{x}y}^{v}\left(y,\nabla_{x}y\right)\nabla_{x}y,\label{eq:flux-formulation-tensor}
\end{equation}
while the modified viscous flux is defined to be the auxiliary variable,
\begin{equation}
\mathcal{\tilde{F}}^{v}\left(y,\sigma\right)=\sigma,\label{eq:flux-formulation-modified-viscous-flux}
\end{equation}
recovering a standard mixed method~\cite{Coc98_SIAM}.

A slight modification of the flux formulation for the case of linear
advection-diffusion or viscous Burgers, where $\mathcal{F}^{v}\left(y,\nabla_{x}y\right)=\epsilon\nabla_{x}y$,
is obtained by setting $G\left(y\right)=\sqrt{\epsilon}$ , which
recovers the formulation advocated by Broersen and Stevenson~\cite{Bro14,Bro15}
and later Demkowicz and Gopalakrishnan~\cite{Dem15_20} in the context
of Discontinuous Petrov-Galerkin (DPG) methods for singularly perturbed
problems~\cite{Roo08}. A similar approach was used in the original
description of the LDG method where nonlinear diffusion coefficients
were considered by Cockburn and Shu~\cite{Coc98_SIAM}.

In the case of compressible Navier-Stokes flow, we take an approach
similar to that of Chan et al.~\cite{Cha14} with the scaling advocated
by Broersen and Stevenson also incorporated. The constitutive tensor
$G\left(y\right)$ is defined such that
\begin{equation}
\left(G\left(y\right)\nabla_{x}y\right)_{i}=\mu_{\infty}^{-\nicefrac{1}{2}}\left(0,\tau_{1i},\ldots,\tau_{d_{x}i},-q_{i}\right),\label{eq:navier-stokes-auxiliary-component}
\end{equation}
where $\mu_{\infty}$ is the freestream dynamic viscosity. The viscous
flux is defined in terms of the auxiliary variable as
\begin{equation}
\mathcal{F}_{i}^{v}\left(y,\sigma\right)=\mu_{\infty}^{\nicefrac{1}{2}}\left(\sigma_{1i},\sigma_{2i},\ldots,\sigma_{d_{x}+1i},\sigma_{i+1j}v_{j}+\sigma_{mi}\right).\label{eq:navier-stokes-viscous-flux-auxiliary-component}
\end{equation}
In this way, the auxiliary variable is defined, up to a factor $\mu_{\infty}^{-\nicefrac{1}{2}}$,
as the viscous stress tensor, $\tau$, given by~(\ref{eq:navier-stokes-viscous-stress-tensor-component})
and thermal heat flux, $q$, given by~(\ref{eq:navier-stokes-thermal-heat-flux-component}).
In contrast to Chan~et.~al.~\cite{Cha14} we do not strongly enforce
symmetry of the viscous stress tensor, $\tau$. However, we may explore
this approach in future work as it could lead to a more computationally
efficient and physically accurate formulation.

\subsubsection{Interior and boundary interfaces\label{subsec:Interior-and-Boundary}}

We assume that $\mathcal{E}$ consists of two disjoint subsets: the
interior interfaces 
\begin{equation}
\left\{ \epsilon_{0}\in\mathcal{E}\middlebar\epsilon_{0}\cap\partial\Omega=\emptyset\right\} \label{eq:interior-interfaces}
\end{equation}
 and exterior interfaces 
\begin{equation}
\left\{ \epsilon_{\partial}\in\mathcal{E}\middlebar\epsilon_{\partial}\subset\partial\Omega\right\} ,\label{eq:boundary-interfaces}
\end{equation}
 so that $\mathcal{E}=\mathcal{E}_{0}\cup\mathcal{E}_{\partial}$.
For interior interfaces $\epsilon_{0}\in\mathcal{E}_{0}$ there exists
$\kappa^{+},\kappa^{-}\in\mathcal{T}$ such that $\epsilon_{0}=\partial\kappa^{+}\cap\partial\kappa^{-}$.
On interior interfaces Equations~(\ref{eq:interface-condition-strong-viscous})
and~(\ref{eq:interface-condition-state-strong-viscous}) are defined
as

\begin{align}
\left\llbracket n\cdot\mathcal{F}\left(y,\sigma\right)\right\rrbracket =n^{+}\cdot\mathcal{F}\left(y^{+},\sigma^{+}\right)+n^{-}\cdot\mathcal{F}\left(y^{-},\sigma^{-}\right)=0, & \textup{ on }\epsilon\qquad\forall\epsilon\in\mathcal{E}_{0},\label{eq:interface-condition-strong-interior}\\
\average{G\left(y\right)}\left\llbracket y\otimes n_{x}\right\rrbracket =\frac{1}{2}\left(G\left(y^{+}\right)+G\left(y^{-}\right)\right)\left(y^{+}\otimes n_{x}^{+}+y^{-}\otimes n_{x}^{-}\right)=0, & \textup{ on }\epsilon\qquad\forall\epsilon\in\mathcal{E}_{0}.\label{eq:interface-condition-state-strong-interior}
\end{align}
where $n^{+},n^{-}$ denote the outward facing normal of $\kappa^{+},\kappa^{-}$
respectively, so that $n^{+}=-n^{-}$. For exterior interfaces

\begin{align}
\left\llbracket n\cdot\mathcal{F}\left(y,\sigma\right)\right\rrbracket =n^{+}\cdot\mathcal{F}\left(y^{+},\sigma^{+}\right)-n^{+}\cdot\mathcal{\mathcal{F}}_{\partial}\left(y^{+},\sigma^{+}\right)=0, & \textup{ on }\epsilon\qquad\forall\epsilon\in\mathcal{E}_{\partial},\label{eq:interface-condition-strong-exterior}\\
\average{G\left(y\right)}\left\llbracket y\otimes n_{x}\right\rrbracket =G_{\partial}\left(y^{+}\right)\left(y^{+}\otimes n_{x}^{+}-y_{\partial}\left(y^{+}\right)\otimes n_{x}^{+}\right)=0, & \textup{ on }\epsilon\qquad\forall\epsilon\in\mathcal{E}_{\partial}.\label{eq:interface-condition-state-strong-exterior}
\end{align}

Here $n^{+}\cdot\mathcal{F}_{\partial}\left(y^{+},\sigma^{+}\right)$
is the imposed normal boundary flux, $G_{\partial}\left(y^{+}\right)$
is the boundary modified constitutive tensor, and $y_{\partial}\left(y^{+}\right)$
is the boundary state, which are functions chosen depending on the
type of boundary condition. Therefore, we further decompose $\mathcal{E}_{\partial}$
into disjoint subsets of inflow and outflow interfaces $\mathcal{E}_{\partial}=\mathcal{E}_{\text{in}}\cup\mathcal{E}_{\text{out}},$
so that at an outflow interface $\epsilon_{\textup{out}}$ the boundary
flux is defined as the interior convective flux, and the boundary
state is defined as the interior state,
\begin{align}
n^{+}\cdot\mathcal{F}_{\partial}\left(y^{+},\sigma^{+}\right)=n^{+}\cdot\mathcal{F}\left(y^{+},\sigma_{\text{out}}=0\right), & \textup{ on }\epsilon\qquad\forall\epsilon\in\mathcal{E}_{\text{out}},\label{eq:normal-boundary-flux-outflow}\\
G_{\partial}\left(y^{+}\right)=G\left(y^{+}\right), & \textup{ on }\epsilon\qquad\forall\epsilon\in\mathcal{E}_{\text{out}},\label{eq:boundary-tensor-outflow}\\
y_{\partial}\left(y^{+}\right)=y^{+}, & \textup{ on }\epsilon\qquad\forall\epsilon\in\mathcal{E}_{\text{out}},\label{eq:boundary-state-outflow}
\end{align}
and therefore Equations~(\ref{eq:interface-condition-strong-exterior})
and~(\ref{eq:interface-condition-state-strong-exterior}) are satisfied
trivially. At an inflow boundary $\epsilon_{\textup{in}}\in\mathcal{E}_{\text{in}}$,
the normal convective boundary flux and boundary state are prescribed
values independent of the interior state $y^{+}$ , while the normal
viscous boundary flux is defined as the interior normal viscous flux,
\begin{align}
n^{+}\cdot\mathcal{\mathcal{F}}_{\partial}\left(y^{+},\sigma^{+}\right)=n^{+}\cdot\mathcal{\mathcal{F}}_{\text{in}}^{c}-n^{+}\cdot\mathcal{\tilde{F}}^{v}\left(y^{+},\sigma^{+}\right), & \textup{ on }\epsilon\qquad\forall\epsilon\in\mathcal{E}_{\text{in}},\label{eq:normal-boundary-flux-inflow}\\
G_{\partial}\left(y^{+}\right)=G\left(y_{\partial}\left(y^{+}\right)\right), & \textup{ on }\epsilon\qquad\forall\epsilon\in\mathcal{E}_{\text{in}},\label{eq:boundary-tensor-inflow}\\
y_{\partial}\left(y^{+}\right)=y_{\text{in}}, & \textup{ on }\epsilon\qquad\forall\epsilon\in\mathcal{E}_{\text{in}}.\label{eq:boundary-state-inflow}
\end{align}

\subsubsection{Weak formulation\label{sec:Weak-formulation}}

A weak formulation in physical space is obtained by integrating the
conservation law~(\ref{eq:conservation-strong-viscous}), the constitutive
law~(\ref{eq:constitutive-strong-viscous}), and the corresponding
interface conditions~(\ref{eq:interface-condition-strong-viscous}),
(\ref{eq:interface-condition-state-strong-viscous}) for each element
and interface against separate test functions: find $\left(y,\sigma\right)\in Y\times\Sigma$
such that
\begin{align}
0= & \;\;\;\;\sum_{\kappa\in\mathcal{T}}\left(\nabla\cdot\mathcal{F}\left(y,\sigma\right)-f,v\right)_{\kappa}\nonumber \\
 & +\sum_{\kappa\in\mathcal{T}}\left(\sigma-G(y)\nabla_{x}y,\tau\right)_{\kappa}\nonumber \\
 & -\sum_{\epsilon\in\mathcal{E}}\left(\left\llbracket n\cdot\mathcal{F}\left(y,\sigma\right)\right\rrbracket ,w_{y}\right)_{\epsilon}\nonumber \\
 & -\sum_{\epsilon\in\mathcal{E}}\left(\average{G\left(y\right)}\left\llbracket y\otimes n_{x}\right\rrbracket ,w_{\sigma}\right)_{\epsilon}\qquad\forall\left(v,\tau,w_{y},w_{\sigma}\right)\in V_{y}\times V_{\sigma}\times W_{y}\times W_{\sigma}.\label{eq:weak-formulation-viscous}
\end{align}
The solution spaces $Y$ and $\Sigma$ are the broken Sobolev spaces,
\begin{eqnarray}
Y & = & \left\{ y\in\left[L^{2}\left(\Omega\right)\right]^{m\hphantom{\times d_{x}}}\bigl|\forall\kappa\in\mathcal{T},\:\:\hphantom{\nabla_{x}\cdot}\left.y\right|_{\kappa}\in\left[H^{1}\left(\kappa\right)\right]^{m}\right\} ,\label{eq:solution-space-vector}\\
\Sigma & = & \left\{ \sigma\in\left[L^{2}\left(\Omega\right)\right]^{m\times d_{x}}\bigl|\forall\kappa\in\mathcal{T},\left.\nabla_{x}\cdot\sigma\right|_{\kappa}\in\left[L^{2}\left(\Omega\right)\right]^{m}\right\} ,\label{eq:solution-space-tensor}
\end{eqnarray}
while the test spaces are defined as $V_{y}=\left[L^{2}\left(\Omega\right)\right]^{m}$
and $V_{\sigma}=\left[L^{2}\left(\Omega\right)\right]^{m\times d}$,
with $W_{y}$ and $W_{\sigma}$ defined to be the corresponding single-valued
trace spaces, cf. Carstensen et al.~\cite{Car16}.

\subsection{Formulation in reference space with variable geometry\label{subsec:variable-geometry}}

Analogous to our previous work~\cite{Cor18}, the grid must be treated
as a variable in order to align discrete grid interfaces with flow
interfaces or more generally to move the grid to resolve under-resolved
flow features. Therefore, we transform the strong formulation~(\ref{eq:conservation-strong-viscous}),~(\ref{eq:constitutive-strong-viscous}),~(\ref{eq:interface-condition-strong-viscous}),~(\ref{eq:interface-condition-state-strong-viscous})
and weak formulation~(\ref{eq:weak-formulation-viscous}) of the
flow equations from physical to reference coordinates in order to
facilitate differentiation with respect to geometry.

\subsubsection{Mapping from reference space\label{subsec:mapping-reference-space}}

We assume that there is a continuous, invertible mapping 
\begin{equation}
u:\hat{\Omega}\rightarrow\Omega,\label{eq:shape-mapping}
\end{equation}
from a reference domain $\hat{\Omega}\subset\mathbb{R}^{d}$ to the
physical domain $\Omega\subset\mathbb{R}^{d}$. We assume that $\hat{\Omega}$
is partitioned by $\hat{\mathcal{T}}$, so that $\overline{\hat{\Omega}}=\cup_{\hat{\kappa}\in\hat{\mathcal{T}}}\overline{\hat{\kappa}}$.
Also, we consider the set of interfaces $\hat{\mathcal{E}}$ consisting
of disjoint interfaces $\hat{\epsilon}$, such that $\cup_{\hat{\epsilon}\in\hat{\mathcal{E}}}\hat{\epsilon}=\cup_{\hat{\kappa}\in\hat{\mathcal{T}}}\partial\hat{\kappa}$.
The mapping $u$ is further assumed to be (piecewise) differentiable
with derivative or Jacobian matrix denoted
\begin{equation}
\left.\nabla u\right|_{\hat{\kappa}}:\hat{\kappa}\rightarrow\mathbb{R}^{d\times d}\qquad\forall\hat{\kappa}\in\hat{\mathcal{T}}.\label{eq:jacobian}
\end{equation}
The cofactor matrix $\left.\cof\left(\nabla u\right)\right|_{\hat{\kappa}}:\hat{\kappa}\rightarrow\mathbb{R}^{d\times d},$
is defined for $\hat{\kappa}\in\hat{\mathcal{T}}$,
\begin{equation}
\cof\left(\nabla u\left(\hat{x}\right)\right)=\det\left(\nabla u\left(\hat{x}\right)\right)\left(\nabla u\left(\hat{x}\right)\right)^{-\top}\qquad\forall\hat{x}\in\hat{\kappa},\label{eq:cofactor}
\end{equation}
where $\left.\det\left(\nabla u\right)\right|_{\hat{\kappa}}:\hat{\kappa}\rightarrow\mathbb{R}$
is the determinant of the Jacobian.

As detailed in our related work~\cite{Cor20_Ref}, assuming that
$y$ and $v$ are functions over reference space, the weak formulation
of a conservation law in physical space can be evaluated in reference
space according to
\begin{equation}
\left(\nabla\cdot\mathcal{F}\left(y\circ u^{-1}\right),v\circ u^{-1}\right)_{\kappa}=\left(\left(\cof\left(\nabla u\right)\nabla\right)\cdot\mathcal{F}\left(y\right),v\right)_{\hat{\kappa}}.\label{eq:weak-form-physical-space-reference-space}
\end{equation}
Likewise, treating $\sigma$ and $\tau$ as functions over reference
space, the constitutive law can be evaluated in reference space according
to
\begin{equation}
\left(\sigma\circ u^{-1}-G(y\circ u^{-1})\nabla_{x}\left(y\circ u^{-1}\right),\tau\circ u^{-1}\right)_{\kappa}=\left(\det\left(\nabla u\right)\sigma-G(y)\left(\cof\left(\nabla u\right)\nabla\right)_{x}y,\tau\right)_{\hat{\kappa}},\label{eq:weak-form-constitutive-law-physical-space-reference-space}
\end{equation}
In order to represent the spatial gradient in a space-time setting
($d=d_{x}+1$), we define $\left(\cof\left(\nabla u\right)\nabla\right)_{x}$
to be the spatial components of $\left(\cof\left(\nabla u\right)\nabla\right)$,
so that if $\left(\cof\left(\nabla u\right)\nabla\right)y:\hat{\kappa}\rightarrow\mathbb{R}^{m\times d}$
then $\left(\cof\left(\nabla u\right)\nabla\right)y:\hat{\kappa}\rightarrow\mathbb{R}^{m\times d_{x}}$,
while in a spatial ($d=d_{x}$) setting $\left(\cof\left(\nabla u\right)\nabla\right)=\left(\cof\left(\nabla u\right)\nabla\right)_{x}$.

The weak formulation of each interface condition can similarly be
evaluated in reference space according to

\begin{equation}
\left(\left\llbracket n\cdot\mathcal{F}\left(y\circ u^{-1}\right)\right\rrbracket ,w_{y}\circ u^{-1}\right)_{\epsilon}=\left(\left\llbracket s\left(\nabla u\right)\cdot\mathcal{F}\left(y\right)\right\rrbracket ,w_{y}\right)_{\hat{\epsilon}}\label{eq:weak-form-flux-jump-physical-space-reference-space}
\end{equation}
and
\begin{equation}
\left(\left(\average{G\left(y\circ u^{-1}\right)}\left\llbracket \left(y\circ u^{-1}\right)\otimes n_{x}\right\rrbracket \right),w_{\sigma}\circ u^{-1}\right)_{\epsilon}=\left(\average{G\left(y\right)}\left\llbracket y\otimes s\left(\nabla u\right)_{x}\right\rrbracket ,w_{\sigma}\right)_{\hat{\epsilon}}\label{eq:weak-form-state-jump-physical-space-reference-space}
\end{equation}
In this setting, $\hat{\epsilon}\in\hat{\mathcal{E}}$, $\epsilon=u\left(\hat{\epsilon}\right)\in\mathcal{E}$,
and $n:\epsilon\rightarrow\mathbb{R}^{d}$ is the unit normal, which
can be evaluated in terms of $u$ according to 
\begin{equation}
n=\left(\frac{s\left(\nabla u\right)}{\left\Vert s\left(\nabla u\right)\right\Vert }\right)\circ u^{-1},\label{eq:unit-normal-in-terms-of-scaled-normal}
\end{equation}
where $s\left(\nabla u\right):\hat{\epsilon}\rightarrow\mathbb{R}^{d}$
is defined as follows. We assume that there exists a parameterization
$\theta_{\hat{\epsilon}}:\hat{D}\rightarrow\hat{\epsilon}$, mapping
from points $\left(\xi_{1},\ldots,\xi_{d-1}\right)$ in parameter
space $\hat{D}\subset\mathbb{R}^{d-1}$ to points $\hat{x}$ on the
reference space interface, such that the reference space tangent plane
basis vectors $\partial_{\xi_{1}}\theta_{\hat{\epsilon}},\ldots,\partial_{\xi_{d-1}}\theta_{\hat{\epsilon}}$
are of unit magnitude. A parameterization of $\epsilon=u\left(\hat{\epsilon}\right)$
is then given by the composition $\theta_{\epsilon}=u\circ\theta_{\hat{\epsilon}}:\hat{D}\rightarrow\epsilon$.
Given $\hat{\epsilon}\in\hat{\mathcal{E}}$, the scaled normal $\left.s\left(\nabla u\right)\right|_{\hat{\epsilon}}:\hat{\epsilon}\rightarrow\mathbb{R}^{d}$
is defined for $\hat{x}\in\hat{\epsilon}$ as the scaled normal of
the tangent plane of $\epsilon$ corresponding to the parameter $\theta_{\hat{\epsilon}}^{-1}\left(\hat{x}\right)$.
If $d=3$ and $\xi,\eta$ denote the parametric coordinates, then
\begin{equation}
\left.s\left(\nabla u\right)\right|_{\hat{\epsilon}}=\left(\partial_{\xi}\theta_{\epsilon}\times\partial_{\eta}\theta_{\epsilon}\right)\circ\theta_{\hat{\epsilon}}^{-1}\label{eq:scaled-face-normal-3d}
\end{equation}
where $\partial_{\xi}\theta_{\epsilon}\times\partial_{\eta}\theta_{\epsilon}$
is the cross product of the tangent plane basis vectors.

A general formula for evaluating the cross product of tangent plane
basis vectors is given by the following: let $\left(\boldsymbol{x}_{1},\ldots,\boldsymbol{x}_{d}\right)$
denote the coordinate directions in $\mathbb{R}^{d}$, and the parameterization
be given in terms of components $\theta_{\epsilon}=\left(\theta_{\epsilon}^{1},\ldots,\theta_{\epsilon}^{d}\right),$
then
\begin{equation}
\partial_{\xi_{1}}\theta_{\epsilon}\times\cdots\times\partial_{\xi_{d-1}}\theta_{\epsilon}=\det\left(\begin{array}{ccc}
\partial_{\xi_{1}}\theta_{\epsilon}^{1} & \cdots & \partial_{\xi_{1}}\theta_{\epsilon}^{d}\\
\vdots & \ddots & \vdots\\
\partial_{\xi_{d-1}}\theta_{\epsilon}^{1} & \cdots & \partial_{\xi_{d-1}}\theta_{\epsilon}^{d}\\
\boldsymbol{x}_{1} & \ldots & \boldsymbol{x}_{d}
\end{array}\right).\label{eq:generalized-cross-product}
\end{equation}
By the chain rule we can express $\partial_{\xi_{i}}\theta_{\epsilon}$
in terms of $\nabla u$,
\begin{equation}
\partial_{\xi_{i}}\theta_{\epsilon}\left(\xi\right)=\nabla u\left(\theta_{\hat{\epsilon}}\left(\xi\right)\right)\cdot\partial_{\xi_{i}}\theta_{\hat{\epsilon}}\left(\xi\right),\label{eq:parameterization-tangent-basis-in-terms-of-gradient}
\end{equation}
so that in general the physical space scaled normal as a function
of $\nabla u$ is
\begin{equation}
\left.s\left(\nabla u\right)\right|_{\hat{\epsilon}}=\left(\partial_{\xi_{1}}\theta_{\epsilon}\times\cdots\times\partial_{\xi_{d-1}}\theta_{\epsilon}\right)\circ\theta_{\hat{\epsilon}}^{-1}.\label{eq:scaled-face-normal-nd}
\end{equation}

In the present work, we have adopted the more standard convention
in the definition of the generalized cross product given by Equation~(\ref{eq:generalized-cross-product}),
which is used to define the generalized scaled normal given by Equation~(\ref{eq:scaled-face-normal-nd}).
This definition differs from Equation~(3.33) of our previous work~\cite{Cor18}
by a factor of $\left(-1\right)^{\left(d-1\right)}$ in order to ensure
that (\ref{eq:generalized-cross-product}) and (\ref{eq:scaled-face-normal-nd})
are positively oriented, cf. Massey~\cite{Mas83}. 

\subsubsection{Strong and weak formulation in reference space\label{subsec:reference-space-formulation}}

The strong form in reference space is

\begin{align}
\left(\cof\left(\nabla u\right)\nabla\right)\cdot\mathcal{F}\left(y,\sigma\right)=0 & \textup{ in }\hat{\kappa}\qquad\forall\hat{\kappa}\in\hat{\mathcal{T}},\label{eq:conservation-strong-reference-viscous}\\
\det\left(\nabla u\right)\sigma-G(y)\left(\cof\left(\nabla u\right)\nabla\right)_{x}y=0 & \textup{ in }\hat{\kappa}\qquad\forall\hat{\kappa}\in\hat{\mathcal{T}},\label{eq:constitutive-strong-reference-viscous}\\
\left\llbracket s\left(\nabla u\right)\cdot\mathcal{F}\left(y,\sigma\right)\right\rrbracket =0 & \textup{ on }\hat{\epsilon}\qquad\forall\hat{\epsilon}\in\hat{\mathcal{E}},\label{eq:interface-condition-strong-reference-viscous}\\
\average{G\left(y\right)}\left\llbracket y\otimes s\left(\nabla u\right)_{x}\right\rrbracket =0 & \textup{ on }\hat{\epsilon}\qquad\forall\hat{\epsilon}\in\hat{\mathcal{E}},\label{eq:interface-condition-state-strong-reference-viscous}\\
b\left(u\right)-u=0 & \textup{ on }\hat{\epsilon}\qquad\forall\hat{\epsilon}\in\hat{\mathcal{E}},\label{eq:geometric-boundary-condition}
\end{align}
where $\nabla u$ is the Jacobian of the mapping from reference to
physical space, $\det\left(\nabla u\right)$ is its determinant, $\cof\left(\nabla u\right)$
is its cofactor matrix, and $s\left(\nabla u\right)$ is the scaled
normal as defined in Section~\ref{subsec:mapping-reference-space}.
Equation~(\ref{eq:geometric-boundary-condition}) imposes geometric
boundary conditions that constrain points to the boundary of the physical
domain via a projection operator $b:U\rightarrow U$, where $U=\left[H^{1}\left(\hat{\Omega}\right)\right]^{d}$is
the $\mathbb{R}^{d}$-valued Sobolev space over $\hat{\Omega}$. Examples
of $b\left(u\right)$ are given in earlier work~\cite{Cor18,Cor19_AVIATION}.
We assume that $Y$ and $\Sigma$, originally defined for functions
over physical space, cf.~(\ref{eq:solution-space-vector}) and~(\ref{eq:solution-space-tensor}),
now consist, respectively, of functions defined in $\mathbb{R}^{m}$-valued
and $\mathbb{R}^{m\times d_{x}}$-valued broken Sobolev spaces over
$\hat{\mathcal{T}}$. We further assume that the test spaces $V_{y}$,
$V_{\sigma},W_{y},W_{\sigma}$ now consist of functions defined over
reference space.

We define a provisional state operator $\tilde{e}:Y\times\Sigma\times U\rightarrow\left(V_{y}\times V_{\sigma}\times W_{y}\times W_{\sigma}\right)^{*}$
for $\left(y,\sigma,u\right)\in Y\times\Sigma\times U$, by
\begin{align}
\tilde{e}\left(y,\sigma,u\right)=\left(v,\tau,w_{y},w_{\sigma}\right)\mapsto & \;\;\;\;\sum_{\hat{\kappa}\in\hat{\mathcal{T}}}\left(\left(\cof\left(\nabla u\right)\nabla\right)\cdot\mathcal{F}\left(y,\sigma\right),v\right)_{\hat{\kappa}}\nonumber \\
 & +\sum_{\hat{\kappa}\in\hat{\mathcal{T}}}\left(\det\left(\nabla u\right)\sigma-G(y)\left(\cof\left(\nabla u\right)\nabla\right)_{x}y,\tau\right)_{\hat{\kappa}}\nonumber \\
 & -\sum_{\hat{\epsilon}\in\hat{\mathcal{E}}}\left(\left\llbracket s\left(\nabla u\right)\cdot\mathcal{F}\left(y,\sigma\right)\right\rrbracket ,w_{y}\right)_{\hat{\epsilon}}\nonumber \\
 & -\sum_{\hat{\epsilon}\in\hat{\mathcal{E}}}\left(\average{G\left(y\right)}\left\llbracket y\otimes s\left(\nabla u\right)_{x}\right\rrbracket ,w_{\sigma}\right)_{\hat{\epsilon}}\label{eq:state-operator-weak-formulation-reference-viscous}
\end{align}
which has a Fréchet derivative defined for perturbation $\left(\delta y,\delta\sigma,\delta u\right)\in Y\times\Sigma\times U$,
and test functions $\left(v,\tau,w_{y},w_{\sigma}\right)\in V_{y}\times V_{\sigma}\times W_{y}\times W_{\sigma}$,
by its partial derivative with respect to the state variable $y$,
\begin{align}
\tilde{e}_{y}\left(y,\sigma,u\right)\delta y=\left(v,\tau,w_{y},w_{\sigma}\right)\mapsto & \;\;\;\;\sum_{\hat{\kappa}\in\hat{\mathcal{T}}}\left(\left(\cof\left(\nabla u\right)\nabla\right)\cdot\left(\mathcal{F}_{y}\left(y,\sigma\right)\delta y\right),v\right)_{\hat{\kappa}}\nonumber \\
 & +\sum_{\hat{\kappa}\in\hat{\mathcal{T}}}\left(-\left(\left(G'(y)\delta y\right)\left(\cof\left(\nabla u\right)\nabla\right)_{x}\delta y+G(y)\left(\cof\left(\nabla u\right)\nabla\right)_{x}\delta y\right),\tau\right)_{\hat{\kappa}}\nonumber \\
 & -\sum_{\hat{\epsilon}\in\hat{\mathcal{E}}}\left(\left\llbracket s\left(\nabla u\right)\cdot\left(\mathcal{F}_{y}\left(y,\sigma\right)\delta y\right)\right\rrbracket ,w_{y}\right)_{\hat{\epsilon}}\nonumber \\
 & -\sum_{\hat{\epsilon}\in\hat{\mathcal{E}}}\left(\average{G'\left(y\right)\delta y}\left\llbracket y\otimes s\left(\nabla u\right)_{x}\right\rrbracket +\average{G\left(y\right)}\left\llbracket \delta y\otimes s\left(\nabla u\right)_{x}\right\rrbracket ,w_{\sigma}\right)_{\hat{\epsilon}},\label{eq:state-operator-state-derivative-weak-formulation-reference-viscous}
\end{align}
its partial derivative with respect to the auxiliary variable $\sigma$,
\begin{align}
\tilde{e}_{\sigma}\left(y,\sigma,u\right)\delta\sigma=\left(v,\tau,w_{y},w_{\sigma}\right)\mapsto & \;\;\;\;\sum_{\hat{\kappa}\in\hat{\mathcal{T}}}\left(\left(\cof\left(\nabla u\right)\nabla\right)\cdot\left(\mathcal{F}_{\sigma}\left(y,\sigma\right)\delta\sigma\right),v\right)_{\hat{\kappa}}\nonumber \\
 & +\sum_{\hat{\kappa}\in\hat{\mathcal{T}}}\left(\det\left(\nabla u\right)\delta\sigma,\tau\right)_{\hat{\kappa}}\nonumber \\
 & -\sum_{\hat{\epsilon}\in\hat{\mathcal{E}}}\left(\left\llbracket s\left(\nabla u\right)\cdot\left(\mathcal{F}_{\sigma}\left(y,\sigma\right)\delta\sigma\right)\right\rrbracket ,w_{y}\right)_{\hat{\epsilon}},\label{eq:state-operator-auxiliary-derivative-weak-formulation-reference-viscous}
\end{align}
and its partial derivative with respect to the geometry variable $u$,
\begin{align}
\tilde{e}_{u}\left(y,\sigma,u\right)\delta u=\left(v,\tau,w_{y},w_{\sigma}\right)\mapsto & \;\;\;\;\sum_{\hat{\kappa}\in\hat{\mathcal{T}}}\left(\left(\left(\cof'\left(\nabla u\right)\nabla\delta u\right)\nabla\right)\cdot\mathcal{F}\left(y,\sigma\right),v\right)_{\hat{\kappa}}\nonumber \\
 & +\sum_{\hat{\kappa}\in\hat{\mathcal{T}}}\left(\left(\det'\left(\nabla u\right)\nabla\delta u\right)\sigma-G(y)\left(\left(\cof'\left(\nabla u\right)\nabla\delta u\right)\nabla\right)_{x}y,\tau\right)_{\hat{\kappa}}\nonumber \\
 & -\sum_{\hat{\epsilon}\in\mathcal{\hat{E}}}\left(\left\llbracket \left(s'\left(\nabla u\right)\nabla\delta u\right)\cdot\mathcal{F}\left(y,\sigma\right)\right\rrbracket ,w_{y}\right)_{\hat{\epsilon}}\nonumber \\
 & -\sum_{\hat{\epsilon}\in\hat{\mathcal{E}}}\left(\average{G\left(y\right)}\left\llbracket y\otimes\left(s'\left(\nabla u\right)\nabla\delta u\right)_{x}\right\rrbracket ,w_{\sigma}\right)_{\hat{\epsilon}}.\label{eq:state-operator-geometry-derivative-weak-formulation-reference-viscous}
\end{align}
The state operator $e:Y\times\Sigma\times U\rightarrow\left(V_{y}\times V_{\sigma}\times W_{y}\times W_{\sigma}\right)^{*}$,
which imposes geometric boundary conditions (\ref{eq:geometric-boundary-condition})
by composing the provisional state operator (\ref{eq:state-operator-weak-formulation-reference-viscous})
with the projection $b\left(u\right)$, is defined by
\begin{equation}
e\left(y,\sigma,u\right)=\tilde{e}\left(y,\sigma,b\left(u\right)\right),\label{eq:state-operator-composition}
\end{equation}
with Fréchet derivative defined, for state $\left(y,\sigma,u\right)\in Y\times\Sigma\times U$
and perturbation $\left(\delta y,\delta\sigma,\delta u\right)\in Y\times\Sigma\times U$,
by
\begin{equation}
e'\left(y,\sigma,u\right)=\left(\delta y,\delta\sigma,\delta u\right)\mapsto\tilde{e}_{y}\left(y,\sigma,b\left(u\right)\right)\delta y+\tilde{e}_{\sigma}\left(y,\sigma,b\left(u\right)\right)\delta\sigma+\tilde{e}_{u}\left(y,\sigma,b\left(u\right)\right)b'\left(u\right)\delta u.\label{eq:state-operator-composition-derivative}
\end{equation}
The state equation in reference coordinates is $e\left(y,\sigma,u\right)=0.$
The corresponding weak formulation in reference coordinates is: find
$\left(y,\sigma,u\right)\in Y\times\Sigma\times U$ such that
\begin{equation}
\left\langle e\left(y,\sigma,u\right),\left(v,\tau,w_{y},w_{\sigma}\right)\right\rangle =0\qquad\forall\left(v,\tau,w_{y},w_{\sigma}\right)\in V_{y}\times V_{\sigma}\times W_{y}\times W_{\sigma},\label{eq:weak-formulation-3-reference-coordinates}
\end{equation}
so that the solution satisfying (\ref{eq:conservation-strong-reference-viscous})
and (\ref{eq:interface-condition-strong-reference-viscous}) weakly
and (\ref{eq:geometric-boundary-condition}) strongly is therefore
given as $\left(y,\sigma,b\left(u\right)\right)\in Y\times\Sigma\times U$.

\subsection{Discretization\label{subsec:Discretization}}

We choose discrete (finite-dimensional) subspaces $Y_{h}\subset Y$,
$\Sigma_{h}\subset\Sigma$, $U_{h}\subset U$, $V_{y,h}\subset V_{y}$,
$V_{\sigma,h}\subset V_{\sigma}$, $W_{y,h}\subset W_{y}$, and $W_{\sigma,h}\subset W_{\sigma}$
to discretize the weak formulation (\ref{eq:weak-formulation-3-reference-coordinates}),
which is restricted to the discrete subspaces via the discrete state
operator,
\begin{equation}
e_{h}:Y_{h}\times\Sigma_{h}\times U_{h}\rightarrow\left(V_{y,h}\times V_{\sigma,h}\times W_{y,h}\times W_{\sigma,h}\right)^{*}\label{eq:discrete-state-operator}
\end{equation}
defined such that $e_{h}\left(y,\sigma,u\right)=e\left(y,\sigma,u\right)$
for all $\left(y,\sigma,u\right)\in Y_{h}\times\Sigma_{h}\times U_{h}$
and the $h$-subscript indicates that discrete subspaces have been
selected.

We use standard piecewise polynomials, cf.~\cite{Har13}, defined
over reference elements. Let $\mathcal{P}_{p}$ denote the space of
polynomials spanned by the monomials $\boldsymbol{x}^{\alpha}$ with
multi-index $\alpha\in\mathbb{N}_{0}^{d}$ , satisfying $\sum_{i=1}^{d}\alpha_{i}\leq p$.
In the case of a simplicial grid,
\begin{eqnarray}
Y_{h} & = & \left\{ y\in Y\middlebar\forall\hat{\kappa}\in\hat{\mathcal{T}},\left.y\right|_{\hat{\kappa}}\in\left[\mathcal{P}_{p}\right]^{m\hphantom{\times d_{x}}}\right\} ,\label{eq:discrete-solution-state-space-simplex}\\
\Sigma_{h} & = & \left\{ \sigma\in\Sigma\middlebar\forall\hat{\kappa}\in\hat{\mathcal{T}},\left.\sigma\right|_{\hat{\kappa}}\in\left[\mathcal{P}_{p}\right]^{m\times d_{x}}\right\} .\label{eq:discrete-solution-flux-space-simplex}
\end{eqnarray}
The polynomial degree of the state space and flux space are in general
distinct. In the present work, we choose $V_{y,h}=Y_{h}$, $V_{\sigma,h}=\Sigma_{h}$,
while $W_{y,h}$, and $W_{\sigma,h}$ are chosen to be the corresponding
single-valued polynomial trace spaces. While the present approach
is a discrete least squares method~\cite{Kei17} with a priori chosen
test spaces, future work will investigate a least squares finite element
formulation~\cite{Boc98,Boc09} with optimal test spaces automatically
generated using the discontinuous Petrov–Galerkin methodology of Demkowicz
and Gopalakrishnan~\cite{Dem10,Dem11,Dem15_20}.

The discrete subspace $U_{h}$ of mappings from reference space to
physical space are also discretized into $\mathbb{R}^{d}$-valued
piecewise polynomials, in the case of a simplicial grid
\begin{equation}
U_{h}=\left\{ u\in U\middlebar\forall\hat{\kappa}\in\hat{\mathcal{T}},\left.u\right|_{\hat{\kappa}}\in\left[\mathcal{P}_{p}\right]^{d}\right\} .\label{eq:discrete-shape-space-simplex}
\end{equation}
The case that the chosen polynomial degree of $U_{h}$ is equal to
that of $Y_{h}$ is referred to as isoparametric. It is also possible
to choose the polynomial degree of $U_{h}$ to be less (sub-parametric)
or greater (super-parametric) than that of $Y_{h}$. 

\subsection{Solver\label{subsec:Solver}}

In general, the dimensionality of the discrete solution space and
discrete residual space do not match. Therefore, the weak formulation
is solved iteratively using unconstrained optimization to minimize
$\frac{1}{2}\left\Vert e_{h}\left(y,\sigma,u\right)\right\Vert ^{2}$,
by seeking a stationary point\footnote{The stationary point~(\ref{eq:discrete-residual-stationary}) and
Newton's method~(\ref{eq:discrete-increment-newton}) were stated
incorrectly in previous work, cf.~\cite{Cor18}, Equations~(77)
and~(80).},

\begin{equation}
e_{h}'\left(y,\sigma,u\right)^{*}e_{h}\left(y,\sigma,u\right)=0.\label{eq:discrete-residual-stationary}
\end{equation}
Given an initialization $\left(y,\sigma,u\right)_{0}$ the solution
is repeatedly updated
\begin{equation}
\left(y,\sigma,u\right)_{i+1}=\left(y,\sigma,u\right)_{i}+\Delta\left(y,\sigma,u\right)_{i}\qquad i=0,1,2,\ldots\label{eq:discrete-iteration}
\end{equation}
until (\ref{eq:discrete-residual-stationary}) is satisfied to a given
tolerance. One approach is to use Newton's method, which is a second-order
method with increment given by,
\begin{equation}
\Delta\left(y,\sigma,u\right)=-\left(\left(e_{h}''\left(y,\sigma,u\right)^{*}\cdot\right)e_{h}\left(y,\sigma,u\right)+e_{h}'\left(y,\sigma,u\right)^{*}e_{h}'\left(y,\sigma,u\right)\right)^{-1}\left(e_{h}'\left(y,\sigma,u\right)^{*}e_{h}\left(y,\sigma,u\right)\right).\label{eq:discrete-increment-newton}
\end{equation}
Alternatively, the Gauss-Newton method neglects second derivatives,
yet recovers the second-order convergence rate of Newton's method
as the residual vanishes and ensures a positive semi-definite matrix,
resulting in an increment given by
\begin{equation}
\Delta\left(y,\sigma,u\right)=-\left(e_{h}'\left(y,\sigma,u\right)^{*}e_{h}'\left(y,\sigma,u\right)\right)^{-1}\left(e_{h}'\left(y,\sigma,u\right)^{*}e_{h}\left(y,\sigma,u\right)\right).\label{eq:discrete-increment-gauss-newtoon}
\end{equation}
We employ a Levenberg-Marquardt method to solve~(\ref{eq:discrete-residual-stationary}),
which augments the Gauss-Newton method~(\ref{eq:discrete-increment-gauss-newtoon})
with a regularization term,
\begin{equation}
\Delta\left(y,\sigma,u\right)=-\left(e_{h}'\left(y,\sigma,u\right)^{*}e_{h}'\left(y,\sigma,u\right)+I_{\lambda}\left(y,\sigma,u\right)\right)^{-1}\left(e_{h}'\left(y,\sigma,u\right)^{*}e_{h}\left(y,\sigma,u\right)\right),\label{eq:discrete-increment}
\end{equation}
where the regularization operator, 
\begin{equation}
I_{\lambda}\left(y,\sigma,u\right):\left(\delta y,\delta\sigma,\delta u\right)\mapsto\left(\lambda_{y}\delta y,\lambda_{\sigma}\delta\sigma,\lambda_{u}\delta u\right),\label{eq:regularization-operator}
\end{equation}
ensures invertibility and therefore positive definiteness of the linear
system of equations. Separate regularization coefficients $\lambda_{y},\lambda_{\sigma},\lambda_{u}\ge0$
are defined for each solution variable. In practice, the state and
auxiliary regularization coefficients can be set to zero, $\lambda_{y}=\lambda_{\sigma}=0$,
while the grid regularization coefficient $\lambda_{u}>0$ must be
positive in order to ensure rank sufficiency and to limit excessive
grid motion. Additional symmetric positive definite operators can
be incorporated into the regularization operator~\cite{Cor19_AVIATION}.
In the present work, we incorporate a linear elastic grid regularization,
which is a symmetric positive definite operator that has the effect
of distributing the grid motion to neighboring elements. The linear
elastic grid regularization is a variation of the Laplacian grid regularization,
$\delta u\mapsto-\lambda_{\Delta u}\left(b'\left(u\right)^{*}\Delta b'\left(u\right)\right)\delta u$,
with $\lambda_{\Delta u}\geq0$, that we employed in previous work~\cite{Cor19_AVIATION}
that offers the added benefit of introducing a compressibility effect
into the grid motion that we have found useful for resolving thin
viscous layers. Other possible regularization strategies include the
weighted elliptic regularization proposed by Zahr et al.~\cite{Zah19,Zah20_SCITECH}.
The resulting linear system of equations is positive definite and
symmetric. In the present work, we employ a sparse direct solver provided
by Eigen~\cite{Eig10}. 

The grid topology may need to be modified by the solver in order to
fit a priori unknown interfaces and ensure element validity while
resolving sharp gradients, for which we employ standard edge refinement
and edge collapse algorithms~\cite{Loh08}. In the present work,
element quality is used as an indicator for local refinement. Elements
that become highly anisotropic as MDG-ICE moves the grid to resolve
thin viscous structures are adaptively split by refining their longest
edge. In the case of nonlinear elements, if the determinant of the
Jacobian at any degree of freedom is negative, we apply a control
that projects the elements to a linear shape representation and locally
refines the element if projecting the cell does not recover a valid
grid. Often, the introduction of the additional grid topology and
resolution is sufficient for the solver to recover a valid grid. The
solver does not currently incorporate any other grid smoothing or
optimization terms based on element quality~\cite{Zah18_AIAA,Zah18,Zah19}.

\section{Results\label{sec:Results}}

We now apply MDG-ICE to compute steady and unsteady solutions for
flows containing sharp, smooth, gradients. Solutions to unsteady problems
are solved using a space-time formulation. Unless otherwise indicated,
the grid is assumed to consist of isoparametric elements, see Section~\ref{subsec:Discretization}.

\subsection{Linear advection-diffusion\label{sec:linear-advection-diffusion-1d}}

We consider steady, one-dimensional linear advection-diffusion described
in Section~\ref{subsec:Linear-advection-diffusion}, subject to the
following boundary conditions

\begin{align}
y\left(x=0\right) & =0,\nonumber \\
y\left(x=1\right) & =1.\label{eq:steady-boundary-layer-boundary-conditions}
\end{align}
The exact solution is given by

\begin{equation}
y\left(x\right)=\frac{1-\exp\left(x\cdot\mathrm{Pe}\right)}{1-\exp\left(\mathrm{Pe}\right)}\label{eq:steady-boundary-layer-exact-solution}
\end{equation}
where $\mathrm{Pe}=\frac{1}{\varepsilon}=\frac{v\ell}{\mu}$ is the
Péclet number, $v$ is the characteristic velocity, $\ell$ is the
characteristic length, and $\mu$ is the mass diffusivity. In this
case, the solution exhibits a boundary layer like profile at $x=1$~\cite{Mas13}.
\begin{figure}
\centering{}%
\begin{tabular}{cc}
\includegraphics[width=0.4\linewidth]{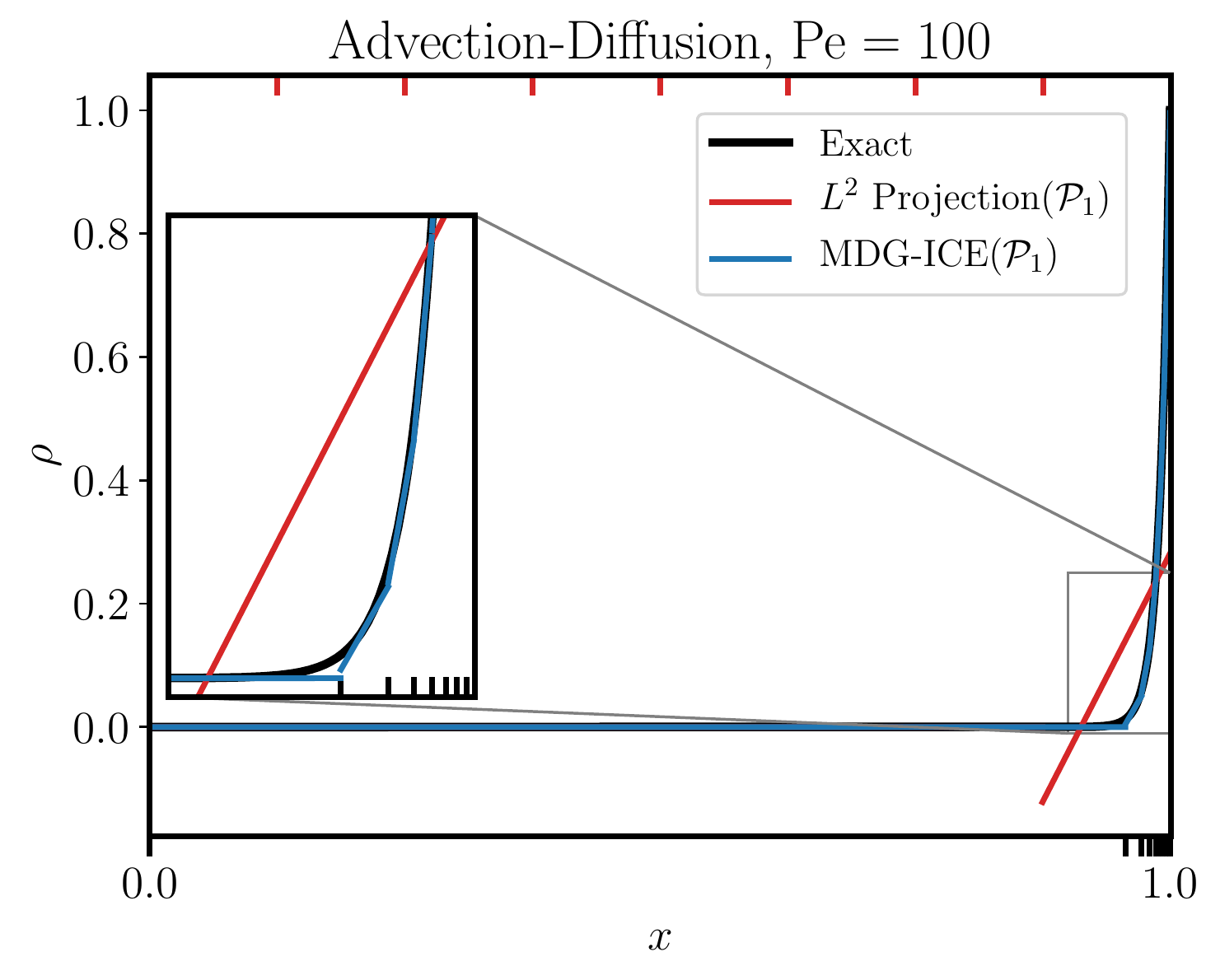} & \includegraphics[width=0.4\linewidth]{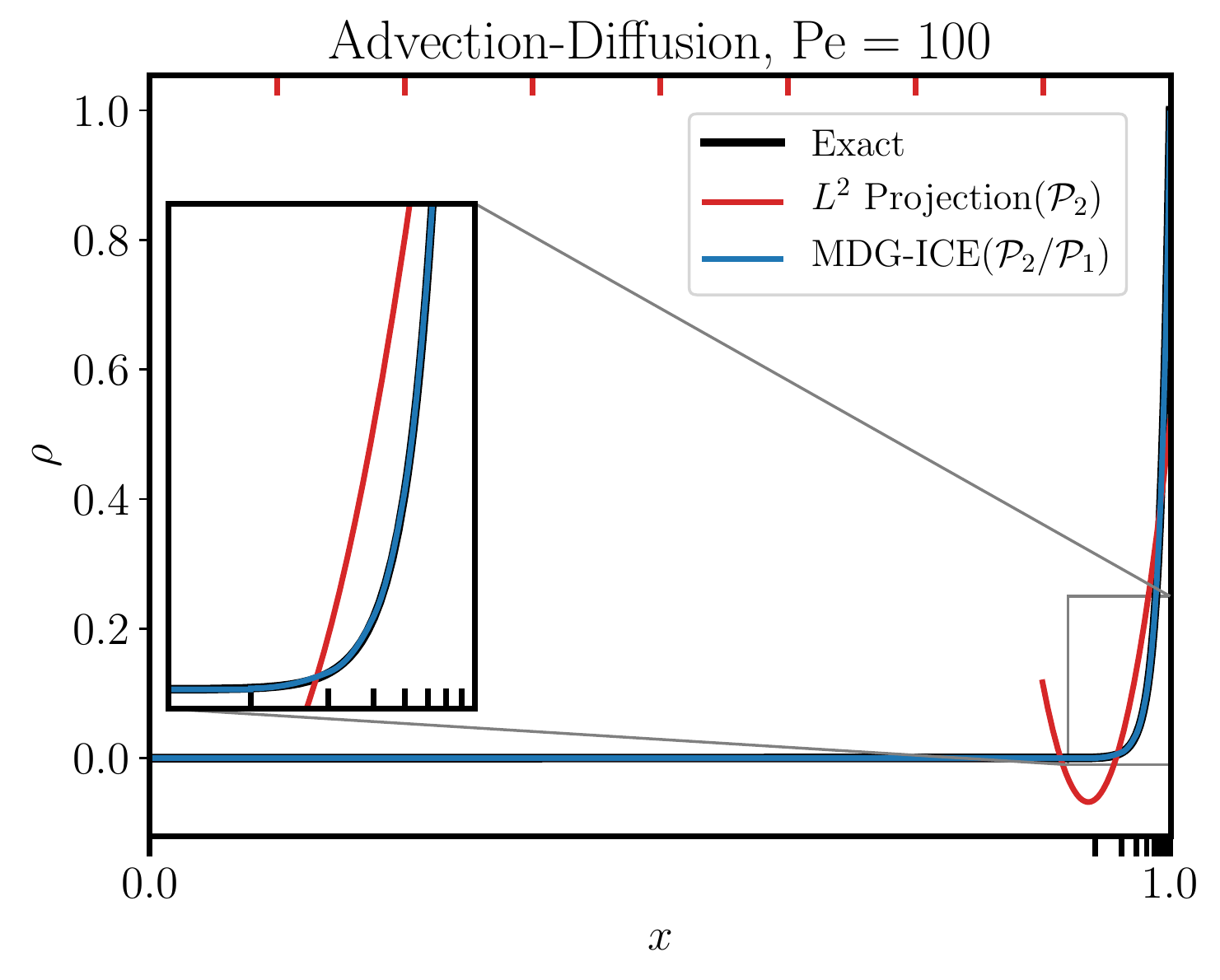}\tabularnewline
\includegraphics[width=0.4\linewidth]{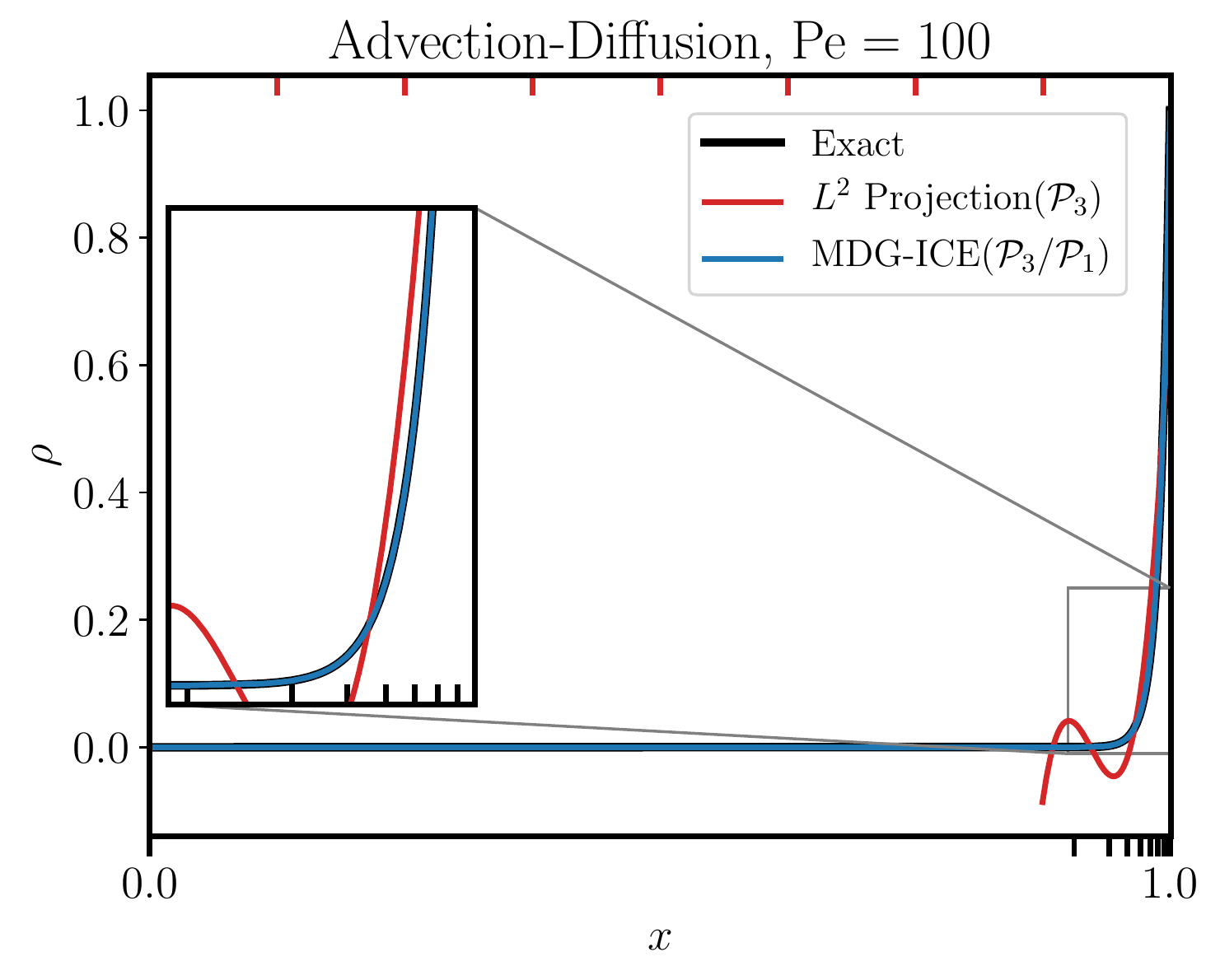} & \includegraphics[width=0.4\linewidth]{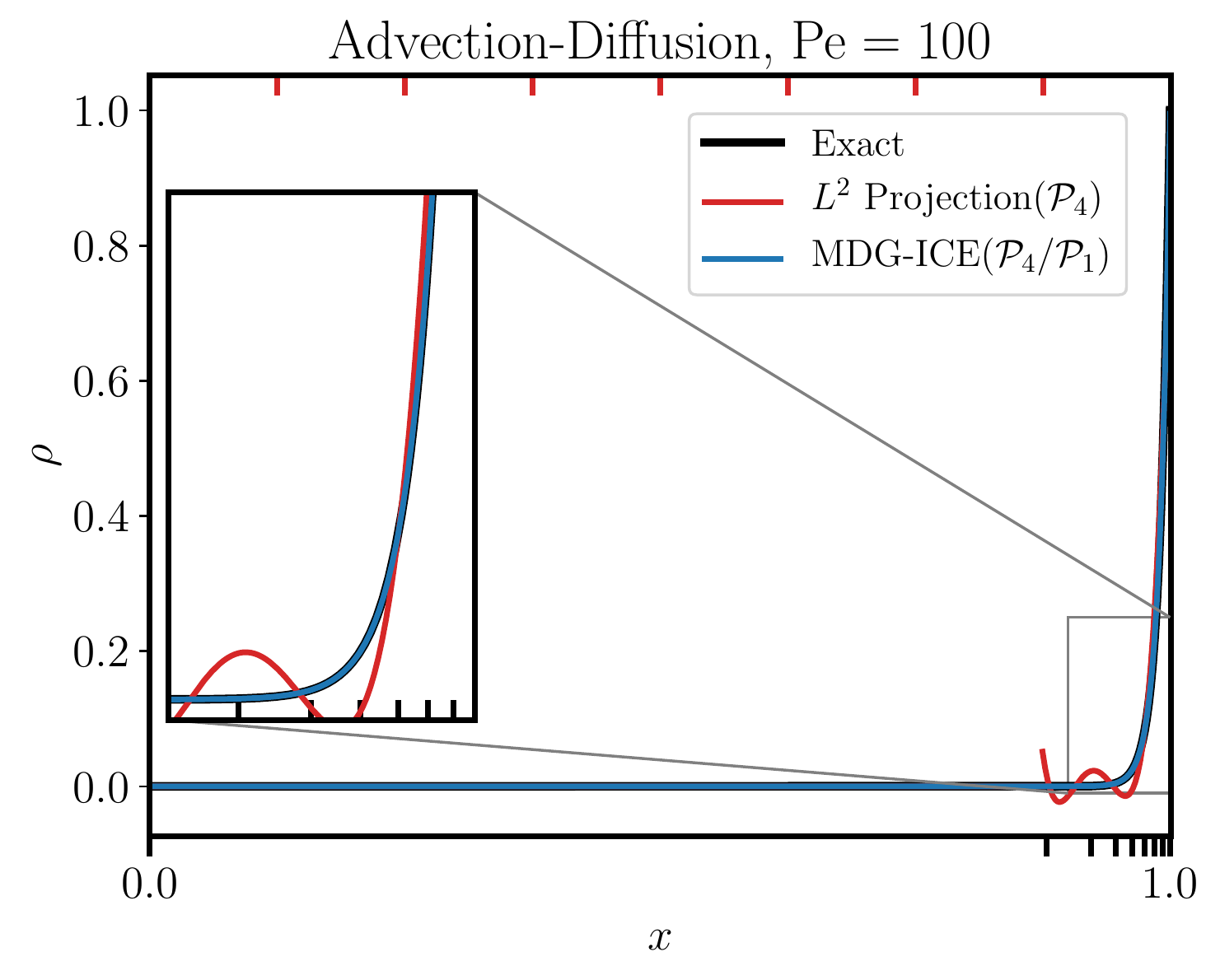}\tabularnewline
\end{tabular}\caption{\label{fig:steady-boundary-layer}Steady, linear advection-diffusion,
$\mathrm{Pe}=100$. The $L^{2}$ projection of the exact solution
onto a uniform grid consisting of 8 linear line cells is compared
for various polynomial degrees to MDG-ICE, which automatically moved
the initially uniform grid, indicated by red tick marks, in order
to resolve the boundary layer profile, resulting in the adapted grid
indicated with black tick marks.}
\end{figure}
\begin{figure}
\centering{}\includegraphics[width=0.6\linewidth]{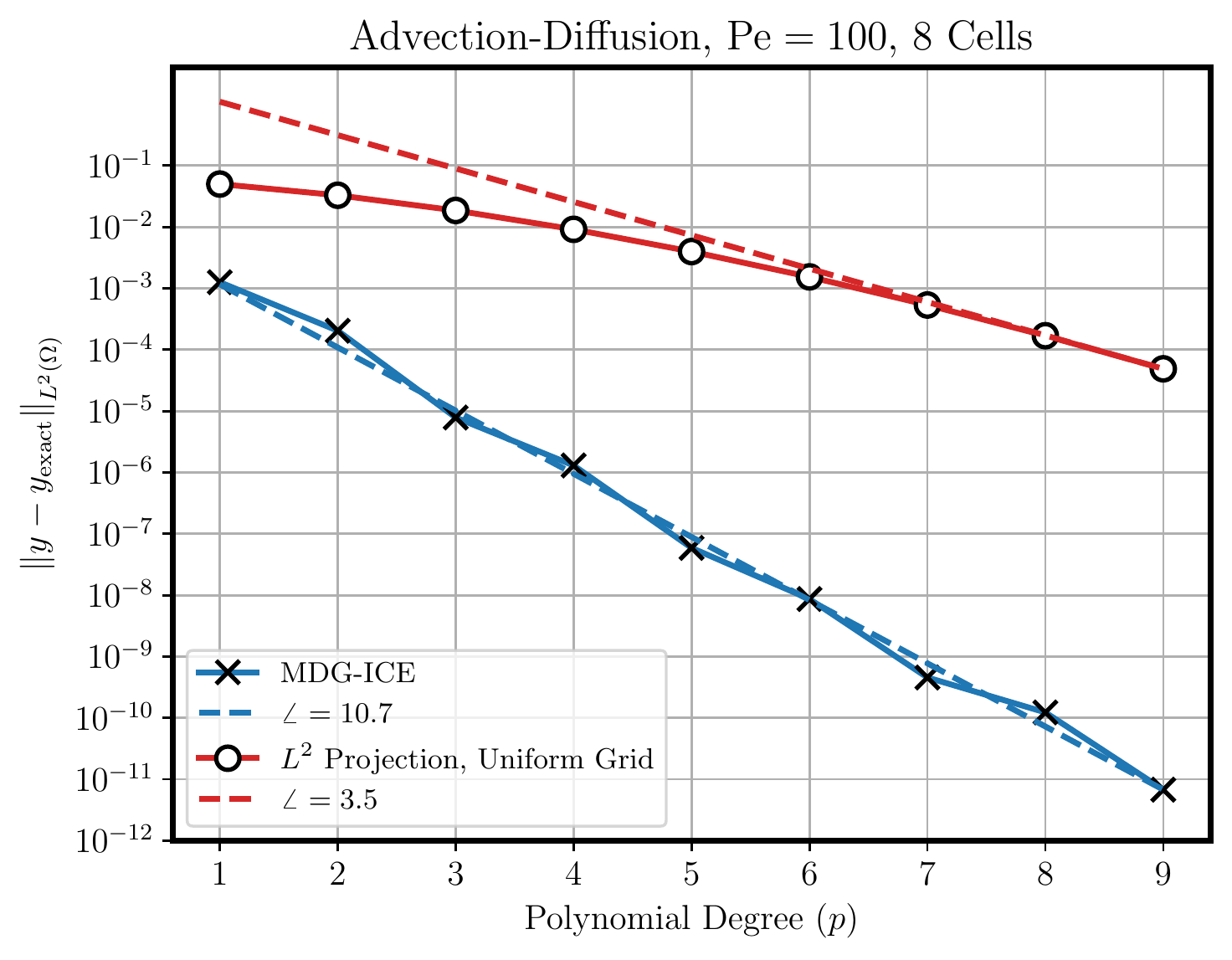}\caption{\label{fig:mdg-ice-steady_boundary_layer_pe_0010_p_refine}Steady,
linear advection-diffusion, $\mathrm{Pe}=100$. The rate of convergence
with respect to polynomial degree on a log-linear plot is shown, comparing
the $L^{2}$ projection onto a uniform grid to MDG-ICE, which automatically
moved the initially uniform grid to resolve the boundary layer profile.
Reference slopes of $10.7$ and $3.5$ are shown, illustrating the
increased rate of convergence achieved using MDG-ICE.}
\end{figure}

Figure~\ref{fig:steady-boundary-layer} shows the MDG-ICE solution
to the linear advection-diffusion problem with exact solution~(\ref{eq:steady-boundary-layer-exact-solution})
for $\mathrm{Pe}=100$ as well as the corresponding $L^{2}$ projection
of the exact solution onto a uniform grid of 8 linear line cells.
The $L^{2}$ projection minimizes the error in the $L^{2}$ norm and
therefore provides an upper bound on the accuracy attainable by methods
based on a static grid, e.g., DG. By moving the grid to resolve the
boundary layer profile, MDG-ICE is able to achieve accurate, oscillation-free,
solutions for a range of polynomial degrees.

Figure~\ref{fig:mdg-ice-steady_boundary_layer_pe_0010_p_refine}
presents the corresponding convergence results with respect to polynomial
degree, i.e., $p$-refinement. The rate of convergence of MDG-ICE
with respect to polynomial degree is compared to the $L^{2}$ projection
of the exact solution onto a uniform grid. These results confirm that
MDG-ICE resolves sharp boundary layers with enhanced accuracy compared
to static grid methods. Even for a $\mathcal{P}_{1}$ approximation,
MDG-ICE provides nearly two orders of magnitude improved accuracy
compared to the best approximation available on a uniform grid, a
gap that only widens at higher polynomial degrees. The MDG-ICE error
is plotted on a log-linear plot, with a reference slope of $10.7$,
indicating spectral convergence. This shows that the $r$-adaptivity
provided by MDG-ICE enhances the effectiveness of $p$-refinement,
even in the presence of initially under-resolved flow features.  This
results demonstrates the enhanced accuracy of MDG-ICE for resolving
initially under-resolved flow features using high-order finite element
approximation in comparison to traditional static grid methods, such
as DG.

\subsection{Space-time Burgers viscous shock formation \label{subsec:space-time-burgers-viscous-shock-formation}}

For a space-time Burgers flow, described in Section~\ref{subsec:One-dimensional-Burgers-flow},
a shock will form at time $t=t_{s}=0.5$ for the following initial
conditions

\begin{equation}
y\left(x,t=0\right)=\frac{1}{2\pi t_{s}}\sin\left(2\pi x\right)+y_{\infty},\label{eq:burgers-shock-formation-initial}
\end{equation}
where $y_{\infty}=0.2$ is the freestream velocity. The space-time
solution was initialized by extruding the temporal inflow condition,
given by Equation~(\ref{eq:burgers-shock-formation-initial}), throughout
the space-time domain.
\begin{figure}
\begin{centering}
\subfloat[\label{fig:burgers-shock-formation-spacetime-inviscid}Inviscid MDG-ICE$\left(\mathcal{P}_{5}/\mathcal{P}_{1}\right)$
space-time solution computed using 200 triangle elements.]{\begin{centering}
\includegraphics[width=0.45\linewidth]{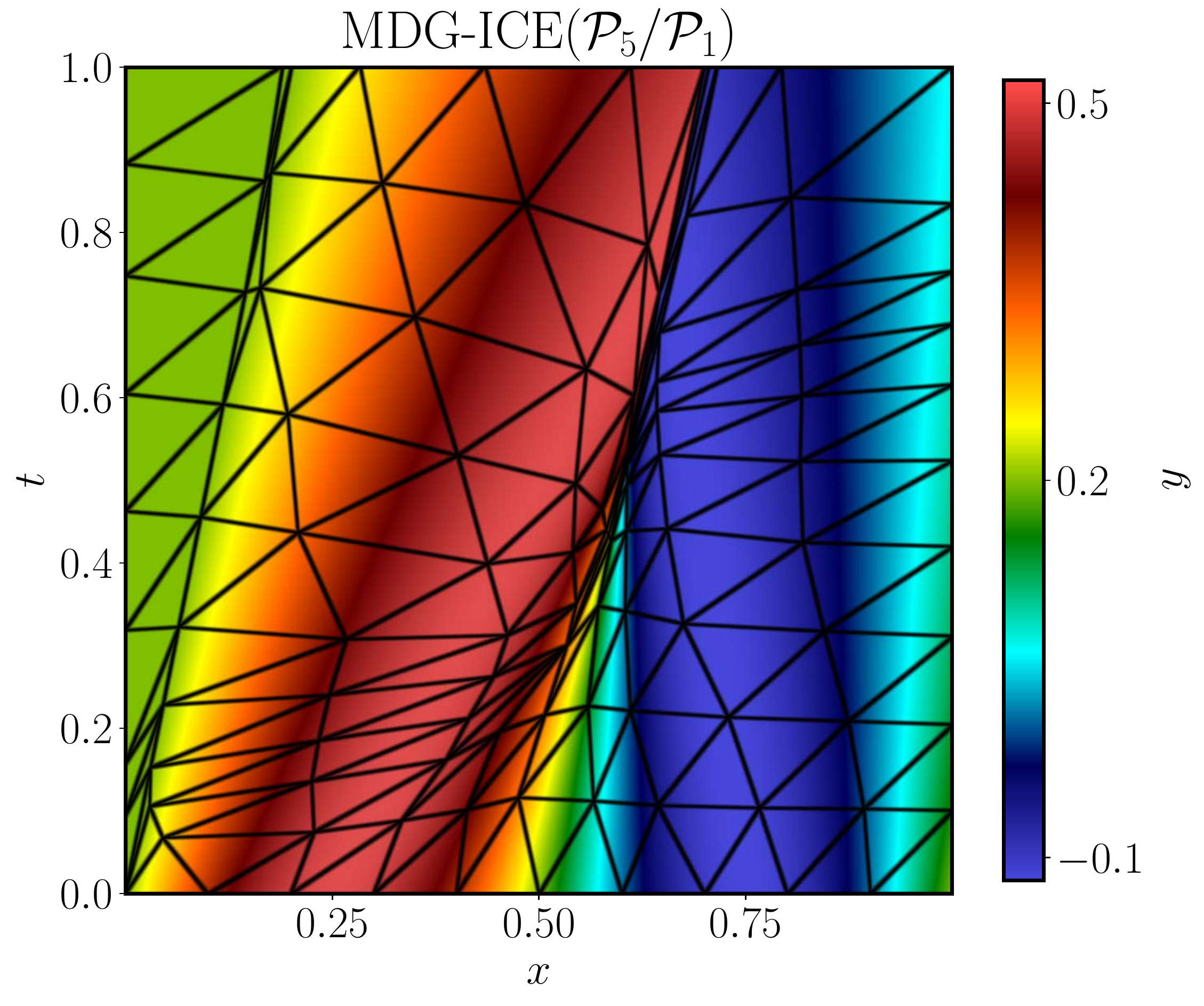}
\par\end{centering}
}
\par\end{centering}
\begin{centering}
\subfloat[\label{fig:burgers-shock-formation-spacetime-1e-3}MDG-ICE$\left(\mathcal{P}_{5}/\mathcal{P}_{1}\right)$
space-time solution for $\epsilon=10^{-3}$ computed using 200 triangle
elements.]{\begin{centering}
\includegraphics[width=0.45\linewidth]{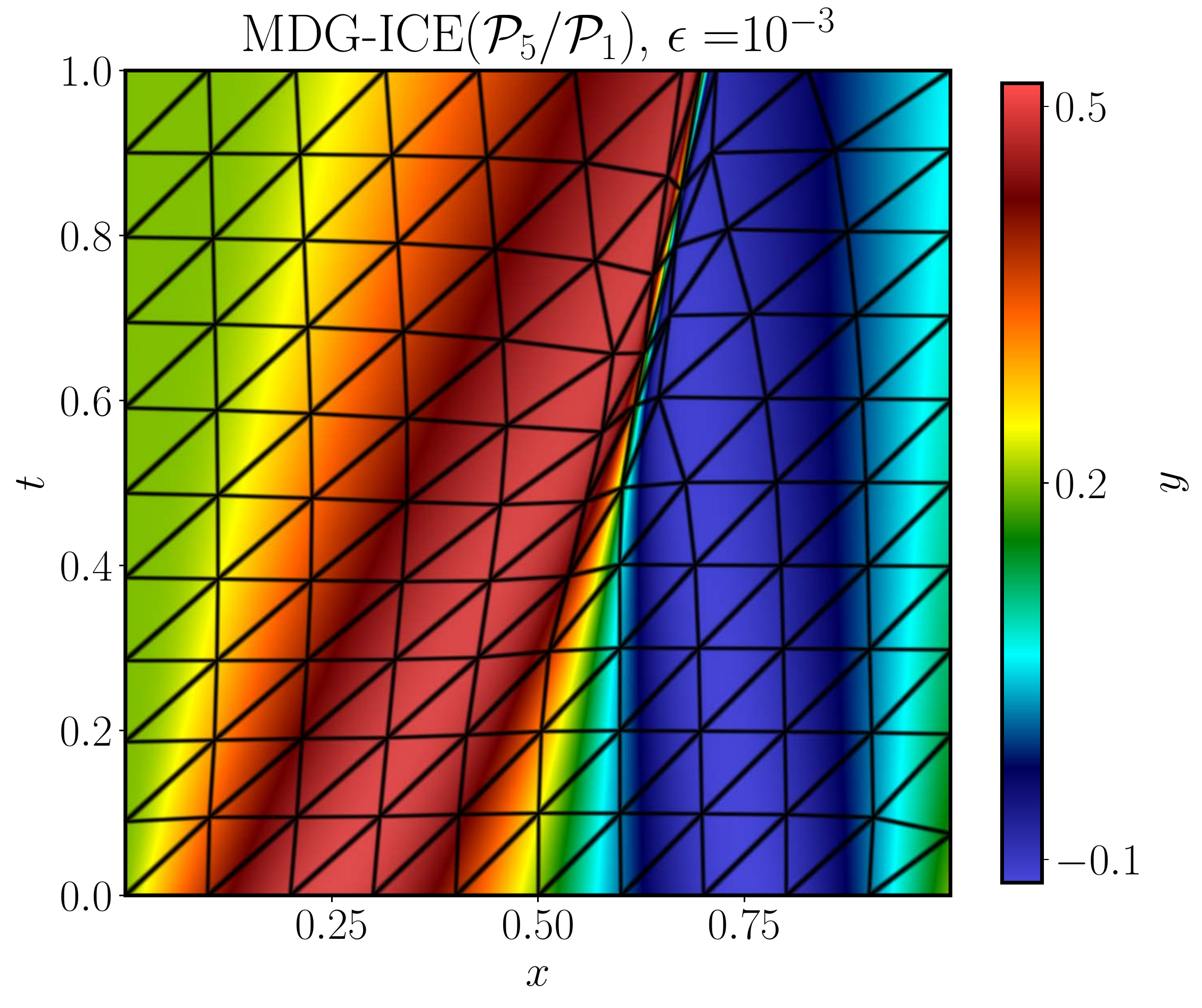}
\par\end{centering}
}
\par\end{centering}
\centering{}\subfloat[\label{fig:burgers-shock-formation-spacetime-1e-4}MDG-ICE$\left(\mathcal{P}_{5}/\mathcal{P}_{1}\right)$
space-time solution for $\epsilon=10^{-4}$ computed using 200 triangle
elements.]{\begin{centering}
\includegraphics[width=0.45\linewidth]{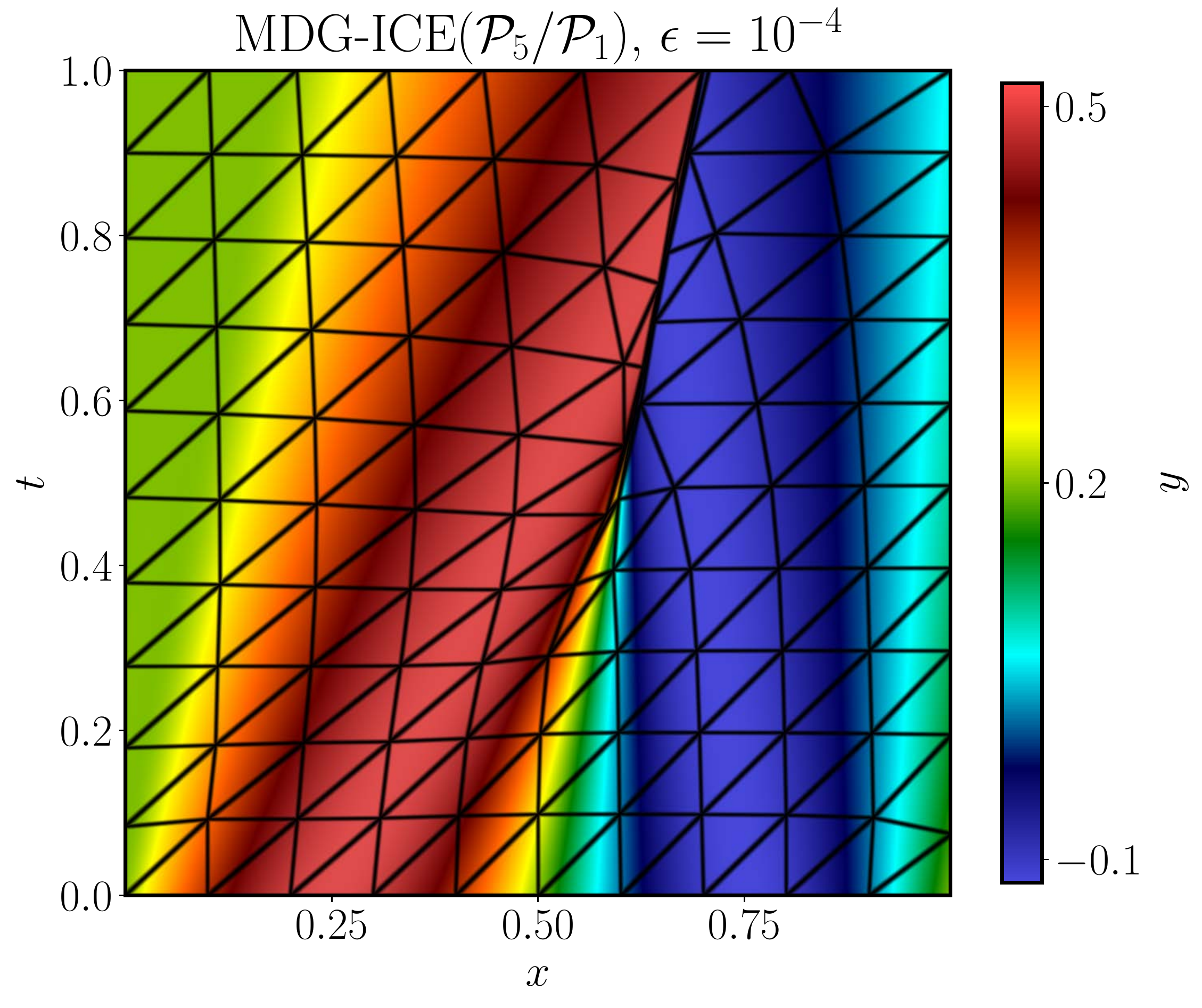}
\par\end{centering}
}\caption{\label{fig:burgers-shock-formation-spacetime}Space-time Burgers shock
formation: inviscid, viscous $\left(\epsilon=10^{-3}\right)$, and
viscous $\left(\epsilon=10^{-4}\right)$ solutions computed using
$\mathcal{P}_{5}$ linear triangle elements without shock capturing.
Instead, the viscous shock was resolved via anisotropic space-time
$r$-adaptivity. The solver was initialized by extruding the initial
condition at $t=0$ in time.}
\end{figure}
\begin{figure}
\centering{}%
\begin{minipage}[c][1\totalheight][t]{0.3\columnwidth}%
\subfloat[\label{fig:burgers-shock-formation-inviscid}Inviscid MDG-ICE$\left(\mathcal{P}_{5}/\mathcal{P}_{1}\right)$
at $t=0$ and $t=1$. The corresponding space-time solution is shown
in Figure~\ref{fig:burgers-shock-formation-spacetime-inviscid}.]{\begin{centering}
\includegraphics[width=0.8\linewidth]{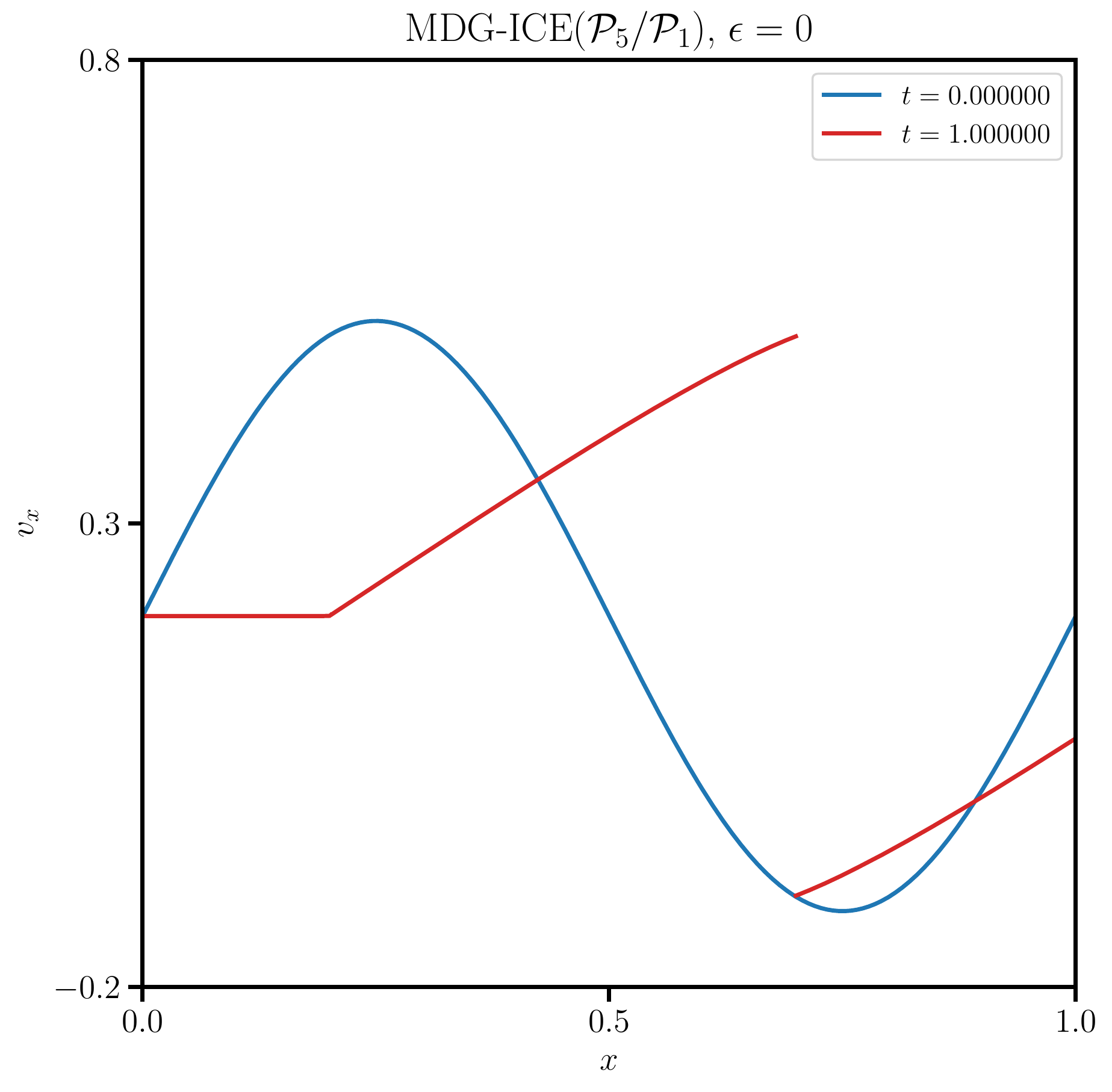}
\par\end{centering}
}%
\end{minipage}%
\begin{minipage}[c][1\totalheight][t]{0.3\columnwidth}%
\subfloat[\label{fig:burgers-shock-formation-1e-3}Viscous MDG-ICE$\left(\mathcal{P}_{5}/\mathcal{P}_{1}\right)$
with $\epsilon=\protect\expnumber 1{-3}$ at $t=0$ and $t=1$ . The
corresponding space-time solution is shown in Figure~\ref{fig:burgers-shock-formation-spacetime-1e-3}.]{\begin{centering}
\includegraphics[width=0.8\linewidth]{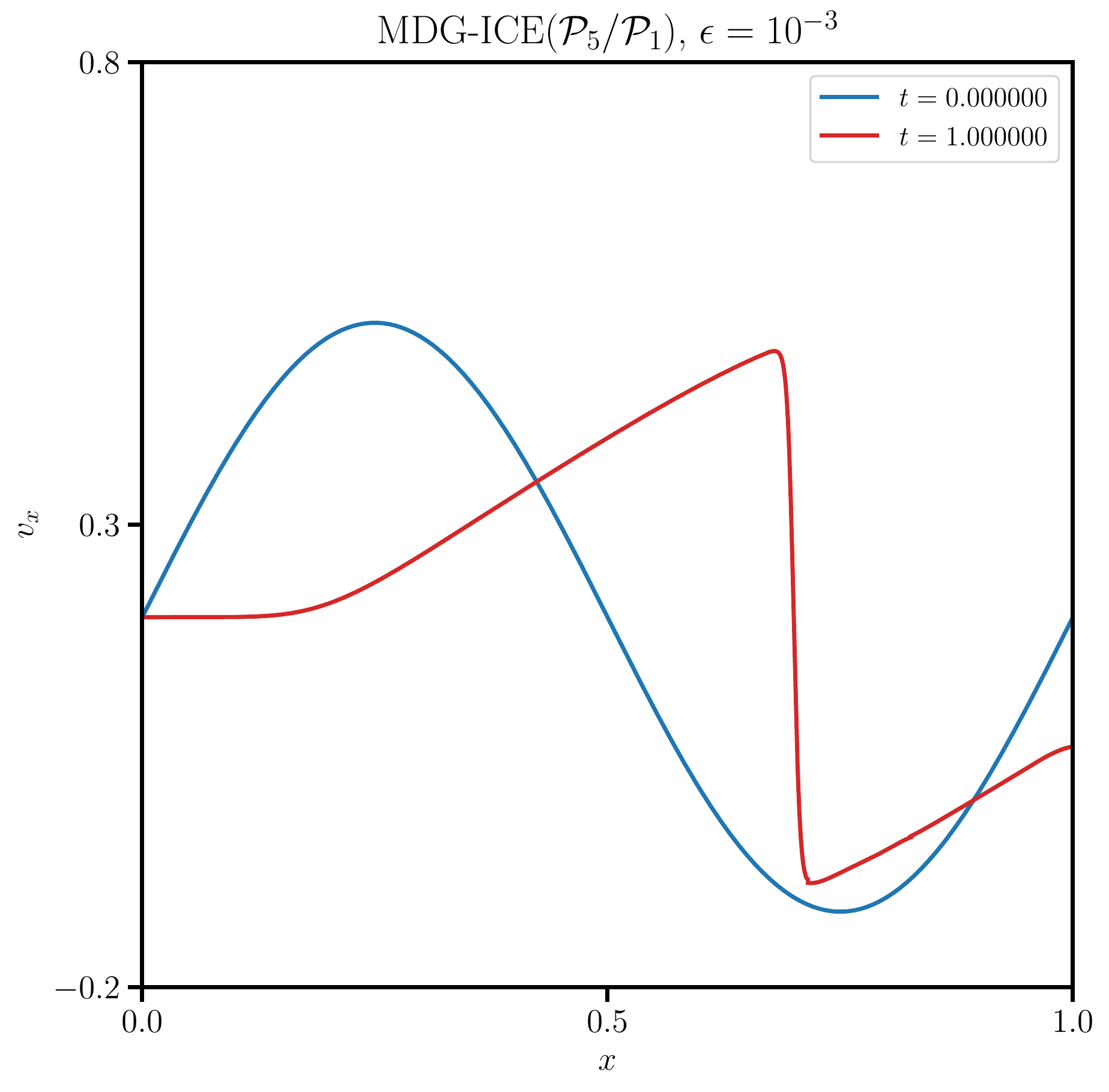}
\par\end{centering}
}%
\end{minipage}%
\begin{minipage}[c][1\totalheight][t]{0.3\columnwidth}%
\subfloat[\label{fig:burgers-shock-formation-1e-4}Viscous MDG-ICE$\left(\mathcal{P}_{5}/\mathcal{P}_{1}\right)$
with $\epsilon=\protect\expnumber 1{-4}$ at $t=0$ and $t=1$ . The
corresponding space-time solution is shown in Figure~\ref{fig:burgers-shock-formation-spacetime-1e-4}.]{\begin{centering}
\includegraphics[width=0.8\linewidth]{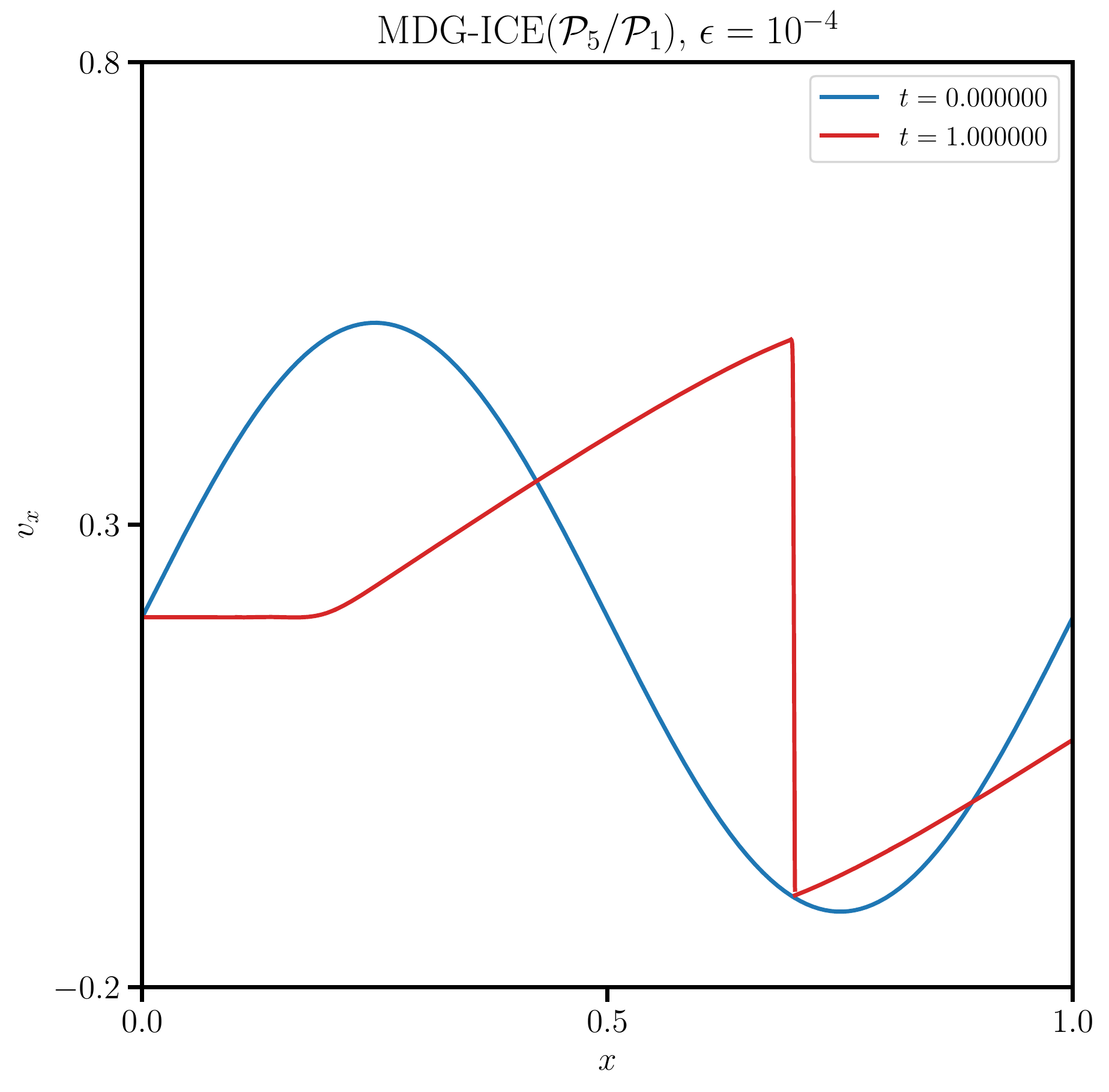}
\par\end{centering}
}%
\end{minipage}\caption{\label{fig:burgers-shock-formation}Burgers shock formation one-dimensional
profiles at $t=0$ and $t=1$: inviscid, viscous $\left(\epsilon=10^{-3}\right)$,
and viscous $\left(\epsilon=10^{-4}\right)$ solutions computed using
$\mathcal{P}_{5}$ linear triangle elements without shock capturing.
The corresponding space-time solutions are shown in Figure~\ref{fig:burgers-shock-formation-spacetime}.}
\end{figure}

Figure~\ref{fig:burgers-shock-formation-spacetime} presents the
space-time Burgers shock formation solutions for an inviscid flow,
a viscous flow with $\epsilon=10^{-3}$, and a viscous flow with $\epsilon=10^{-4}$.
Figure~\ref{fig:burgers-shock-formation} presents the corresponding
one-dimensional profiles at $t=0$ and $t=1$. The space-time solutions
were initialized by extruding the inflow condition, given by Equation~(\ref{eq:burgers-shock-formation-initial}),
at $t=0$ in time. The initial simplicial grid was generated by converting
a uniform $10\times10$ quadrilateral grid into triangles.

In the inviscid case, MDG-ICE fits the point of shock formation and
tracks the shock at the correct speed. In addition to the shock, the
inviscid flow solution has a derivative discontinuity that MDG-ICE
also detects and tracks at the correct speed of $0.2$. For the two
viscous flow cases, $\epsilon=10^{-3}$ and $\epsilon=10^{-4}$, MDG-ICE
accurately resolves each viscous shock as a sharp, yet smooth, profile
by adjusting the grid geometry, without modifying the grid topology.
This case demonstrates the inherent ability of MDG-ICE to achieve
anisotropic space-time $r$-adaptivity for unsteady flow problems.

\subsection{Mach 5 viscous bow shock\label{subsec:Viscous-Bow-Shock}}

\begin{figure}
\begin{centering}
\subfloat[\label{fig:Viscous-Bow-Shock-DG-Mesh}392 linear triangle cells]{\begin{centering}
\includegraphics[width=0.9\linewidth]{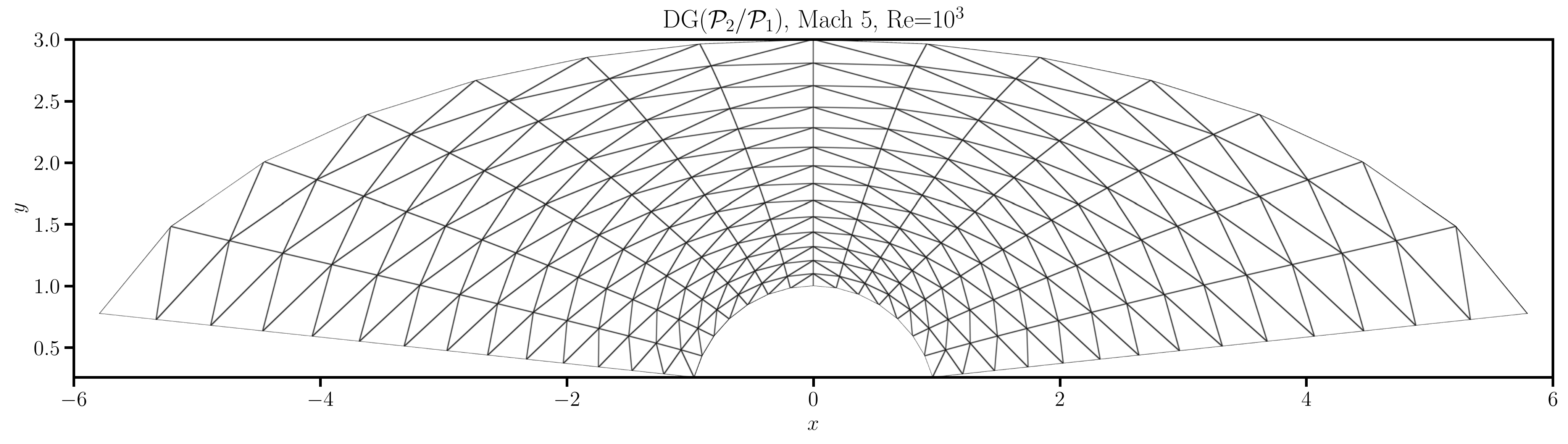}
\par\end{centering}
}
\par\end{centering}
\begin{centering}
\subfloat[\label{fig:Viscous-Bow-Shock-DG-Temperature}392 linear $\mathcal{P}_{2}$
triangle cells]{\begin{centering}
\includegraphics[width=0.9\linewidth]{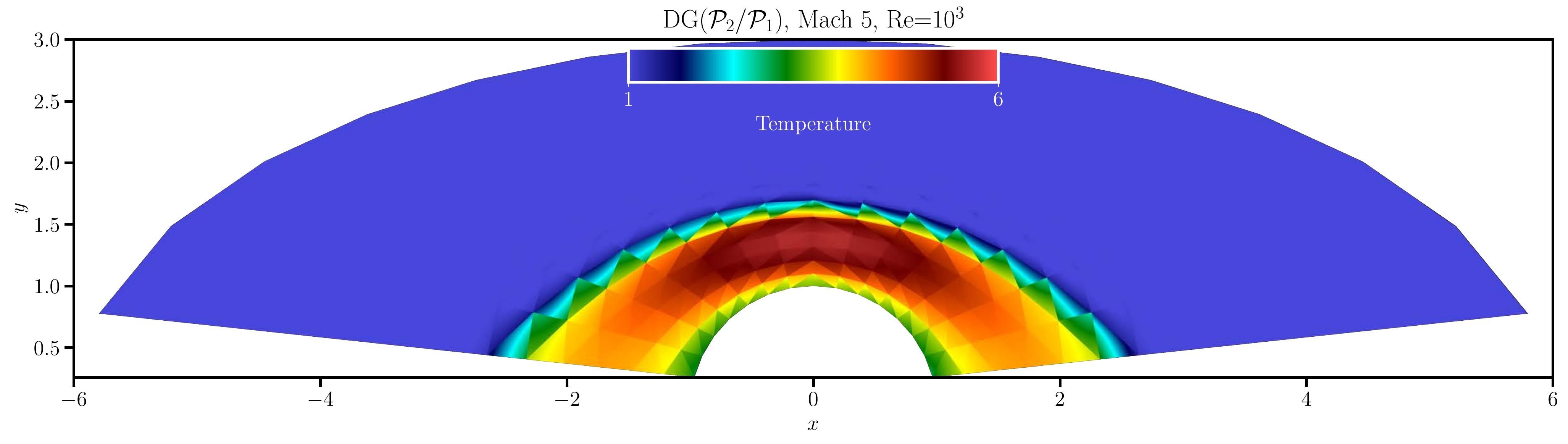}
\par\end{centering}
}
\par\end{centering}
\caption{\label{fig:Viscous-Bow-Shock-DG}The initial linear grid and temperature
field corresponding to a shock captured DG($\mathcal{P}_{2}/\mathcal{P}_{1}$)
solution for the viscous Mach 5 bow shock at $\mathrm{Re}=10^{3}$.}
\end{figure}
\begin{figure}
\begin{centering}
\subfloat[\label{fig:Viscous-Bow-Shock-Mesh}400 isoparametric $\mathcal{P}_{4}$
triangle cells]{\begin{centering}
\includegraphics[width=0.9\linewidth]{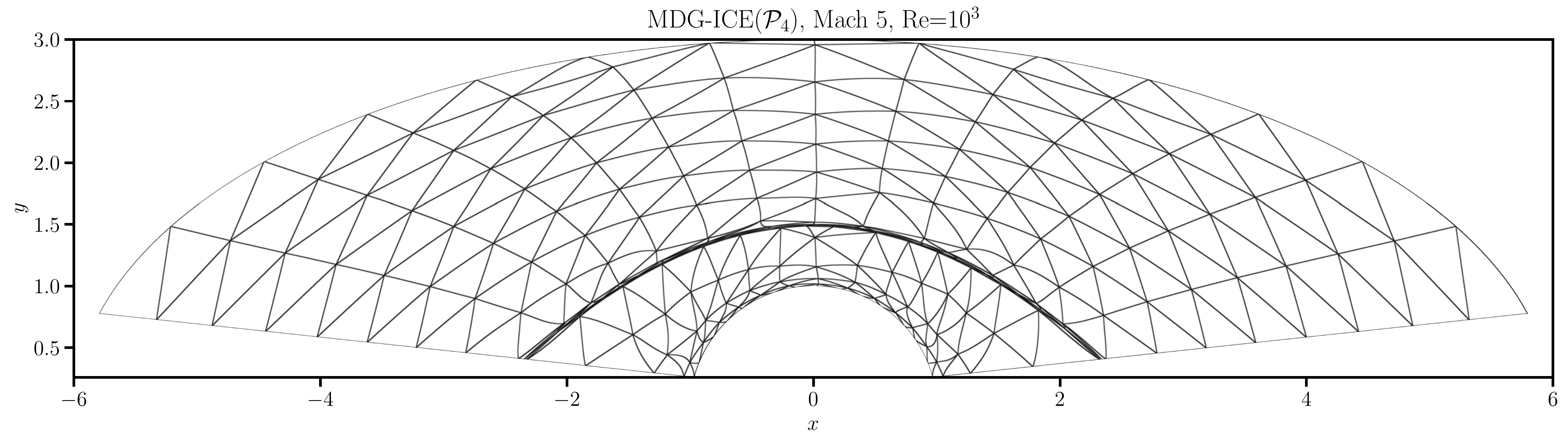}
\par\end{centering}
}
\par\end{centering}
\begin{centering}
\subfloat[\label{fig:Viscous-Bow-Shock-Temperature}400 isoparametric $\mathcal{P}_{4}$
triangle cells]{\begin{centering}
\includegraphics[width=0.9\linewidth]{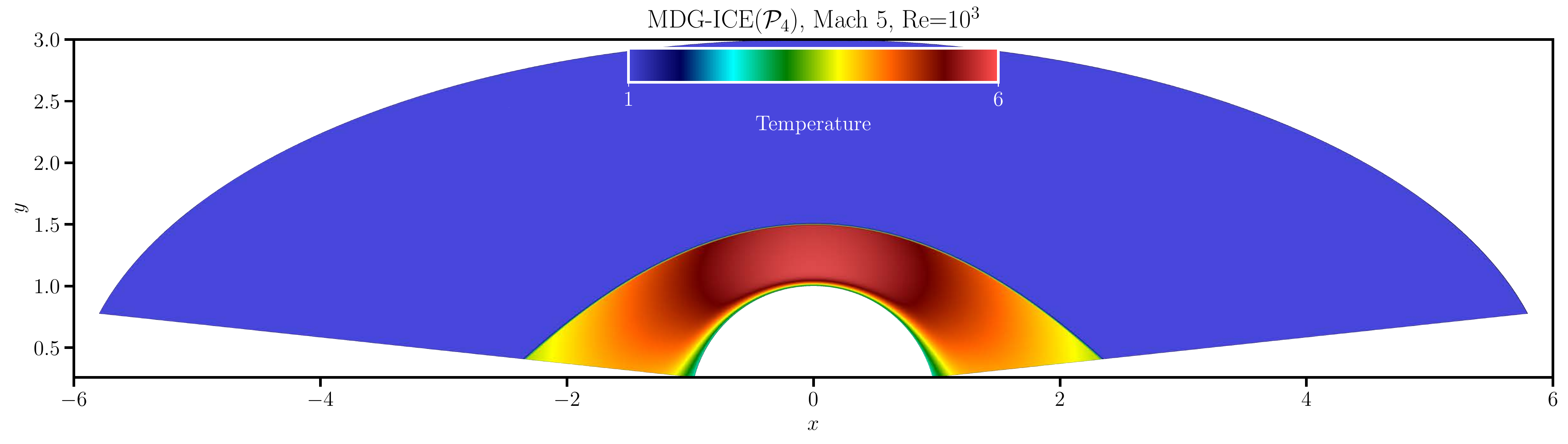}
\par\end{centering}
}
\par\end{centering}
\begin{centering}
\subfloat[\label{fig:Viscous-Bow-Shock-Mach}400 isoparametric $\mathcal{P}_{4}$
triangle cells]{\begin{centering}
\includegraphics[width=0.9\linewidth]{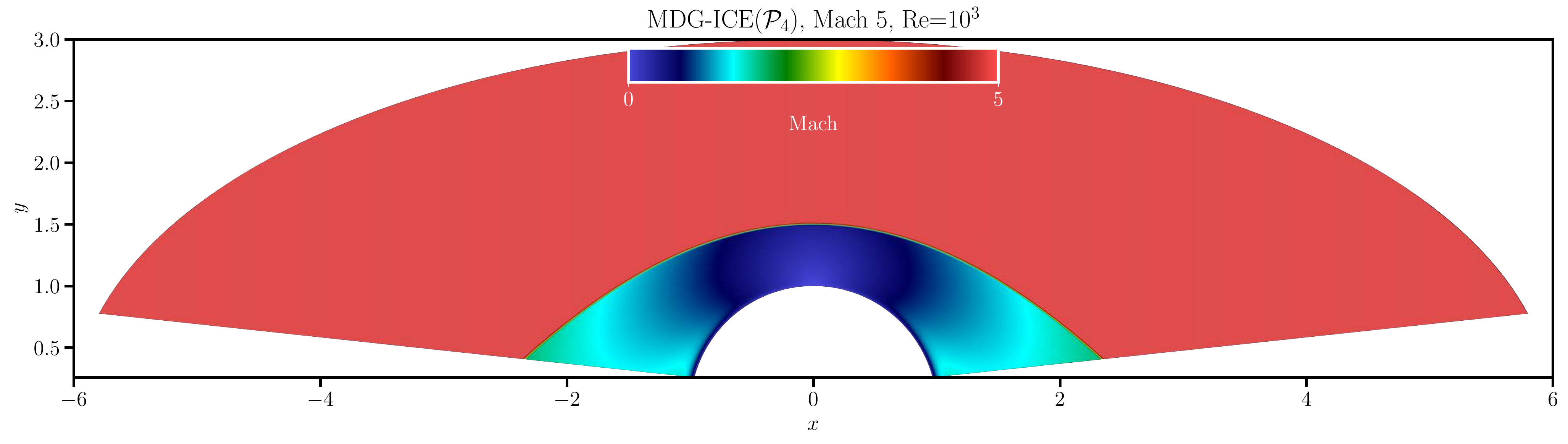}
\par\end{centering}
}
\par\end{centering}
\centering{}\subfloat[\label{fig:Viscous-Bow-Shock-Pressure}400 isoparametric $\mathcal{P}_{4}$
triangle cells]{\begin{centering}
\includegraphics[width=0.9\linewidth]{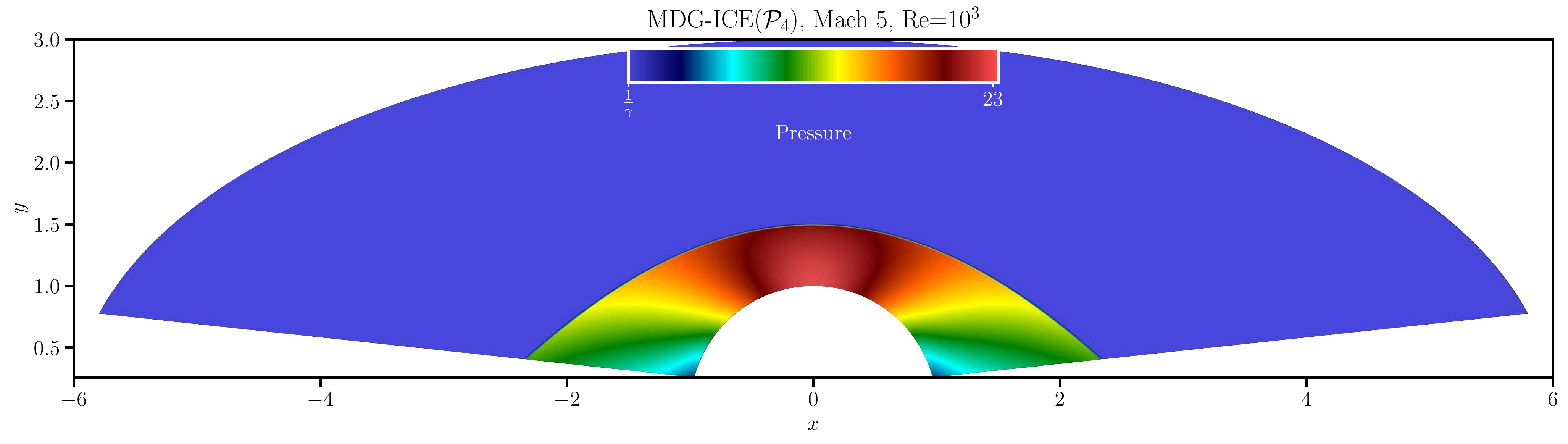}
\par\end{centering}
}\caption{\label{fig:Viscous-Bow-Shock-MDG-ICE}The MDG-ICE solution computed
using $\mathcal{P}_{4}$ isoparametric triangle elements for the viscous
Mach 5 bow shock at $\mathrm{Re}=10^{3}$. The MDG-ICE grid was initialized
by projecting the linear triangle grid shown in Figure~\ref{fig:Viscous-Bow-Shock-DG-Mesh}
to the closest point on the boundary of the domain. The MDG-ICE field
variables were initialized by cell averaging the interpolated the
DG($\mathcal{P}_{2}/\mathcal{P}_{1}$) solution shown in Figure~\ref{fig:Viscous-Bow-Shock-DG-Temperature}.
The MDG-ICE flux variables were initialized to zero for consistency
with the initial piecewise constant field variables. The location
of the shock along the line $x=0$ was computed as $y=1.49995$ for
a stand-off distance of $0.49995$.}
\end{figure}
\begin{figure}
\subfloat[\label{fig:Bow-Shock-1D-Temperature}The temperature sampled along
$x=0$. The exact temperature at the stagnation point, $T=2.5$, is
marked with the symbol $\times$. ]{\begin{centering}
\includegraphics[width=0.3\linewidth]{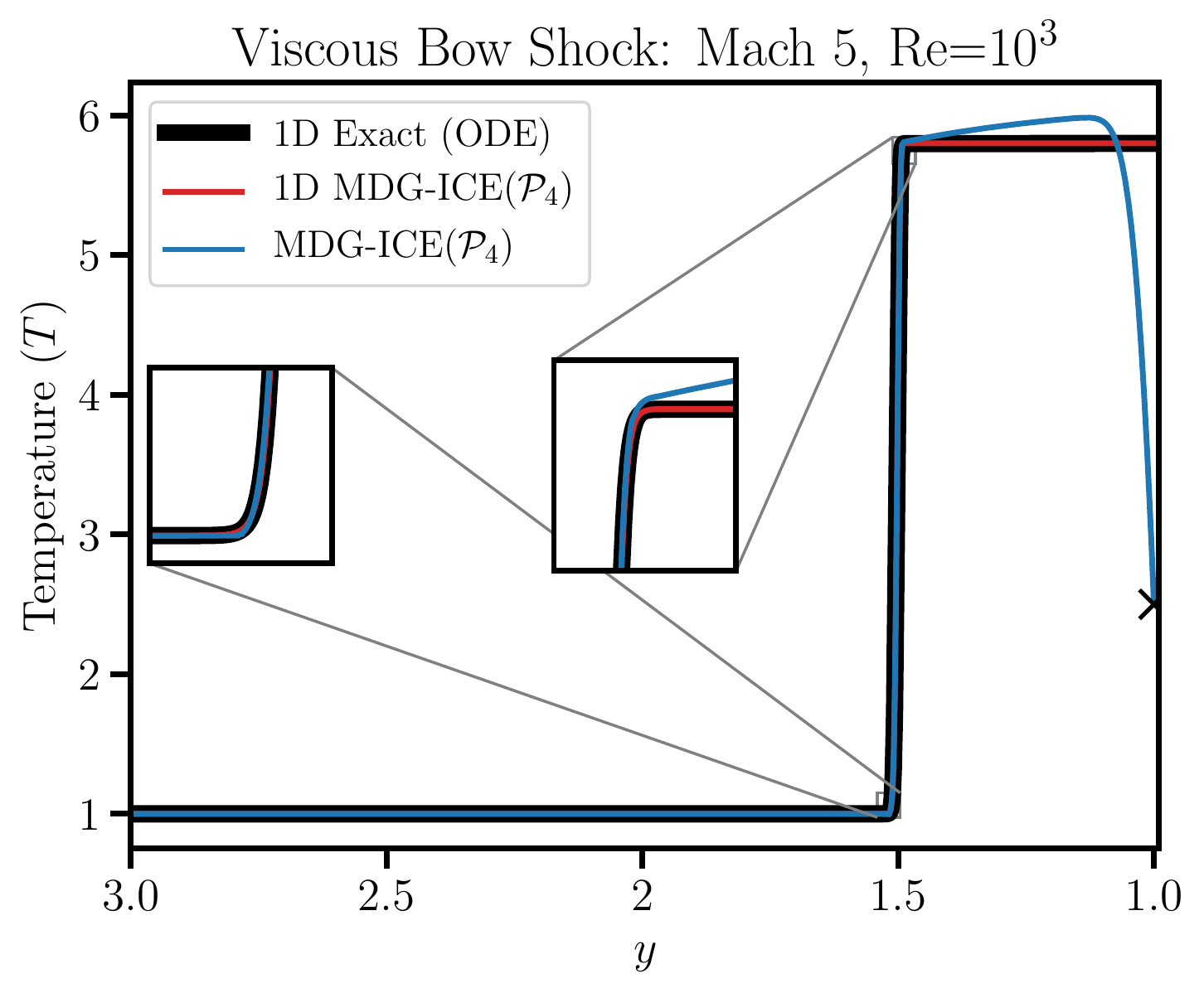}
\par\end{centering}
}\hfill{}\subfloat[\label{fig:Bow-Shock-1D-VelocityX}The normal velocity, $v_{n}=0$,
sampled along $x=0$. The exact normal velocity at the stagnation
point, $v_{n}=0$, is marked with the symbol $\times$. ]{\begin{centering}
\includegraphics[width=0.3\linewidth]{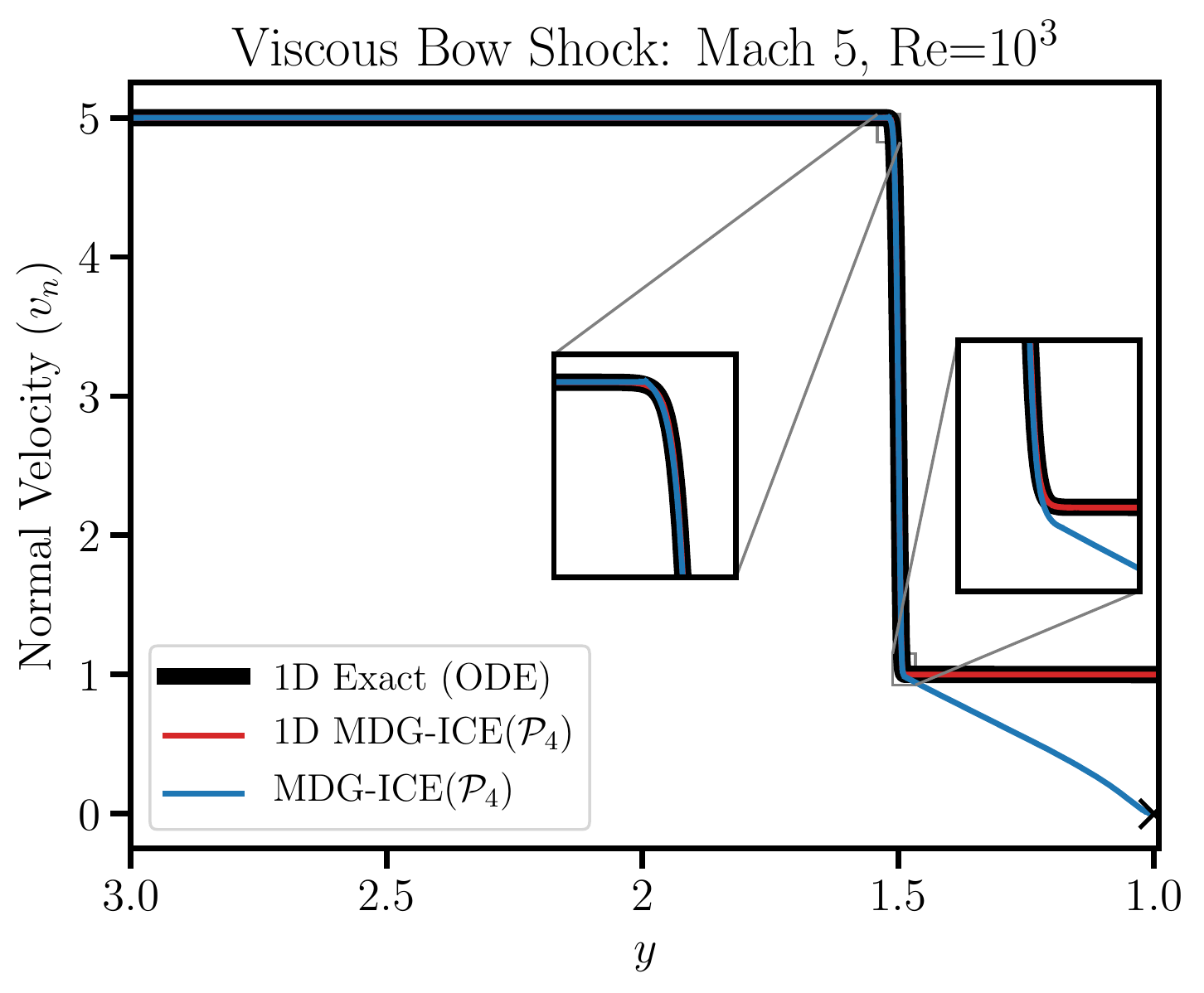}
\par\end{centering}
}\hfill{}\subfloat[\label{fig:Bow-Shock-1D-Pressure}The pressure, $p$, sampled along
$x=0$. The exact pressure at the stagnation point for an inviscid
flow, $p\approx23.324$, is marked with the symbol $\times$.]{\begin{centering}
\includegraphics[width=0.3\linewidth]{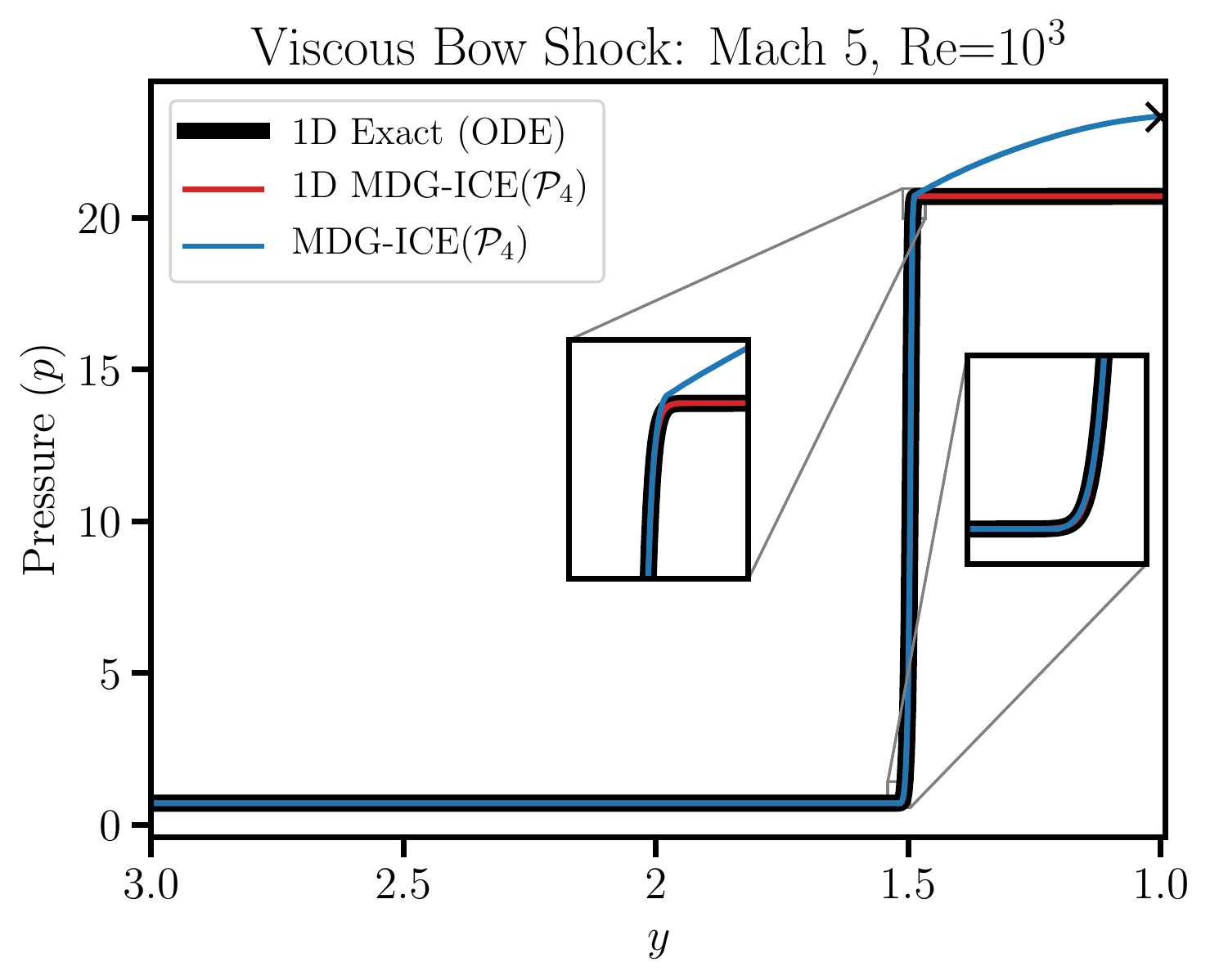}
\par\end{centering}
}

\subfloat[\label{fig:Bow-Shock-1D-Density}The density, $\rho$, sampled along
$x=0$. The density at the stagnation point, computed using the stagnation
pressure corresponding to an inviscid flow, $\rho\approx13.061389724919298$,
is marked with the symbol $\times$. ]{\begin{centering}
\includegraphics[width=0.3\linewidth]{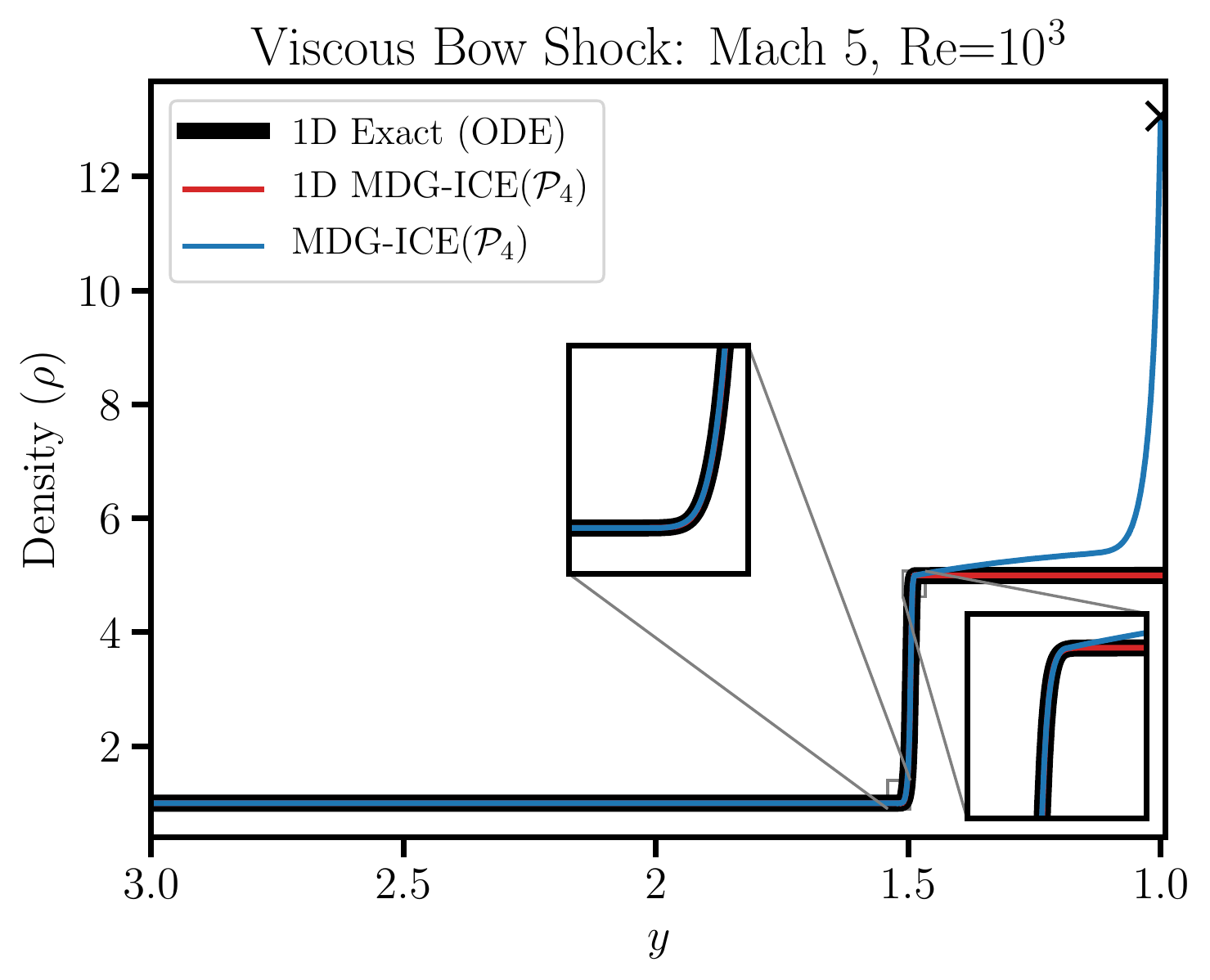}
\par\end{centering}
}\hfill{}\subfloat[\label{fig:Bow-Shock-1D-FluxXMomentumX}The normal component of the
normal viscous stress tensor, $\tau_{nn}$, sampled along $x=0$.
]{\begin{centering}
\includegraphics[width=0.3\linewidth]{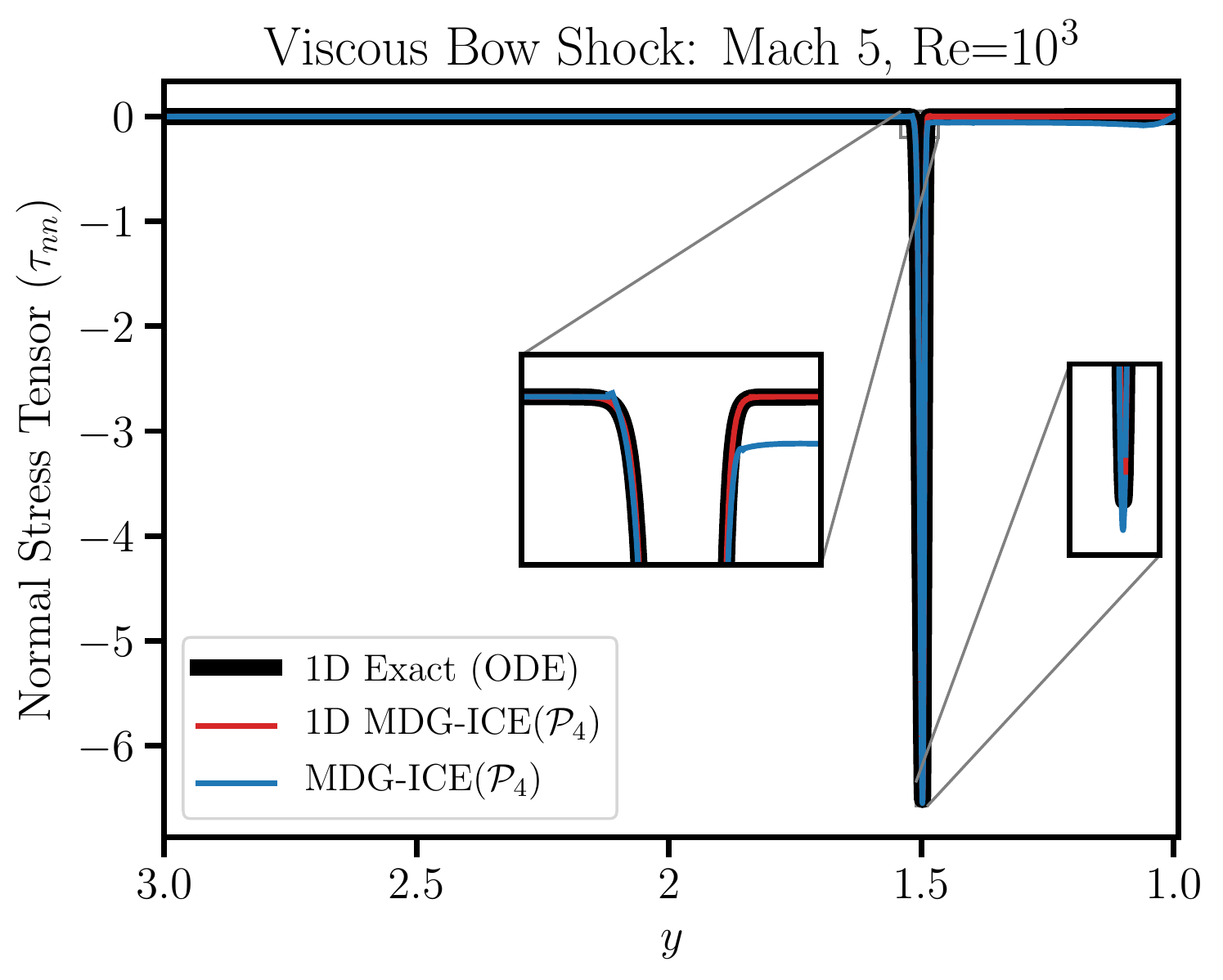}
\par\end{centering}
}\hfill{}\subfloat[\label{fig:Bow-Shock-1D-FluxXEnergyStagnationDensity}The normal thermal
heat flux, $q_{n}$, sampled along $x=0$. ]{\begin{centering}
\includegraphics[width=0.3\linewidth]{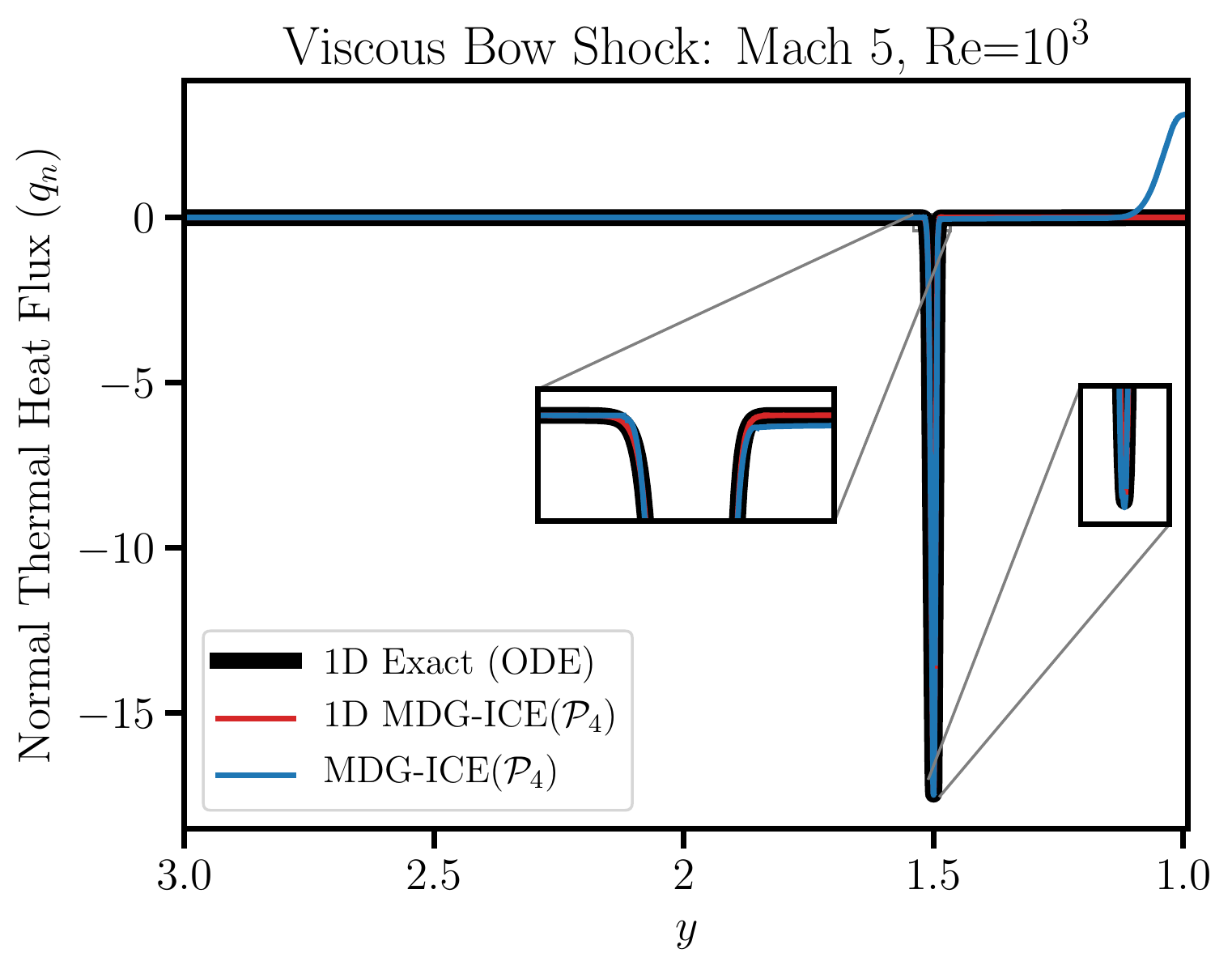}
\par\end{centering}
}

\caption{\label{fig:Bow-Shock-1D}Centerline profiles of temperature and normal
velocity for the viscous Mach 5 bow shock at $\mathrm{Re}=10^{3}$
computed with MDG-ICE($\mathcal{P}_{4}$) compared to ODE and MDG-ICE($\mathcal{P}_{4}$)
approximations of the exact solution for the corresponding one-dimensional
viscous shock. The one-dimensional MDG-ICE($\mathcal{P}_{4}$) approximation
was computed using 16 isoparametric line cells. The location of the
shock was computed as $y=1.49995$ for a stand-off distance of $0.49995$.}
\end{figure}
\begin{figure}
\subfloat[\label{fig:Bow-Shock-1D-PressureCoefficient}The pressure coefficient,
$C_{p}$, sampled at each degree of freedom on the surface of the
cylinder. The exact pressure coefficient at the stagnation point for
an inviscid flow, $C_{p}\approx1.8087699607027568$, is marked with
the symbol $\times$. ]{\begin{centering}
\includegraphics[width=0.45\linewidth]{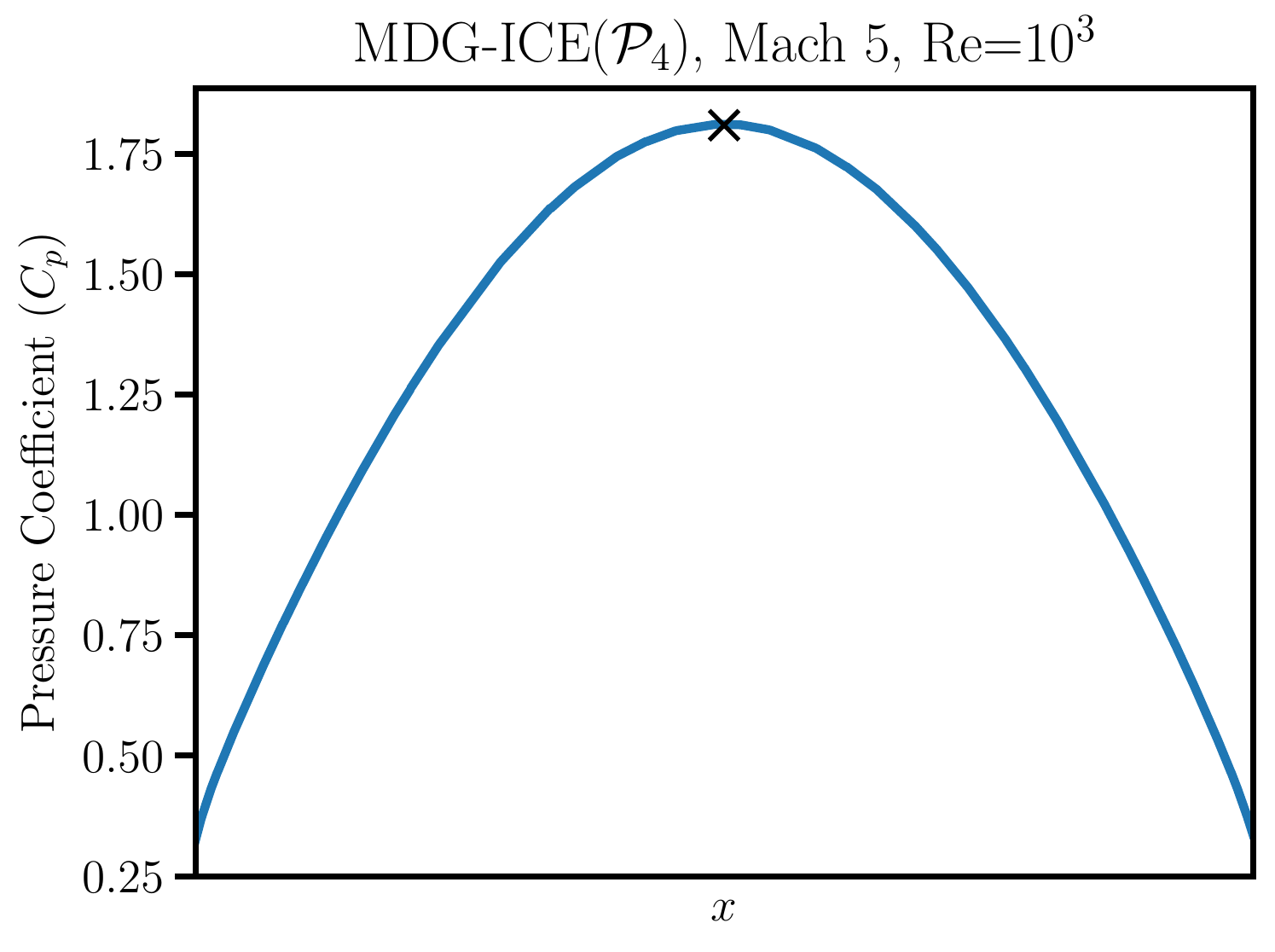}
\par\end{centering}
}\hfill{}\subfloat[\label{fig:Bow-Shock-1D-StantonNumber}The Stanton number sampled
at each degree of freedom on the surface of the cylinder. The computed
Stanton number at the stagnation point is marked with the symbol
$\times$.]{\begin{centering}
\includegraphics[width=0.45\linewidth]{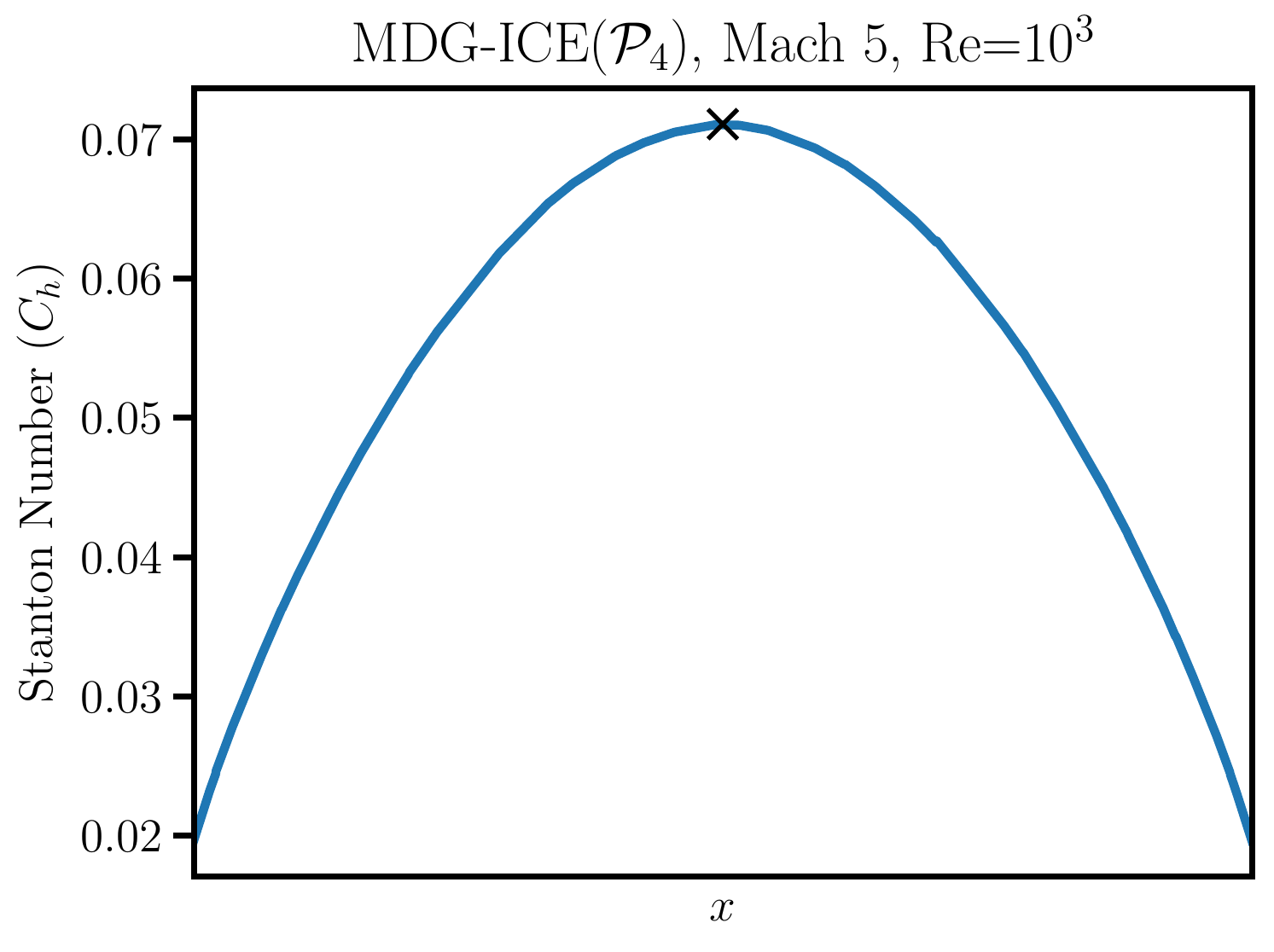}
\par\end{centering}
}

\caption{\label{fig:Bow-Shock-1D-PressureCoefficient-StantonNumber}Pressure
coefficient and Stanton number for the viscous Mach 5 bow shock at
$\mathrm{Re}=10^{3}$ computed with MDG-ICE($\mathcal{P}_{4}$) .}
\end{figure}
\begin{figure}
\begin{centering}
\subfloat[\label{fig:Viscous-Bow-Shock-Mesh-Re-010000}527 isoparametric $\mathcal{P}_{4}$
triangle cells]{\begin{centering}
\includegraphics[width=0.9\linewidth]{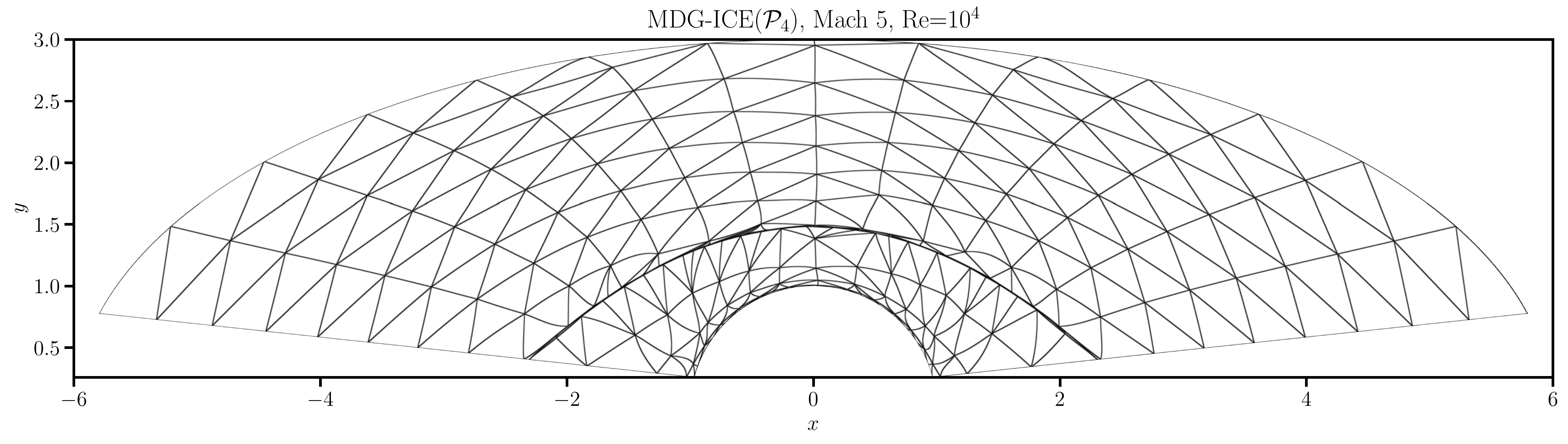}
\par\end{centering}
}
\par\end{centering}
\begin{centering}
\subfloat[\label{fig:Viscous-Bow-Shock-Temperature-Re-010000}527 isoparametric
$\mathcal{P}_{4}$ triangle cells]{\begin{centering}
\includegraphics[width=0.9\linewidth]{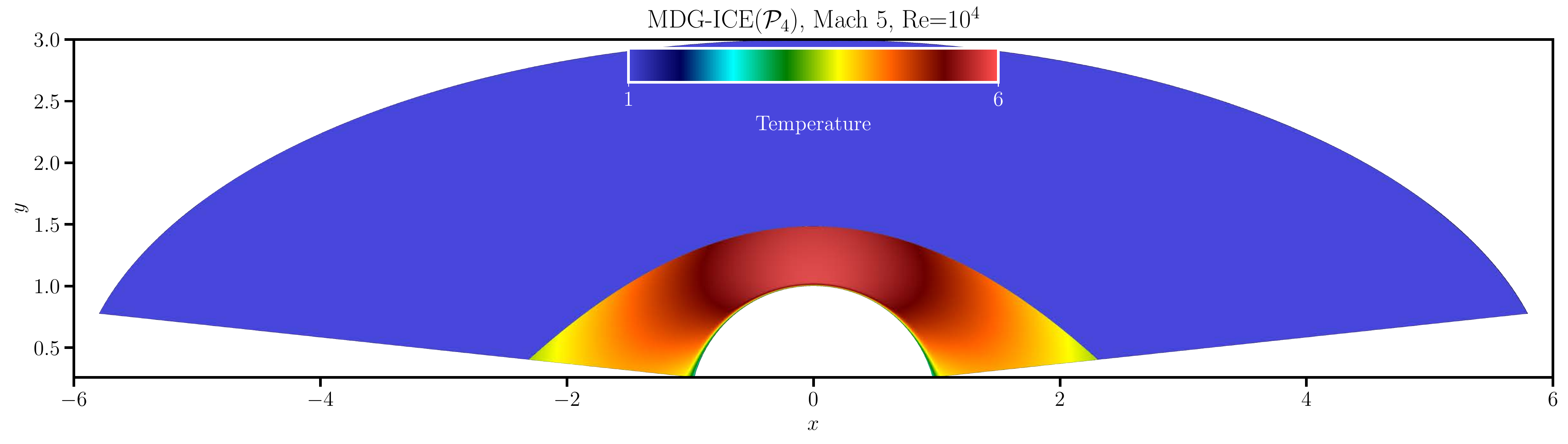}
\par\end{centering}
}
\par\end{centering}
\begin{centering}
\subfloat[\label{fig:Viscous-Bow-Shock-Mach-Re-010000}527 isoparametric $\mathcal{P}_{4}$
triangle cells]{\begin{centering}
\includegraphics[width=0.9\linewidth]{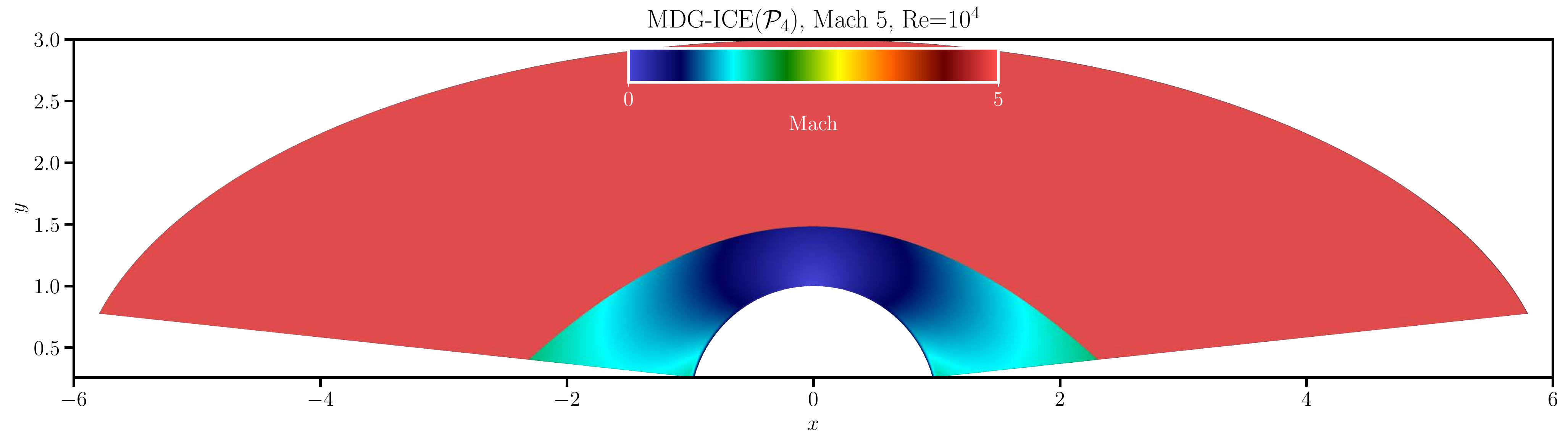}
\par\end{centering}
}
\par\end{centering}
\begin{centering}
\subfloat[\label{fig:Viscous-Bow-Shock-Pressure-Re-010000}527 isoparametric
$\mathcal{P}_{4}$ triangle cells]{\begin{centering}
\includegraphics[width=0.9\linewidth]{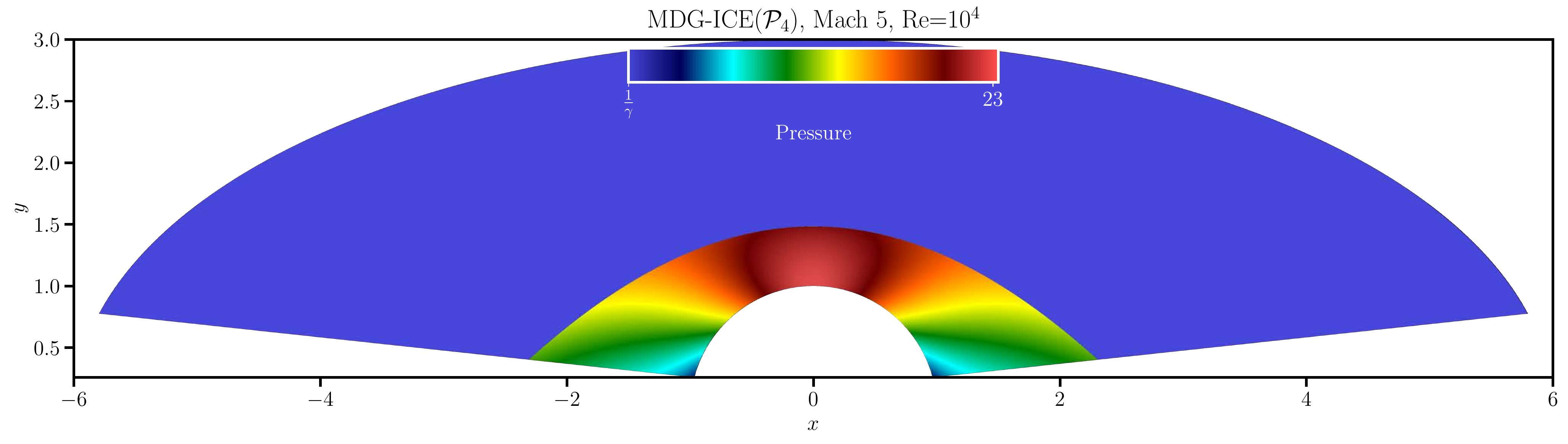}
\par\end{centering}
}
\par\end{centering}
\centering{}\caption{\label{fig:Viscous-Bow-Shock-MDG-ICE-Re-010000}The MDG-ICE solution
computed using $\mathcal{P}_{4}$ isoparametric triangle elements
for the viscous Mach 5 bow shock at $\mathrm{Re}=10^{4}$. The MDG-ICE
solution was initialized from the MDG-ICE solution at $\mathrm{Re}=10^{3}$
shown in Figure~\ref{fig:Viscous-Bow-Shock-DG-Mesh}. The location
of the shock along the line $x=0$ was computed as $y=1.48437$ for
a stand-off distance of $0.48437$.}
\end{figure}
\begin{figure}
\begin{centering}
\subfloat[\label{fig:Viscous-Bow-Shock-Mesh-Re-100000}768 isoparametric $\mathcal{P}_{4}$
triangle cells]{\begin{centering}
\includegraphics[width=0.9\linewidth]{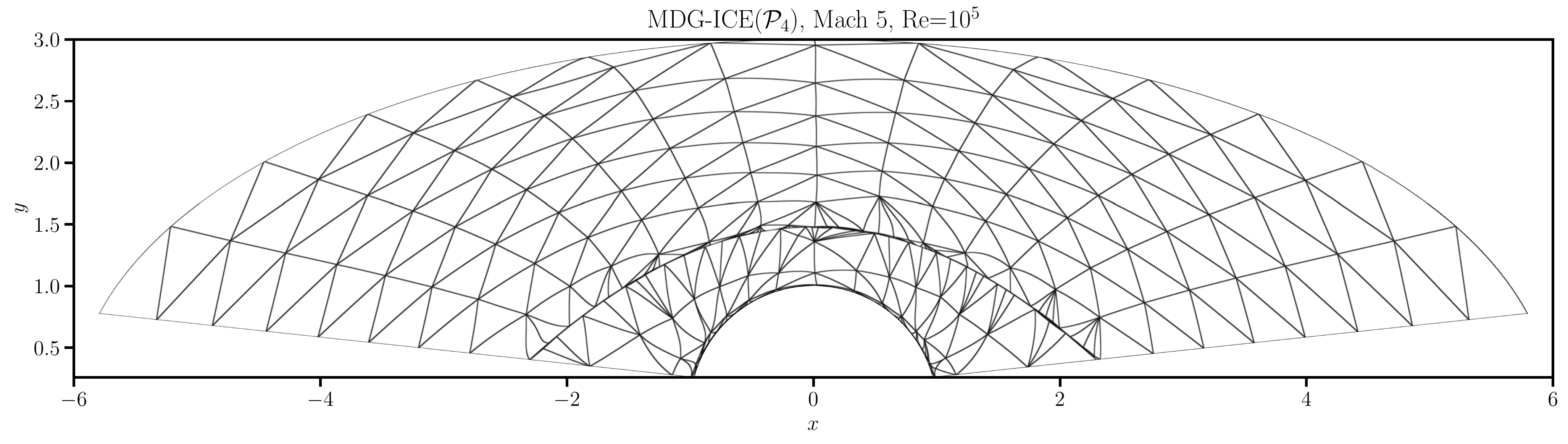}
\par\end{centering}
}
\par\end{centering}
\begin{centering}
\subfloat[\label{fig:Viscous-Bow-Shock-Temperature-Re-100000}768 isoparametric
$\mathcal{P}_{4}$ triangle cells]{\begin{centering}
\includegraphics[width=0.9\linewidth]{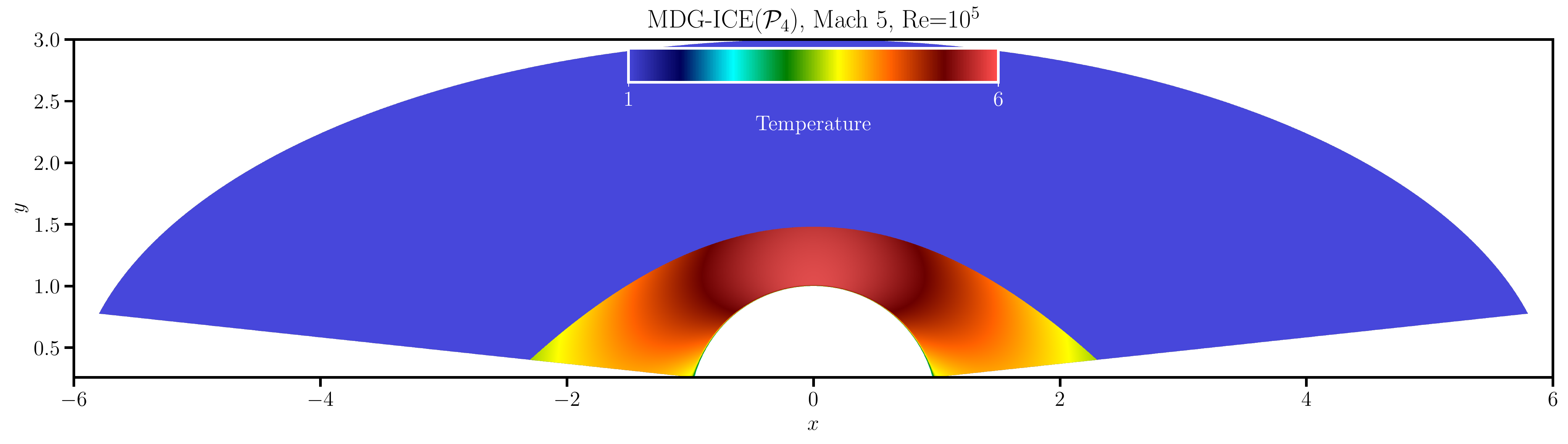}
\par\end{centering}
}
\par\end{centering}
\begin{centering}
\subfloat[\label{fig:Viscous-Bow-Shock-Mach-Re-100000}768 isoparametric $\mathcal{P}_{4}$
triangle cells]{\begin{centering}
\includegraphics[width=0.9\linewidth]{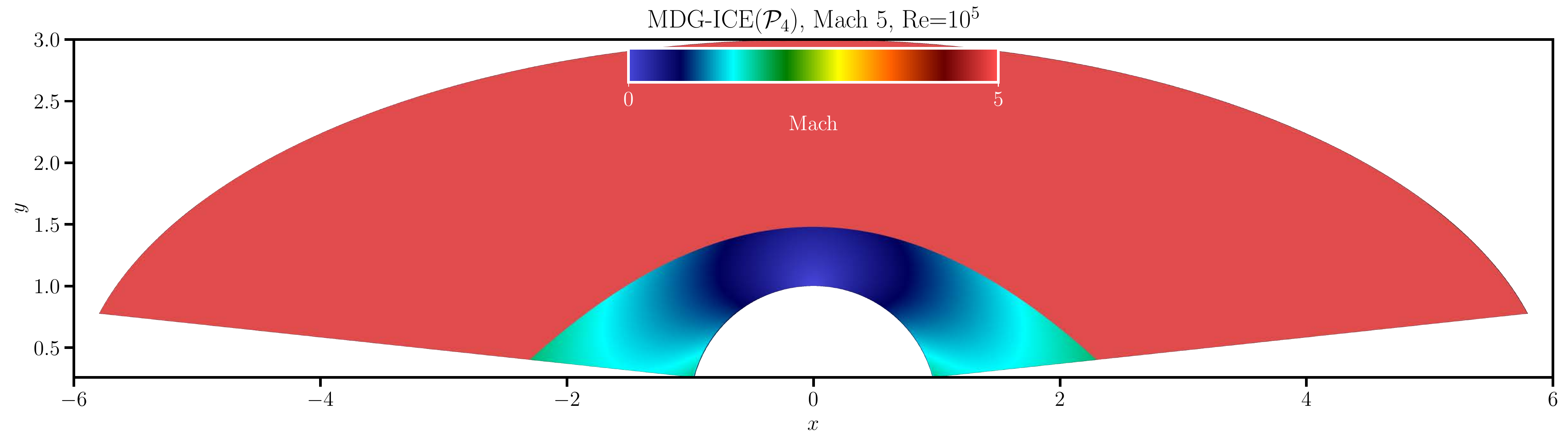}
\par\end{centering}
}
\par\end{centering}
\centering{}\subfloat[\label{fig:Viscous-Bow-Shock-Pressure-Re-100000}768 isoparametric
$\mathcal{P}_{4}$ triangle cells]{\begin{centering}
\includegraphics[width=0.9\linewidth]{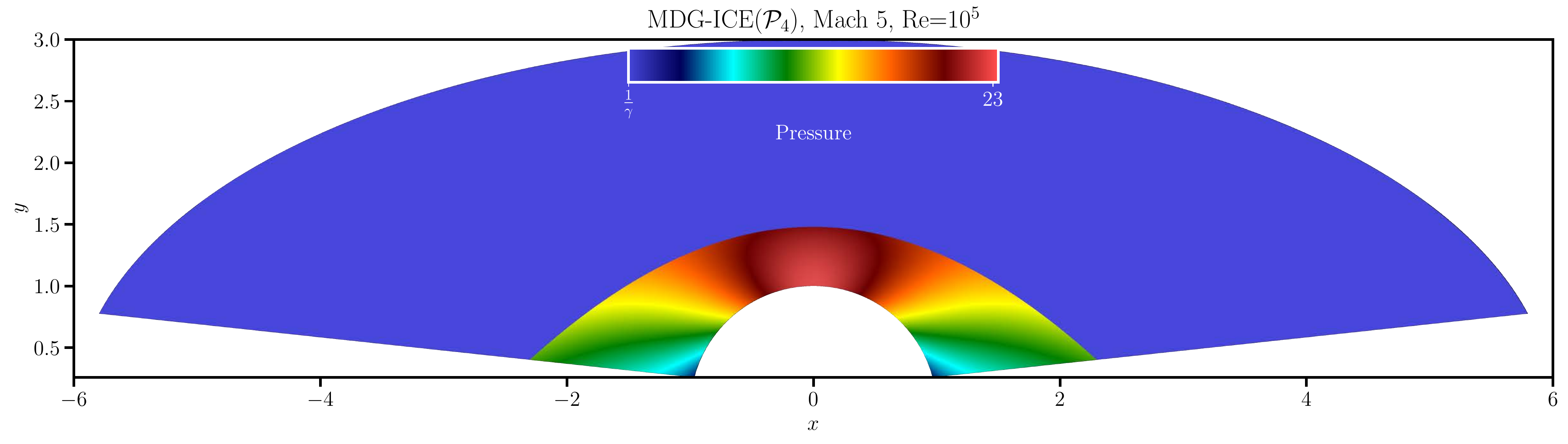}
\par\end{centering}
}\caption{\label{fig:Viscous-Bow-Shock-MDG-ICE-Re-100000}The MDG-ICE solution
computed using $\mathcal{P}_{4}$ isoparametric triangle elements
for the viscous Mach 5 bow shock at $\mathrm{Re}=10^{5}$. The MDG-ICE
solution was initialized from the MDG-ICE solution at $\mathrm{Re}=10^{4}$
shown in Figure~\ref{fig:Viscous-Bow-Shock-DG-Mesh}. The location
of the shock along the line $x=0$ was computed as $y=1.4809125$
for a stand-off distance of $0.4809125$.}
\end{figure}
\begin{figure}
\centering{}%
\begin{minipage}[c][1\totalheight][t]{0.3\columnwidth}%
\subfloat[\label{fig:Viscous-Bow-Shock-Temperature-Re-00001000-zoom}The final
grid and temperature fields of the MDG-ICE$\left(\mathcal{P}_{4}\right)$
solution computed using 400 $\mathcal{P}_{4}$ isoparametric triangle
elements for the viscous Mach 5 bow shock at $10^{3}$ $\mathrm{Re}$
shown in Figure~\ref{fig:Viscous-Bow-Shock-Mesh} and Figure~\ref{fig:Viscous-Bow-Shock-Temperature}
respectively.]{\begin{centering}
\begin{tabular}{c}
\includegraphics[width=0.8\linewidth]{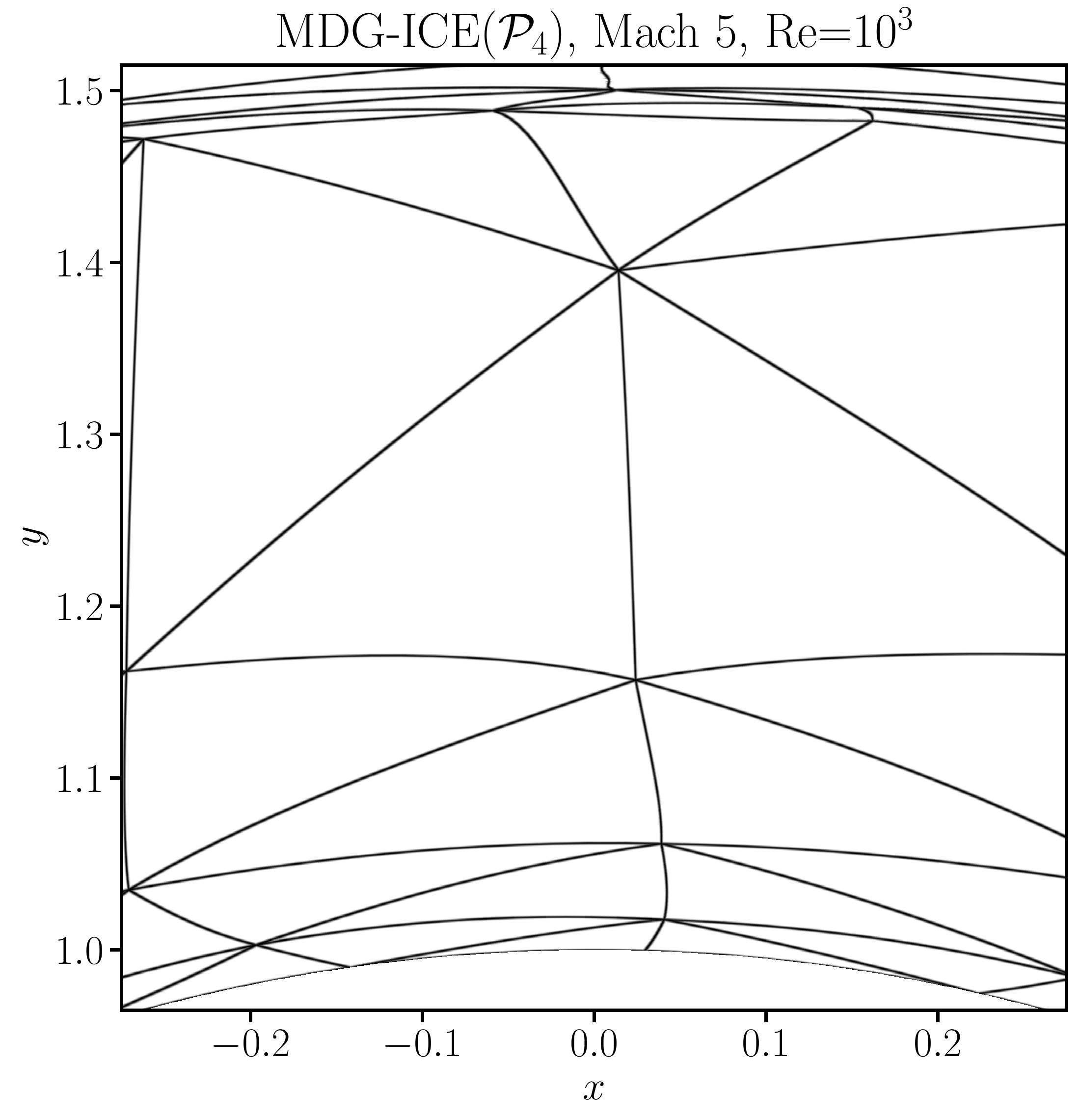}\tabularnewline
\includegraphics[width=0.8\linewidth]{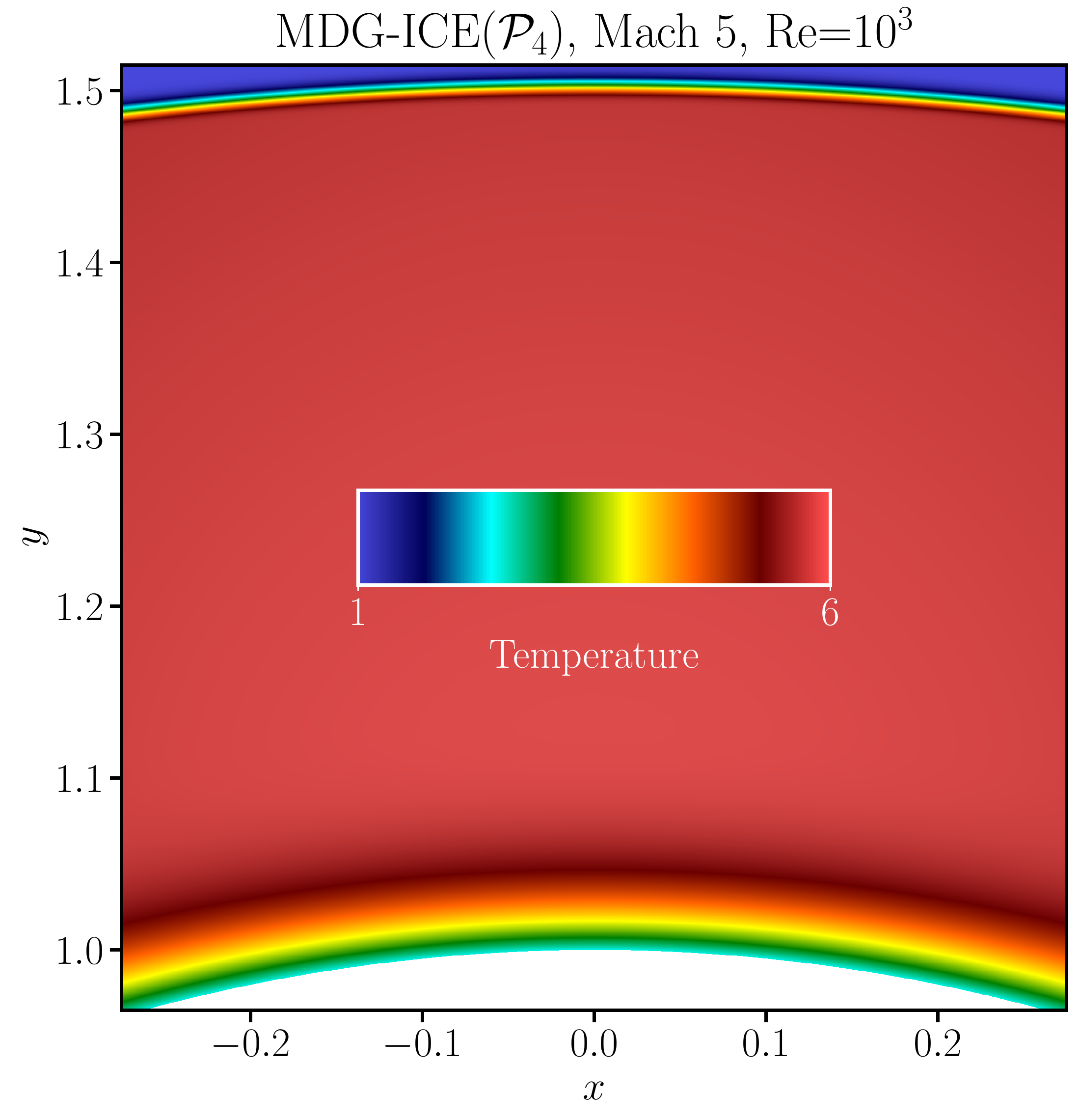}\tabularnewline
\end{tabular}
\par\end{centering}
}%
\end{minipage}%
\begin{minipage}[c][1\totalheight][t]{0.3\columnwidth}%
\subfloat[\label{fig:Viscous-Bow-Shock-Temperature-Re-00010000-zoom}The final
grid and temperature fields of the MDG-ICE$\left(\mathcal{P}_{4}\right)$
solution computed using 527 $\mathcal{P}_{4}$ isoparametric triangle
elements for the viscous Mach 5 bow shock at $10^{4}$ $\mathrm{Re}$
shown in Figure~\ref{fig:Viscous-Bow-Shock-Mesh-Re-010000} and Figure~\ref{fig:Viscous-Bow-Shock-Temperature-Re-010000}.]{\begin{centering}
\begin{tabular}{c}
\includegraphics[width=0.8\linewidth]{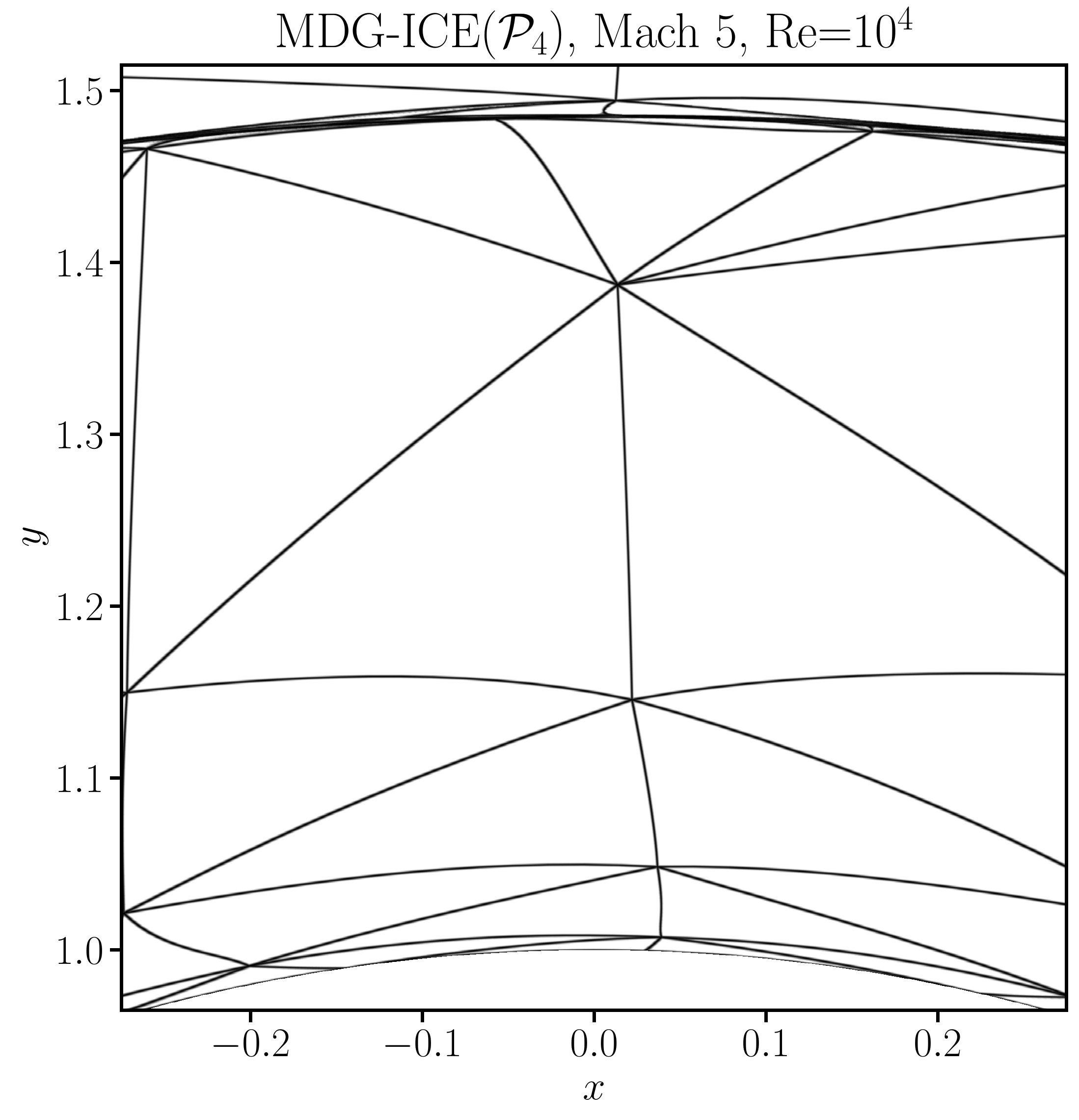}\tabularnewline
\includegraphics[width=0.8\linewidth]{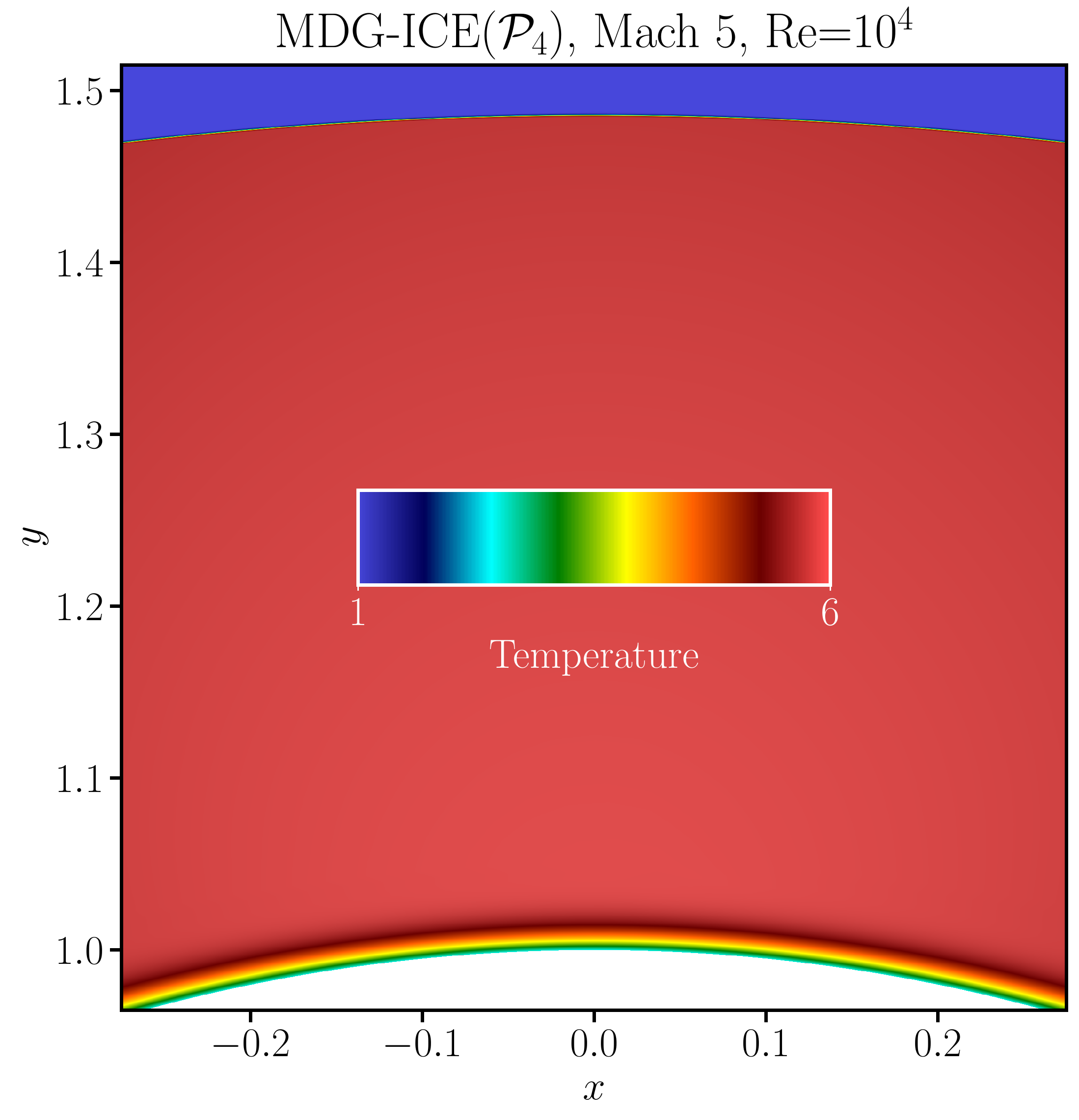}\tabularnewline
\end{tabular}
\par\end{centering}
}%
\end{minipage}%
\begin{minipage}[c][1\totalheight][t]{0.3\columnwidth}%
\subfloat[\label{fig:Viscous-Bow-Shock-Temperature-Re-00100000-zoom}The final
grid and temperature fields of the MDG-ICE$\left(\mathcal{P}_{4}\right)$
solution computed using 768 $\mathcal{P}_{4}$ isoparametric triangle
elements for the viscous Mach 5 bow shock at $10^{5}$ $\mathrm{Re}$
shown in Figure~\ref{fig:Viscous-Bow-Shock-Mesh-Re-100000} and Figure~\ref{fig:Viscous-Bow-Shock-Temperature-Re-100000}.]{\begin{centering}
\begin{tabular}{c}
\includegraphics[width=0.8\linewidth]{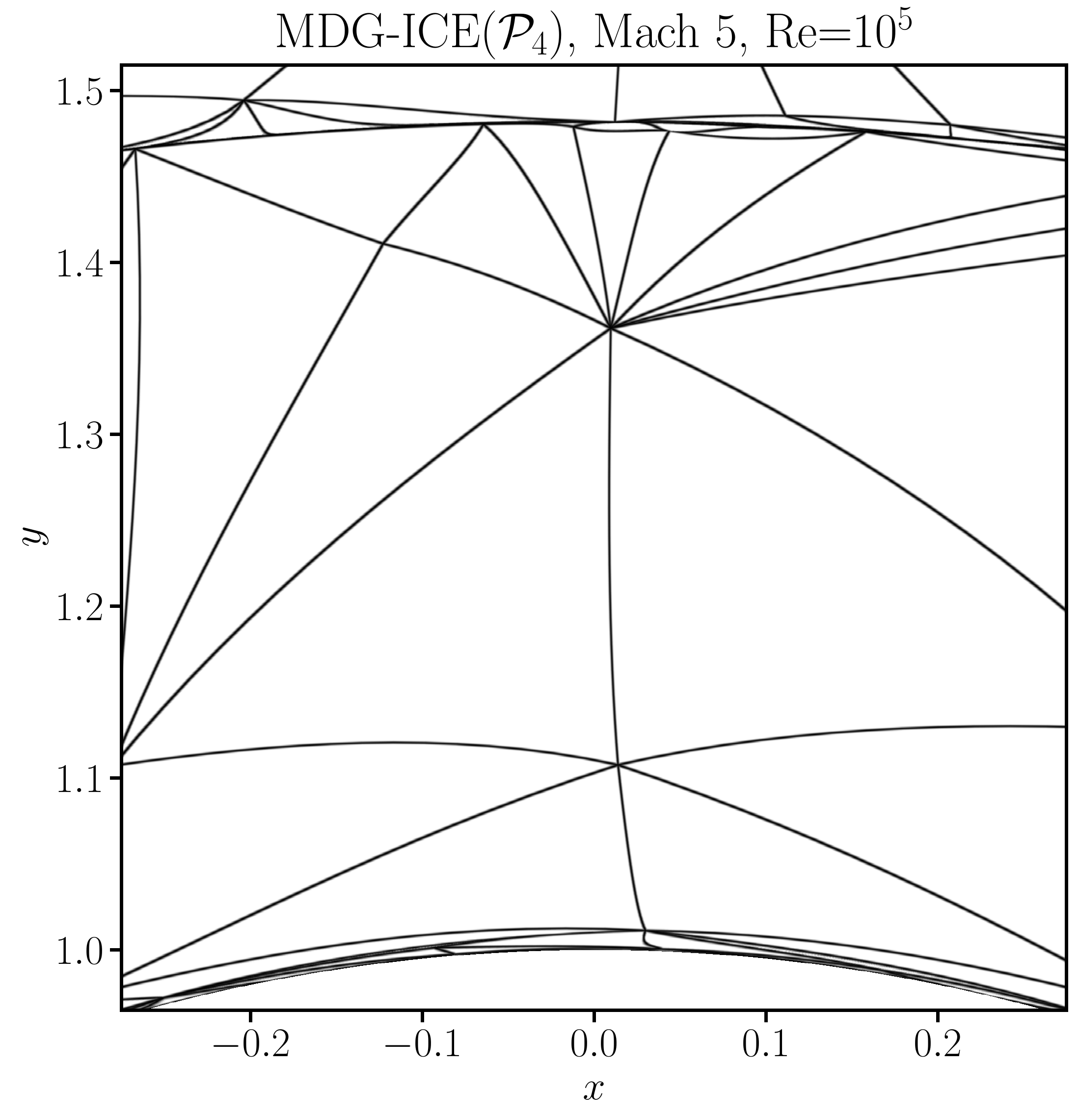}\tabularnewline
\includegraphics[width=0.8\linewidth]{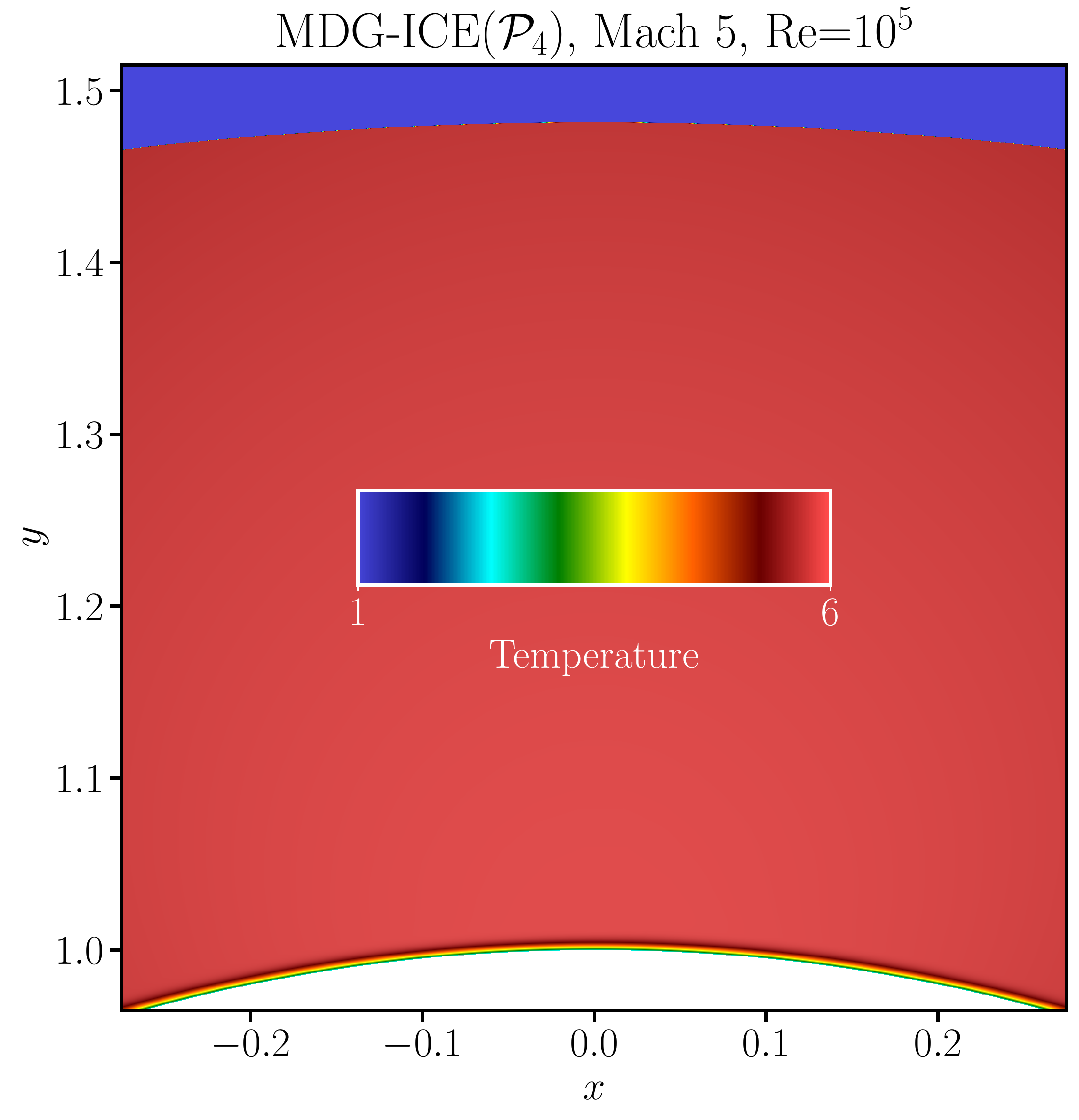}\tabularnewline
\end{tabular}
\par\end{centering}
}%
\end{minipage}\caption{\label{fig:Viscous-Bow-Shock-MDG-ICE-zoom}The final grid and temperature
fields corresponding to the MDG-ICE solution computed using $\mathcal{P}_{4}$
isoparametric triangle elements for the viscous Mach 5 bow shock at
$\mathrm{Re}=10^{3}$, $\mathrm{Re}=10^{4}$, and $\mathrm{Re}=10^{5}$.
Local edge refinement was used to adaptively split highly anisotropic
elements within the viscous structures as they were resolved by MDG-ICE.}
\end{figure}
The viscous MDG-ICE discretization is applied to a compressible Navier-Stokes
flow, described in Section~(\ref{subsec:Compressible-Navier-Stokes-flow}),
and used to approximate the solution to a supersonic viscous flow
over a cylinder in two dimensions. The solution is characterized by
the Reynolds number $\mathrm{Re}$, and the freestream Mach number
$\mathrm{M}_{\infty}$. The Reynolds number is defined as,

\begin{equation}
\mathrm{Re}=\frac{\rho vL}{\mu},\label{eq:reynolds-number}
\end{equation}
where $L$ is the characteristic length. The freestream Mach number
is defined $\mathrm{M}_{\infty}=\frac{v_{\infty}}{c_{\infty}}$ 
\begin{equation}
\mathrm{M}_{\infty}=\frac{v_{\infty}}{c_{\infty}},\label{eq:mach-number}
\end{equation}
where $v_{\infty}$ is the freestream velocity, $c_{\infty}=\sqrt{\nicefrac{\gamma P_{\infty}}{\rho_{\infty}}}$
is the freestream speed of sound, $P_{\infty}$ is the freestream
pressure, and $\rho_{\infty}$ is the freestream density. In this
work we consider Mach 5 flows at Reynolds numbers of $10^{3}$, $10^{4}$,
and $10^{5}$ traveling in the $\left(0,-1\right)$ direction, i.e.,
from top to bottom in Figure~\ref{fig:Viscous-Bow-Shock-DG} , Figure~\ref{fig:Viscous-Bow-Shock-MDG-ICE},
Figure~\ref{fig:Viscous-Bow-Shock-MDG-ICE-Re-010000}, and Figure~\ref{fig:Viscous-Bow-Shock-MDG-ICE-Re-100000}.
Supersonic inflow and outflow boundary conditions are applied at
the ellipse and outflow planes respectively. An isothermal no-slip
wall is specified at the surface of the cylinder of radius $r=1$
centered at the origin. The temperature at the isothermal wall is
given as $T_{\mathrm{wall}}=2.5T_{\infty}$, where $T_{\infty}$ is
the freestream temperature.

Figure~\ref{fig:Viscous-Bow-Shock-MDG-ICE} presents the MDG-ICE($\mathcal{P}_{4}$)
solution at $\mathrm{Re}=10^{3}$ computed on a grid of $400$ isoparametric
triangle elements. The MDG-ICE($\mathcal{P}_{4}$) solution was initialized
by interpolating the DG($\mathcal{P}_{2}$) field variables and the
corresponding grid. Figure~\ref{fig:Viscous-Bow-Shock-DG-Mesh} shows
the grid consisting of 392 linear triangle elements that was used
to initialize the MDG-ICE($\mathcal{P}_{4}$) by interpolating the
field variables. The high-order isoparametric boundary faces were
initialized by projecting to closest point on the boundary of the
domain via the boundary operator~(\ref{eq:geometric-boundary-condition}).
The MDG-ICE($\mathcal{P}_{4}$) field variables were initialized by
cell averaging the interpolated DG($\mathcal{P}_{2}$) solution shown
in Figure~\ref{fig:Viscous-Bow-Shock-DG-Temperature}. The MDG-ICE
auxiliary variable, with spatial components given by~(\ref{eq:navier-stokes-auxiliary-component}),
were initialized to zero for consistency with the initial piecewise
constant field variables. As the MDG-ICE solution converged, the previously
uniformly distributed points were relocated in order to resolve the
viscous shock. This resulted in a loss of resolution downstream, which
was conveniently handled via local edge refinement where highly anisotropic
elements were adaptively split as the viscous structures were resolved.
Although this was unnecessary for maintaining a valid solution, i.e.,
field variables that are finite and grid composed of only cells with
a positive determinant of the Jacobian, we found it sufficient for
maintaining a reasonable grid resolution, as compared to the initial
grid resolution downstream of the shock. The location of the shock
along the line $x=0$, estimated as the location corresponding to
the minimum normal heat flux, was computed as $y=1.49995$, giving
a stand-off distance of $0.49995$.

In one dimension, the viscous shock is described by a system of ordinary
differential equations that can be solved numerically, cf.~\cite{Cha10,Mas13}
for details. We use this solution to verify that the viscous MDG-ICE
formulation predicts the correct viscous shock profile when diffusive
effects are prominent, i.e., at low Reynolds number. Figure~\ref{fig:Bow-Shock-1D}
presents a comparison of an approximation of the exact solution for
a one-dimensional viscous shock to the centerline profiles of the
Mach 5 bow shock at $\mathrm{Re}=10^{3}$ for the following variables:
temperature, $T$, normal velocity, $v_{n}$, pressure, $p$, density,
$\rho$, normal component of the normal viscous stress tensor, $\tau_{nn}$,
and normal heat flux, $q_{n}$, where the normal is taken to be in
the streamwise direction. 

As expected, the one-dimensional profiles deviate from the two-dimensional
bow shock centerline profiles downstream of the viscous shock. The
one-dimensional solution assumes the viscous and diffusive fluxes
are zero outside of the shock. This is not the case for the two-dimensional
bow shock geometry where the blunt body and corresponding boundary
layer produce gradients in the solution downstream of the shock. For
density, in which case the diffusive flux is zero, the jump across
the viscous shock is directly comparable to the one-dimensional solution.
We also directly compare the exact solution to a one-dimensional viscous
shock profile computed by MDG-ICE($\mathcal{P}_{4}$) using 16 isoparametric
line cells. Figure~\ref{fig:Bow-Shock-1D-Density} shows that MDG-ICE
accurately reproduces the exact shock structure of density profile
with only a few high-order anisotropic curvilinear cells.

For reference, the exact and approximate values at the stagnation
point are marked on the centerline plots with the symbol $\times$.
An approximate value corresponding to the inviscid solution was used
when the exact value for the viscous flow was unavailable, e.g., the
stagnation pressure marked in Figure~\ref{fig:Bow-Shock-1D-Pressure}.
Although the analytic stagnation pressure for a inviscid flow neglects
viscous effects, it is not expected to differ significantly from the
value corresponding to the viscous solution for the problem considered
here, as shown in Table~1 in the work of Williams et al.~\cite{Wil16}.

We also report the pressure coefficient and the Stanton number, sampled
at the degrees of freedom, on the cylindrical, no-slip, isothermal
surface. The pressure coefficient at the surface is defined as

\begin{equation}
C_{p}=\frac{p-p_{\infty}}{\frac{1}{2}\rho_{\infty}v_{\infty}^{2}},\label{eq:pressure-coefficient}
\end{equation}
where $p_{\infty}$, $\rho_{\infty}$, and $v_{\infty}$ are the freestream
pressure, density, and velocity respectively. The Stanton number at
the surface is defined as

\begin{equation}
C_{h}=\frac{q_{\mathrm{n}}}{c_{p}\rho_{\infty}v_{\infty}\left(T_{t,\infty}-T_{\mathrm{wall}}\right)},\label{eq:Stanton-Number}
\end{equation}
where $q_{n}$ is the normal heat flux, $T_{\mathrm{wall}}$ is the
wall temperature, and $T_{t,\infty}$ is the freestream stagnation
temperature. In Figure~\ref{fig:Bow-Shock-1D-PressureCoefficient}
the pressure coefficient on the cylindrical surface is plotted and
the exact pressure coefficient at the stagnation point for an inviscid
flow, $C_{p}\approx1.8087699607027568$, is marked with the symbol
$\times$.

In order to compute solutions at higher Reynolds numbers, continuation
of the freestream viscosity, $\mu_{\infty}$ was employed. Figure~\ref{fig:Viscous-Bow-Shock-MDG-ICE-Re-010000}
presents the MDG-ICE($\mathcal{P}_{4}$) solution at $\mathrm{Re}=10^{4}$
computed on a grid of $527$ isoparametric triangle elements, which
was initialized from the $\mathrm{Re}=10^{3}$ MDG-ICE($\mathcal{P}_{4}$)
solution. Figure~\ref{fig:Viscous-Bow-Shock-MDG-ICE-Re-010000} presents
the MDG-ICE($\mathcal{P}_{4}$) solution at $\mathrm{Re}=10^{5}$
computed on a grid of $768$ isoparametric triangle elements, which
was initialized from the $\mathrm{Re}=10^{4}$ MDG-ICE($\mathcal{P}_{4}$)
solution. As in the $\mathrm{Re}=10^{3}$ case, local edge refinement
was used to adaptively split highly anisotropic elements within the
viscous structures as they were resolved by MDG-ICE. At these higher
Reynolds numbers, local refinement was necessary to maintain a valid
grid. In addition to the splitting of highly anisotropic cells, local
edge refinement was also applied to cells in which the determinant
of the Jacobian became non-positive.

Figure~\ref{fig:Viscous-Bow-Shock-MDG-ICE-zoom} compares the MDG-ICE($\mathcal{P}_{4}$)
solutions directly downstream of the viscous shock at $\mathrm{Re}=10^{3}$,
$\mathrm{Re}=10^{4}$, and $\mathrm{Re}=10^{5}$. By adapting the
grid to the flow field, MDG-ICE is able to simultaneously resolve
the thin viscous structure over a range of Reynolds numbers. As the
MDG-ICE solution converges, elements within regions that contain strong
gradients become highly anisotropic and warp nonlinearly to conform
to both the curved shock geometry and efficiently resolve the flow
around the curved blunt body. Thus, unlike a posteriori anisotropic
mesh adaptation, MDG-ICE achieves high-order anisotropic curvilinear
$r$-adaptivity as an intrinsic part of the solver. Furthermore, MDG-ICE
automatically repositions the nodes in order to resolve the flow field
over different length scales as the Reynolds number is increased from
$10^{3}$ to $10^{5}$. As such, MDG-ICE overcomes another challenge
associated with a posteriori anisotropic mesh adaptation which produces
regions of excessive refinement on the scale of the coarse mesh cell
size and therefore must rely on grid coarsening to limit the region
of refinement to the more appropriate length scale corresponding to
the feature under consideration.  

\section{Conclusions and future work\label{sec:Conclusions}}

The Moving Discontinuous Galerkin Method with Interface Condition
Enforcement (MDG-ICE) has been applied to viscous flow problems, involving
both linear and nonlinear viscous fluxes, where it was shown to detect
and resolve previously under-resolved flow features. In the case of
linear advection-diffusion, MDG-ICE adapted the grid to resolve the
initially under-resolved boundary layer, thereby achieving spectral
convergence and a more accurate solution than the best possible approximation
on a uniform static grid, which is given by the $L^{2}$ projection
of the exact solution. Unsteady flows were computed using a space-time
formulation where viscous structures were automatically resolved via
anisotropic space-time $r$-adaptivity. High speed compressible Navier-Stokes
solutions for a viscous Mach 5 bow shock at a Reynolds numbers of
$10^{3}$, $10^{4}$, and $10^{5}$ were presented. The viscous MDG-ICE
formulation was shown to produce the correct viscous shock profile
in one dimension for a Mach 5 flow at $\mathrm{Re}=10^{3}$. The one-dimensional
viscous shock profile was compared to the centerline profile of the
two-dimensional MDG-ICE solution where it was shown to accurately
compute both the shock profile and boundary layer profile simultaneously
using only a few high-order anisotropic curved cells within each region
and thus overcoming an ongoing limitation of anisotropic mesh adaptation.
Local edge refinement was used to adaptively split highly anisotropic
elements within the viscous structures as they were resolved by MDG-ICE.
Finally, MDG-ICE is a consistent discretization of the governing equations
that does not introduce low-order errors via artificial stabilization
or limiting and treats the discrete grid as a variable.

It should be noted that the internal structure of a viscous shock
may not be adequately described by the compressible Navier-Stokes
equations due to non-equilibrium effects~\cite{Pha89}, an issue
surveyed by Powers et al.~\cite{Pow15}. While in the present work
MDG-ICE was shown to provide highly accurate solutions to the compressible
Navier-Stokes equations, future work will apply MDG-ICE to an improved
multi-scale model that incorporates physical effects more adequately
described by the kinetic theory of gases~\cite{Kes10}. MDG-ICE is
a promising method to apply within such a framework due to its ability
to isolate regions in which enhanced physical modeling is required.

In future work, we will also develop a least-squares MDG-ICE formulation
with optimal test functions by applying the DPG methodology of Demkowicz
and Gopalakrishnan~\cite{Dem10,Dem11,Dem15_20}. Using this approach
we will demonstrate high-order convergence for both linear and nonlinear
problems. We also plan to mitigate the need for local $h$-refinement
by considering alternative methods for maintaining grid validity.
We will explore adaptively increasing the order of the local polynomial
approximation for cells within thin internal and boundary layers.
For instance, Chan et al.~\cite{Cha10,Cha14} used a combination
of $h$ and $p$ refinement to resolve viscous shocks. In their adaptive
strategy, $h-$refinement is used until the local grid resolution
is on the order of the viscous scale, at which point $p$-refinement
is used to further enhance accuracy. Additionally, scaling the regularization
by the inverse of volume of the cell, an approach used by Zahr et
al.~\cite{Zah19,Zah20_SCITECH} may also be effective for maintaining
grid validity at higher Reynolds numbers. Ultimately, we plan on maintaining
grid validity by incorporating smoothing, or untangling, into the
projection operator, Equation~(\ref{eq:geometric-boundary-condition}),
which enforces the geometric boundary conditions.

\section*{Acknowledgements}

This work is sponsored by the Office of Naval Research through the
Naval Research Laboratory 6.1 Computational Physics Task Area.

\bibliography{citations}

\end{document}